\documentclass[12pt]{amsart}
\usepackage[utf8]{inputenc}
\usepackage{pifont}

\usepackage{amsmath,amsthm,amssymb,mathscinet,bbm}
\usepackage{parskip}
\usepackage{geometry}
\usepackage{xcolor} 
\usepackage{xfrac, nicefrac}
\usepackage{exscale, relsize}

\usepackage{mathtools} 
\usepackage[hidelinks]{hyperref}

\author{Akash Singha Roy}
\address{Department of Mathematics \\ University of Georgia \\ Athens, GA 30602}
\email{akash01s.roy@gmail.com}


\renewcommand\phi\varphi
\usepackage{accents}

\renewcommand{\pod}[1]{\allowbreak\mathchoice
  {\if@display \mkern 18mu\else \mkern 8mu\fi (#1)}
  {\if@display \mkern 18mu\else \mkern 8mu\fi (#1)}
  {\mkern4mu(#1)}
  {\mkern4mu(#1)}
}
\usepackage{graphicx}
\DeclareMathAlphabet{\curly}{U}{rsfs}{m}{n}

\newcommand{\1}{\mathbbm{1}}

\newcommand\Z{\mathbb{Z}}
\allowdisplaybreaks

\usepackage{cleveref}
\crefname{section}{§}{§§}
\Crefname{section}{§}{§§}

\newcommand\NatNos{\mathbb N}

\newtheorem{thm}{Theorem}[section]
\newtheorem{cor}[thm]{Corollary}
\newtheorem{prop}[thm]{Proposition}
\newtheorem{lem}[thm]{Lemma}

\theoremstyle{remark}
\newtheorem*{rmk}{Remark}

\newcommand\err{\mathcal E}
\newcommand\reals{\mathbb R}

\newcommand\Ree{\mathrm{Re}}

\newcommand\bbm{\mathbbm 1}

\newcommand\chibara{\overline\chi(a)}

\newtheorem{innercustomgeneric}{\customgenericname}
\providecommand{\customgenericname}{}
\newcommand{\newcustomtheorem}[2]{
  \newenvironment{#1}[1]
  {
   \renewcommand\customgenericname{#2}
   \renewcommand\theinnercustomgeneric{##1}
   \innercustomgeneric
  }
  {\endinnercustomgeneric}
}

\newcustomtheorem{customthm}{Theorem}
\newcustomtheorem{customlemma}{Lemma}

\newcommand\sm\setminus

\newcommand\C{\mathbb C}
\newcommand\Res{\text{Res}}

\newcommand\chibar{\overline\chi}

\newcommand\Imm{\mathrm{Im}}

\newcommand\ds{\mathrm{d}s}

\newcommand\GammajBar{\overline{\Gamma_j}}

\newcommand\Dcz{\mathcal D(c_0)}
\newcommand\DczBigReduced{\Dcz \sm (-\infty, 1/\nu]}
\newcommand\DczSmallReduced{\Dcz \sm (-\infty, \eta_e/\nu]}
\newcommand\Elq{\mathcal L_q}
\newcommand\Elqt{\mathcal L_q(t)}
\newcommand\Log{\mathrm{Log}}
\newcommand\Fsnu{\mathcal F(s\nu)}
\newcommand\Fprimesnu{\mathcal F'(s\nu)}
\newcommand\chiNeChiZeroChie{\chi \ne \chi_0, \, \chi_e}
\newcommand\TStschiZ{\mathcal T^*(s, \chi_0)}
\newcommand\TschiZ{\mathcal T(s, \chi_0)}
\newcommand\TStschie{\mathcal T^*(s, \chi_e)}
\newcommand\Tschie{\mathcal T(s, \chi_e)}
\newcommand\Tschi{\mathcal T(s, \chi)}
\newcommand\Lprimesnuchi{L'(s\nu, \chi)}
\newcommand\Lsnuchi{L(s\nu, \chi)}

\newcommand\Hfunc{{\mathcal H}}

\newcommand\Hfuncs{\Hfunc(s)}
\newcommand\FsnuLogDeriv{\frac{\Fprimesnu}{\Fsnu}}
\newcommand\dz{\mathrm{d}z}
\newcommand\alphachiZ{\alpha_{\chi_0}}
\newcommand\alphachie{\alpha_{\chi_e}}
\newcommand\lambdaqElqt{\lambda_q\,\Elqt}
\newcommand\LambdaqElqt{\Lambda_q\,\Elqt}
\newcommand\TwoLambdaqLogqPower{(2\lambda_q \log q)^{\lambda_q + |\Ree(\alphachiZ)| + |\Ree(\alphachie)|}}

\newcommand\CalM{\mathcal M}
\newcommand\sigmanu{\sigma_\nu}
\newcommand\lambdaqlogq{\lambda_q\log q}
\newcommand\Lambdaqlogq{\Lambda_q\log q}
\newcommand\mcA{\mathcal A}
\newcommand\li{\mathrm{li}}
\newcommand\lambdaSt{\lambda^*}

\newcommand\Dqx{\mathcal D_q(x)}
\newcommand\DQx{\mathcal D_Q(x)}
\newcommand\Eqx{\mathcal E_q(x)}

\newcommand\omegaan{\omega_a(n)}
\newcommand\Omegaan{\Omega_a(n)}
\newcommand\pqa{p_{q, a}}
\newcommand\zPoweromegaan{z^{\omega_a(n)}}
\newcommand\zPowerOmegaan{z^{\Omega_a(n)}}

\newcommand\Gzs{\mathcal G_z(s)}

\newcommand\GzTils{\widetilde{\mathcal G}_z(s)}
\newcommand\Cjtil{\widetilde C_j}
\newcommand\CZeroTil{\widetilde C_0}
\newcommand\ceillogqlogell{\lceil \log q/\log \ell \rceil}

\newcommand\ellceilpower{\ell^{\ceillogqlogell}}

\numberwithin{equation}{section}
\begin{document} 
\title[LSD method for $L$-functions]{The Landau-Selberg-Delange method for products of Dirichlet $L$-functions, and applications, I} 
\begin{abstract}
The Landau--Selberg--Delange method gives precise asymptotic formulas for the partial sums $\sum_{n \le x} \, a_n$ of a Dirichlet series $\sum_n \, a_n/n^s$ that behaves like a complex power of the Riemann zeta function. However, situations often arise when the Dirichlet series behaves like a product of complex powers of several Dirichlet $L$--functions to a modulus $q$. In such situations, one often requires sharp asymptotic formulas for the partial sums $\sum_{n \le x} \, a_n$ that 
apply in much wider ranges of $q$ than permitted by known forms of the Landau--Selberg--Delange method. In this manuscript, we address this problem, giving new estimates on $\sum_{n \le x} \, a_n$ in ranges of $q$ that are (in most applications) much wider than attainable from previous results. Our results also weaken certain hypotheses on the size of  $\{a_n\}_n$. As applications of our main theorems, we extend Landau's classical results on the distribution of integers with prime factors restricted to progressions, and improve upon Chang, Martin and Nguyen's work on the distribution of the least invariant factors and least primary factors of multiplicative groups. We also extend the classical Sathe--Selberg theorem and study the local laws of the functions $\Omega_a(n)$ and $\omega_a(n)$, that count (with and without multiplicity, respectively), the number of prime divisors of $n$ lying in the progression $a$ mod $q$.      
\end{abstract}
\maketitle
\section{Introduction}
The subject of mean values of multiplicative functions is one of the most active areas of investigation in analytic number theory \cite{Ike31, Wie32, Ing42, Del54, wirsing61, wirsing67, halasz68, halasz71, MV01, GT17, GHS18, GHS19}. 
The Landau-Selberg-Delange ``LSD" method, which saw its beginnings in \cite{landau1908, selberg1954, Delange71LSD},  
forms one of the most powerful results in this subject. 
One of the most general formulations of this was obtained by Tenenbaum 
\cite[Theorem II.5.2]{tenenbaum15}, who 
gives for a sequence $\{a_n\}_n \subset \C$, an asymptotic formula of the shape 
\begin{equation}\label{eq:TenennumLSD}
\sum_{n \le x} \, a_n \, = \, \frac{x}{(\log x)^{1-\alpha}} \, \sum_{j=0}^N \, \frac{\kappa_j}{(\log x)^j} \, + \, O\left(\text{Error Terms}\right),
\end{equation}
under the hypothesis that the Dirichlet series $\sum_{n=1}^\infty \, a_n/n^s = \zeta(s)^\alpha \, G(s)$ for some reasonably well--behaved analytic function $G(s)$. 
 
He also assumed that the coefficients $\{a_n\}_n$ are termwise absolutely bounded by a sequence $\{b_n\}_n$  whose Dirichlet series  satisfies similar analytic  properties. 
Owing to its high precision and uniformity, Tenenbaum's formulation of the LSD method is widely applicable in a variety of problems in analytic number theory. {There has also been a lot of work extending this in various directions, by  Granville--Koukoulopoulos \cite{GK19}, Tenenbaum--de la Breteche \cite{TB21}, Chang--Martin \cite{CM20}, Cui--Wu \cite{CuiWu2014}, 
and by Phaovibul in his doctoral dissertation \cite{Phaovibul2015}.} 

However, situations often arise when the Dirichlet series $\sum_{n=1}^\infty a_n/n^s$ behaves like a 
function of the form $\prod_{\chi \bmod q} \, L(s \nu, \chi)^{\alpha_\chi}$, where the product is over all Dirichlet characters $\chi$ modulo $q$, and where $\nu>0$ is a fixed parameter. In such situations, several applications require estimates on the partial sums $\sum_{n \le x} \, a_n$, which hold true \textit{uniformly} as the modulus $q$ varies in a wide range,  
one of the typical ranges being the ``Siegel--Walfisz type" range $[1, (\log x)^K]$. One of the first main difficulties here is that there are potentially \textit{two essential singularities} in the zero free region that need to be worried about, one coming from the principal character $\chi_0$ and one from the exceptional character $\chi_e$. One might try to place this in the setting of Tenenbaum's result (or in the setting of any of its improvements alluded to above) 
by writing $\prod_{\chi \bmod q} \, L(s \nu, \chi)^{\alpha_\chi} \, = \, \zeta(s)^{\alpha_{\chi_0}} \cdot G(s)$, but doing this gives rise to at least three issues:
\begin{enumerate}
\item[(i)] The best available bounds on $G(s) = \prod_{\chi \ne \chi_0} \, L(s \nu, \chi)^{\alpha_\chi}$ 
often grow way too rapidly with $q$. 

\item[(ii)] The possibility of a Siegel zero modulo $q$ {further} complicates the analytic behavior of $G(s)$. (In particular, $G(s)$ itself may have an essential singularity in the zero free region.)  

\item[(iii)] The parameter $\nu$ (crucial in applications) needs additional care {beyond} scaling maneuvers. 
\end{enumerate}
Out of these, (i) and (ii) are the most major issues, and they cause the ``Error Terms" in  \eqref{eq:TenennumLSD} 
to blow up drastically as the modulus $q$ grows, thereby \textit{severely} impeding uniformity in $q$. The aforementioned extensions of his work are also unable to deal with any of these issues.   

Concrete examples of these phenomena occur in Theorems 3.4 and 3.6, as well as in Proposition 5.4 of 
\cite{CM20}, where the uniformity in the moduli $q$ is only up to small powers of $\log x$, even with a more explicit version of \cite[Theorem II.5.2]{tenenbaum15} that the authors obtain. 
(In fact, the uniformity in \cite[Theorems 3.4 and 3.6]{CM20} may even drop to a power of $\log \log x$ in some commonly--occurring situations.) 
More widespread examples arise in residue-class distribution problems of arithmetic functions: For instance, 
\cite{LPR21} studies the distribution of Euler's totient $\phi(n)$ and the sum--of--divisors $\sigma(n)$ in residue classes to varying \textit{prime} moduli $p$. Here an application of known forms of the LSD method 
can only give a range $p \le (\log \log x)^2$.  
For composite moduli, the ranges obtained 
are even worse. As such, all previous works \cite{PSR22, PSR22A, PSR23, PSR23A, SR23Sigma, SR23Add, SR23A, SR23B} studying the equidistribution of arithmetic  functions to growing moduli resort to a variety of 
``ad--hoc" methods. 

In this work, we find a new version of the Landau--Selberg--Delange method, which extends 
Tenenbaum's formulation in \cite{tenenbaum15} to situations when the Dirichlet series $\sum_n \, a_n/n^s$ behaves like the product $\prod_{\chi \bmod q} \, L(s \nu, \chi)^{\alpha_\chi}$. In several applications where this used to be unattainable via previous literature,  
our results give genuine asymptotic formulas in ranges of $q$ that are at least as wide as the ``Siegel--Walfisz" range $[1, (\log x)^K]$ (for any fixed $K>0$), 
and  
become even wider as better bounds on the Siegel--zero become available.   
We also give  
results (Theorem  \ref{thm:LFuncLSDVariant1} and Corollary \ref{cor:LFuncLSDVariant1LogPowerInterval} below) that assume some natural growth conditions on  $\{a_n\}_n$ on average, in place of the existence of the sequence $\{b_n\}_n$ that was required for \eqref{eq:TenennumLSD}. As applications of our main results, we can:
\begin{itemize}
\item Extend the classical works of Landau \cite{landau1908} on the distribution of integers with prime factors restricted to arithmetic progressions. 

\item 
Improve upon the works of Chang, Martin and Nguyen \cite{CM20, MN24}, on the distributions of the least invariant factors and least primary factors of unit groups. In these applications, it is also useful that our general results do not require $\{a_n\}_n$ to be a multiplicative function. 
\item Extend the works of Sathe \cite{sathe1953} and Selberg \cite{selberg1954}, and study the distribution of positive integers with a given number of prime factors from an arithmetic progression.  
\end{itemize} 
Before giving the statements of our main results, we summarize some of the ideas we use to deal  with issues (i)-(iii) mentioned above. 
First, we introduce \eqref{eq:lambdaq,mujDef} a parameter $\lambda_q$  which controls the correlation of the vector $(\alpha_\chi)_\chi$ with the vectors $(\chi(a))_\chi$ for all coprime residues $a$ mod $q$; in several applications, $\lambda_q$ is quite small.  
Second, inspired by work of Scourfield \cite{scourfield84} and other classical ideas, we study the logarithmic derivative of $\Fsnu \coloneqq \prod_\chi \, L(s\nu, \chi)^{\alpha_\chi}$ by using an ``inner contour shift" to give a series representation for $\mathcal F'(s\nu)/\Fsnu$ in terms of the zeros of $\{L(s, \chi)\}_\chi$. The additional feature here is that the terms of such a series are smoothed by a rapidly decaying function,  allowing us to bound $\mathcal F'(s\nu)/\Fsnu$ via classical zero density estimates. By constructing certain auxiliary functions, we then bound $|\Fsnu|$, which in turn allow us to control the residual integrals after  
performing a second (``outer") contour shift that is more reminiscent of 
known forms of the LSD method. 
Moreover, the contours  
we consider are modifications of some classical contours that account for the presence of a possible Siegel--zero modulo $q$ as well as for the parameter $\nu$. We deal with the other complications introduced by $\nu$ via certain   
averaging arguments.    
\subsection{The main results: General statements}\label{subsec:GeneralMainResults}~\\
We write complex numbers $s$ as $\sigma+it$, where $\sigma = \Ree(s)$ and $t = \Imm(s)$. Fix $\boldsymbol{c_0 \in (0, 1/3)}$ such that for any integer $q \ge 3$, the product $\prod_{\chi \bmod q} \, L(s, \chi)$ has no zero in the \textsf{classical zero--free region} $\{\sigma+it: \sigma>1-c_0/\log(q(|t|+1))\}$ except at most a simple real zero $\boldsymbol{\eta_e}$ \textsf{(the Siegel zero)} associated to a real character $\boldsymbol{\chi_e}$ \textsf{(the exceptional character)}. Fix $\boldsymbol{\nu>0}$ and  
let $\boldsymbol{\Omega: \reals_{\ge 0} \rightarrow \reals_{\ge 0}}$ be a \textbf{non--decreasing} function. Given complex numbers $\boldsymbol{\{a_n\}_{n=1}^\infty}$ and $\boldsymbol{\{\alpha_\chi\}_{\chi}}$ (with the $\alpha_\chi$ indexed at the Dirichlet characters $\chi$ mod $q$), we say that $\{a_n\}_n$ has \textsf{\textbf{property} $\boldsymbol{\mathcal P(\nu, \{\alpha_\chi\}_\chi; c_0, \Omega)}$} if 
\begin{equation}\label{eq:PropertP}  
\sum_{n=1}^\infty \, \frac{a_n}{n^s} = \Fsnu G(s) \text{ for all   }s\text{ with }\Ree(s)>\frac1\nu, ~ \text{ where }~~ \boldsymbol{\Fsnu \coloneqq \prod_{\chi}\, L(s\nu, \chi)^{\alpha_\chi}},
\end{equation} 
and where $\boldsymbol{G(s)}$ is a function that analytically continues into the region $\{s=\sigma+it: \sigma>\nu^{-1}(1-c_0/\log(q(|t\nu|+1))\}$ and satisfies $\boldsymbol{|G(s)| \le \Omega(|t|)}$ in this region. (Note that this region is a dilate of the classical zero--free region by the factor of $\nu$, which will be helpful to work with.) 

In what follows, we define 
\begin{equation}\label{eq:lambdaq,mujDef}
\boldsymbol{\lambda_q \coloneqq 1 + \max\limits_{\substack{a \bmod q}}~\left|\sum_\chi \, \alpha_\chi \cdot \chi(a)\right|}~~\text{ and }~~\boldsymbol{\mu_j \coloneqq \frac1{j!} \cdot \frac{\mathrm d^j}{\mathrm{d}s^j}\Bigg\vert_{s=1/\nu}\frac{\Fsnu G(s)}s\left(s-\frac1\nu\right)^{\alphachiZ}},
\end{equation}
so that $\mu_j$ is the $j$-th coefficient of the power series of the function $s^{-1}\Fsnu G(s) (s-1/\nu)^{\alphachiZ}$ around $1/\nu$. (The discussion in subsection \cref{subsec:AnalyticContinuations} shows that this function does analytically continue into some neighborhood of $1/\nu$.) Moreover, writing $\alpha_\chi = 
\sum_{\psi \bmod q} \, \alpha_\psi \cdot \phi(q)^{-1} \sum_{a \bmod q}  \, \chibara \psi(a)$ via orthogonality and interchanging sums, we obtain the following important bound
\begin{equation}\label{eq:alphachilambdaq}
|\alpha_\chi| \le \lambda_q~\text{ and }|\beta_\chi| \le \lambda_q \text{ for all characters $\chi$ mod $q$.}  
\end{equation}
Our first main result gives a Landau--Selberg--Delange type estimate on $\sum_{n \le x} a_n$, assuming an average growth condition on $\{a_n\}_n$ on dyadic intervals $(x, 2x]$. 
\begin{thm}\label{thm:LFuncLSDVariant1}
Assume that $\{a_n\}_n$ has property $\mathcal P(\nu, \{\alpha_\chi\}_\chi; c_0, \Omega)$, and that for all $x>1$, we have  
\begin{equation}\label{eq:Variant1GrowthConditionNew}
\sum_{x< n \le 2x}~|a_n| ~\le~ \kappa \, {x^{1/\nu}}, 
\quad \text{ where }\kappa \ge 2\text{ is independent of }x.
\end{equation}
Fix any $K_0>0$. The following  
hold uniformly in all $x \ge q \ge e^{4+5/3\nu}$, 
$h \in (0, x/2]$, $N \in \Z_{\ge 0}$, 
and in $\{\alpha_\chi\}_\chi \subset \C$ with $\max\{|\alphachiZ|, |\alphachie|\} \le K_0$;  
the implied constants depend only on $c_0, \nu$ and $K_0$. 

$\mathrm{\boldsymbol{(1)}}$ If the exceptional zero $\eta_e$ exists and satisfies $1-c_0/10\lambdaqlogq < \eta_e < 1-3\nu/\log x$,  
then  
\begin{multline}\label{eq:LFuncLSDVar1SiegelZero}\allowdisplaybreaks
\sum_{n \le x} a_n - \frac{x^{1/\nu}}{(\log x)^{1-\alphachiZ}} \sum_{j=0}^N~ 
\frac{\mu_j(\log x)^{-j}}{\Gamma(\alphachiZ-j)}~\ll~ \sum_{x < n \le x+h} \, |a_n| \, + \, \kappa \cdot \frac{x^{1+1/\nu} \, \log x}{Th}\\ 
+ (2 \lambda_q \log q)^{\lambda_q+2K_0} \, x^{1/\nu} \left\{\Omega(T)\, \frac{(\log(e\nu T))^{1+\lambda_q}}T ~+~  
\Omega(1) \, \frac{N! \, \big(71 (1+\nu) \cdot (1-\eta_e)^{-1}\big)^{N+K_0+1}}{(\log x)^{1-|\Ree(\alphachiZ)|} \cdot \min\{x/h, (\log x)^{N+1}\}}\right\}
\end{multline}
where $\displaystyle{T \coloneqq {(q\nu)^{-1/2}} \exp\left(0.5\sqrt{\log^2(q\nu) + {c_0\log x}/{\nu\lambda_q}}\right)}$. 

$\mathrm{\boldsymbol{(2)}}$ If $\eta_e$ does not exist or satisfies $\eta_e \le 1- c_0/10\lambdaqlogq$, then for $q < x^{c_0/80\nu\lambda_q}$, we have   
\begin{multline}\label{eq:LFuncLSDVar1NoSiegelZero}\allowdisplaybreaks
\sum_{n \le x} a_n - \frac{x^{1/\nu}}{(\log x)^{1-\alphachiZ}} \sum_{j=0}^N  
\frac{\mu_j(\log x)^{-j}}{\Gamma(\alphachiZ-j)}~\ll ~ \sum_{x < n \le x+h} \, |a_n| \, + \, \kappa \cdot \frac{x^{1+1/\nu} \, \log x}{Th}\\ 
+ (2 \lambda_q \log q)^{\lambda_q} \, x^{1/\nu} \left\{\Omega(T) \, \frac{(\log(e\nu T))^{1+\lambda_q}}T ~+~  
\Omega(1) \, \frac{N! \, \big(2000 (1+\nu) c_0^{-1} \cdot \lambda_q \log q\big)^{N+K_0+1}}{(\log x)^{1-|\Ree(\alphachiZ)|} \cdot \min\{x/h, (\log x)^{N+1}\}}\right\}, 
\end{multline}
where {$T \coloneqq (q\nu)^{-1/2}\exp\left(0.5\sqrt{\log^2(q\nu) + {c_0\log x}/\boldsymbol{20}\nu\lambda_q}\right)$. 
Moreover if $\eta_e$ does not exist, 
then \eqref{eq:LFuncLSDVar1NoSiegelZero} holds uniformly for $q < x^{c_0/\boldsymbol{8}\nu\lambda_q}$, we take $T = (q\nu)^{-1/2}\exp\left(0.5\sqrt{\log^2(q\nu) + {c_0\log x}/{\boldsymbol{4}\nu\lambda_q}}\right)$, we do not need to assume that $|\alphachie| \le K$, and $2000\nu$ can be replaced by $200\nu$.}
\end{thm} 
Theorem \ref{thm:LFuncLSDVariant1Gen} gives more general versions of this result. Under the Landau--Siegel zeros conjecture, 
Theorem \ref{thm:LFuncLSDVariant1} holds uniformly for $q \le x^{(1-\epsilon)c_0/\nu \lambda_q}$, with  
$\epsilon \in (0, 1)$ being an absolute constant. With some technical work, we should be able to improve this to allow 
\textit{any} fixed $\epsilon \in (0, 1)$. If $\eta_e$ \textit{does} exist, then better bounds on it lead to better ranges of uniformity in $q$ (up to the limit $x^{(1-\epsilon)c_0/\nu\lambda_q}$). Also note that in several applications ({including all the ones we provide}), the parameter $|\lambda_q|$ is absolutely bounded, often just by $1$. In such cases, Siegel's Theorem $1-\eta_e \gg_\epsilon~q^{-\epsilon}$ gives a range of $q \le (\log  x)^K$ (for any fixed $K>0$) in the above result, while a bound of the form $\eta_e \le 1-c(\theta_0)/(\log q)^{\theta_0}$ for some fixed $\theta_0 \in (0, 1]$ and some constant $c(\theta_0)>0$, such as in  
the {upcoming work of Yitang Zhang} \cite{zhang2022}, gives a range $q \le \exp\big((c(\theta_0) \log x/3\nu)^{1/\theta_0}\big)$. Further, note that the lower bound restriction $q \ge e^{4+5/3\nu}$ is an unserious technicality: For $q < e^{4+5/3\nu} \ll 1$, known forms of the Landau--Selberg--Delange method apply directly anyway. 

Our hypothesis $\mathcal P(\nu, \{\alpha_\chi\}_\chi; c_0, \Omega)$ generalizes hypothesis ``$\mathcal P(z; c_0, \delta, M)$" assumed in Tenenbaum \cite[Theorem II.5.2]{tenenbaum15}, not only by extending it to this setting of Dirichlet $L$-functions, but also by allowing a reasonably arbitrary growth condition in $G(s)$. Note also that we do not make any assumption on $\alpha_\chi$ for $\chi \ne \chi_0, \chi_e$. The additional generality in  
$\nu$ requires additional arguments (compared to the most natural case $\nu=1$) but has important applications. With these additional arguments, we can also directly generalize \cite[Theorem II.5.2]{tenenbaum15} from the case when the Dirichlet series behaves like a (complex) power of $\zeta(s)$ to when it behaves like a power of $\zeta(s\nu)$. 

When $\nu \in \NatNos$, \eqref{eq:Variant1GrowthConditionNew} follows if $\{a_n\}_n$ is supported on $\nu$--full numbers and bounded on average over intervals of the form $(x, 2x]$ 
(in particular, if $\nu=1$, then \eqref{eq:Variant1GrowthConditionNew} is just saying that $\{a_n\}_n$ is bounded on average on such intervals). In applications, \textit{much stronger} growth conditions than \eqref{eq:Variant1GrowthConditionNew} are often available: For instance, one of the most typical kinds intervals where strong growth conditions are easily available are intervals of length $x/(\log x)^A$. In such situations, we have the following useful consequence of Theorem \ref{thm:LFuncLSDVariant1}.  
\begin{cor}\label{cor:LFuncLSDVariant1LogPowerInterval}
Assume that $\{a_n\}_n$ has property $\mathcal P(\nu, \{\alpha_\chi\}_\chi; c_0, \Omega)$, and that for some {$A \ge 1$},  
\begin{equation}\label{eq:Variant1GrowthConditionSpecial}
\sum_{x< n \le x+x/(\log x)^A}~|a_n| ~\le~ \kappa_A \cdot \frac{x^{1/\nu}}{(\log x)^A}~ \text{ for all }x>1, \quad \text{ where }\kappa_A \ge 2\text{ is independent of }x.
\end{equation}    
Then the assertions of Theorem \ref{thm:LFuncLSDVariant1} hold exactly as stated, but with $\kappa \coloneqq \kappa_A \cdot 2^{A+1/\nu}$ and $h \coloneqq x/(\log x)^A$. In particular, $\min\{x/h, (\log x)^{N+1}\} \asymp (\log x)^{\min\{A, N+1\}}$,  and we have 
$$\sum_{x < n \le x+h} \, |a_n| \, \ll \, \kappa_A \cdot \frac{x^{1/\nu}}{(\log x)^A} \quad\text{ and }\quad \kappa \cdot \frac{x^{1+1/\nu} \, \log x}{Th} \, \ll \,2^A \, \kappa_A  \cdot \frac{x^{1/\nu}\log x}T.$$
\end{cor}
Although we just assumed  \eqref{eq:Variant1GrowthConditionSpecial} for \textit{some} $A$, in practice,  \eqref{eq:Variant1GrowthConditionSpecial} is often available for \textit{all} (sufficiently large) $A$. When this happens, 
$(\log x)^{1-|\Ree(\alphachiZ)|} \cdot \min\{x/h, (\log x)^{N+1}\} =  
(\log x)^{N+2-|\Ree(\alphachiZ)|}$, giving clean savings in the power of $\log x$ in \eqref{eq:LFuncLSDVar1SiegelZero} and \eqref{eq:LFuncLSDVar1NoSiegelZero}. (Here the uniformity of Corollary \ref{cor:LFuncLSDVariant1LogPowerInterval} in  
$A$ and $\kappa_A$, which comes from the uniformity of Theorem  \ref{thm:LFuncLSDVariant1} in the parameters $h$ and $\kappa$, is significant.)

In \cite[Theorem II.5.2]{tenenbaum15}, it is 
assumed that $\Omega(t) = \CalM (1+t)^{1-\delta_0}$ for some $\CalM>0$,  $\delta_0 \in (0, 1]$, and that there is a sequence 
$\{b_n\}_n$ upper--bounding the $|a_n|$ (termwise), whose Dirichlet series $\sum_n  \, b_n/n^s$ has analytic properties similar to those of $\sum_n \, a_n/n^s$. In Theorem \ref{thm:LFuncLSDVariant1}, we replaced these hypotheses by \eqref{eq:Variant1GrowthConditionNew} which is more natural in some applications. 
However, sometimes the following more direct generalization of \cite[Theorem II.5.2]{tenenbaum15} is easier to work with. 
\begin{thm}\label{thm:LFuncLSDVariant2}
Assume that $\{a_n\}_n \subset \C$ has property $\mathcal P(\nu, \{\alpha_\chi\}_\chi; c_0, \Omega)$, and that there exists a sequence $\{b_n\}_n \subset \reals_{\ge 0}$ satisfying $|a_n| \le b_n$ for all $n$,  
such that $\{b_n\}_n$ has property $\mathcal P(\nu, \{\beta_\chi\}_\chi; c_0, \Omega)$ for some $\{\beta_\chi\}_\chi \subset \C$. Further, suppose $\Omega(t) = \CalM (1+t)^{1-\delta_0}$ for some $\CalM>0$ and $\delta_0 \in (0, 1]$. 

Fix $K_0>0$ and define $\Lambda_q \coloneqq 1 + \max\limits_{\substack{a \bmod q}}~\max\left\{\left|\sum_\chi \,  \alpha_\chi \cdot \chi(a)\right|, \left|\sum_\chi \,\beta_\chi \cdot \chi(a)\right|\right\}$. The following estimates hold uniformly in all $x \ge 4$, $N \in \Z_{\ge 0}$, $q \ge e^{4+5/3\nu}$, and in all $\{\alpha_\chi\}_\chi, \{\beta_\chi\}_\chi \subset \C$ with $\max\{|\alphachiZ|, |\alphachie|, |\beta_{\chi_0}|, |\beta_{\chi_e}|\} \le K_0$; 
the implied constants depend only on $c_0, \nu, \delta_0$ and $K_0$. 

$\mathrm{\boldsymbol{(1)}}$ If the exceptional zero $\eta_e$ exists and satisfies $1-c_0/10\Lambdaqlogq < \eta_e < 1-3\nu/\log x$, then  
\begin{multline}\label{eq:LFuncLSDVar2SiegelZero}\allowdisplaybreaks
\sum_{n \le x} a_n - \frac{x^{1/\nu}}{(\log x)^{1-\alphachiZ}} \sum_{j=0}^N~ 
\frac{\mu_j(\log x)^{-j}}{\Gamma(\alphachiZ-j)} ~\ll~ \, \frac{(2 \Lambda_q \log q)^{\Lambda_q+|\Ree(\alphachiZ)|+|\Ree(\alphachie)|} \cdot N! \, (71 \nu)^N \cdot \CalM \, x^{1/\nu}}{(1-\eta_e)^{N+1+|\Ree(\alphachie)|} \cdot (\log x)^{N+2-|\Ree(\alphachiZ)|}}\\ +
\left(\frac{3 \Lambda_q \log q}{\delta_0^{1/2}}\right)^{\frac{3\Lambda_q}2+3K_0}  
\cdot \, \frac{\CalM \, (\log x)^{K_0/2}}{(1-\eta_e)^{5K_0/2}}\left\{\frac{x^{1/\nu}}{ T^{\delta_0/4}} + x^{(2+\eta_e)/3\nu} \, T^{\delta_0/4}\right\}, 
\end{multline}
where $\displaystyle{T \coloneqq {(q\nu)^{-1/2}} \, \exp\left(0.5\sqrt{\log^2(q\nu) + {2c_0\log x}/{\nu \delta_0 \Lambda_q}}\right)}$. 

$\mathrm{\boldsymbol{(2)}}$ If $\eta_e$ does not exist or satisfies $\eta_e \le 1- c_0/10\Lambdaqlogq$, then for $q < x^{c_0/80\nu\Lambda_q}$, we have   
\begin{multline}\label{eq:LFuncLSDVar2NoSiegelZero}\allowdisplaybreaks
\sum_{n \le x} a_n - \frac{x^{1/\nu}}{(\log x)^{1-\alphachiZ}} \sum_{j=0}^N~ 
\frac{\mu_j(\log x)^{-j}}{\Gamma(\alphachiZ-j)} ~\ll~ \, \frac{ (2\Lambda_q \log q)^{\Lambda_q + |\Ree(\alphachiZ)|} \cdot N! \, (2000 \nu c_0^{-1})^N  \cdot \CalM \, x^{1/\nu}}{(\Lambda_q \log q)^{-N-1} \cdot (\log x)^{N+2-|\Ree(\alphachiZ)|}}\\ + 
\left(\frac{3 \Lambda_q \log q}{\delta_0^{1/2}}\right)^{\frac{3\Lambda_q}2+3K_0} \cdot \, {\CalM \, x^{1/\nu} \, (\log x)^{K_0/2}} \, \exp\left(\frac{\delta_0}8\left\{\log(q\nu) -  \sqrt{\frac{c_0\log x}{10\nu\delta_0 \Lambda_q} + \log^2(q\nu)}\right\}\right).  
\end{multline}
Moreover if $\eta_e$ does not exist, then \eqref{eq:LFuncLSDVar2NoSiegelZero} holds uniformly for $q < x^{c_0/\boldsymbol{8}\nu\Lambda_q}$, we do not need to assume that $|\alphachie| \le K$, and $10\nu$ and $2000 \nu$ can be replaced by $\nu$ and $200\nu$ respectively.
\end{thm}
{The remarks made after the statement of Theorem \ref{thm:LFuncLSDVariant1} continue to hold.} We now highlight some of the  applications of Theorem  \ref{thm:LFuncLSDVariant1}, Corollary \ref{cor:LFuncLSDVariant1LogPowerInterval} and \ref{thm:LFuncLSDVariant2}, that we had alluded to previously. 
\subsection{The distribution of integers with prime factors restricted to progressions}\label{subsec:LandauRestrictedPrimeFactors}~\\
Let $U_q \coloneqq (\Z/q\Z)^\times$ denote the multiplicative group/unit group mod $q$. Consider an arbitrary subset $\mcA$ of $U_q$, and let $\mathcal N(x; q, \mcA)$ denote the number of positive integers $n \le x$ such that any prime dividing $n$ lies in one of the residue classes in $\mcA$. Estimating $\mathcal N(x; q, \mcA)$ is a classical problem back to the Landau's seminal work \cite{landau1908} on sums of squares, which also forms one of the earliest beginnings of the Landau--Selberg--Delange method. 

While Landau's work can deal with the case of \textit{fixed} $q$ and \textit{fixed} $\mcA$, it is of interest to ask for extensions of Landau's results when $q$ is allowed to grow with $x$ rapidly, and $\mcA$ is allowed to vary with $q$. Such a result has potential applications as well, as we shall see below. {However, there seems to be no literature attaining any  uniformity in $q$ or $\mcA$}, -- except for work of Chang and Martin \cite{CM20}, -- 
which allows uniformity in $q \le (\log x)^{1/2-\epsilon}$ and gives a huge error term when $\mcA$ is reasonably away from its extreme values $1$ or $\phi(q)$. 
As our first application of Theorem  \ref{thm:LFuncLSDVariant1},   
we give asymptotics on $\mathcal N(x; q, \mcA)$ that allow \textit{complete uniformity} in the modulus $q$ upto \textit{any} fixed power of $\log x$ (the full ``Siegel--Walfisz range") as well as complete uniformity in \textit{all} subsets $\mcA$ of $U_q$. Assuming the Landau--Siegel zeros conjecture, we can improve the range of $q$ up to $\exp(c' \sqrt{\log x})$ for some constant $c'$. ({Note that the range of $q$ in \cite{CM20} is not improvable conditionally either.})  
\begin{thm}\label{thm:App1PrimeFactorsinAPs} 
Fix $K>0$ and $\epsilon_0 \in (0, 1)$. There exists a constant $c_1$ depending only on $\epsilon_0$ and $K$ such that 
the following estimate holds uniformly in $x \ge 4$, $N\in \Z_{\ge 0}$, in moduli $q \le (\log x)^K$ and in $\mathcal A \subset U_q$; the implied constants are allowed to depend only on $c_0, c_1, K$ and $\epsilon_0$.  
\begin{multline}\allowdisplaybreaks\label{eq:App1PrimeFactorsinAPs}
\mathcal N(x; q, \mcA) \, - \, \frac{x}{(\log x)^{1-\frac{|\mcA|}{\phi(q)}}} \sum_{j=0}^N \, \frac{k_j \, (\log x)^{-j}}{\Gamma\left(\frac{|\mcA|}{\phi(q)}-j\right)} \, \ll \,\frac{N!\,(142 \, c_1^{-1})^N \cdot x}{(\log x)^{(N+2)(1-\epsilon_0) - \frac{|\mcA|}{\phi(q)}}} \, + \, x\exp\left(-\sqrt{\frac{c_0\log x}{32}}\right).
\end{multline}
Here the $\{k_j\}_j \subset \C$ depend only on $q$ and $\mcA$, and have been made explicit in \eqref{eq:App1PrimeFactorsinAPkj} and \eqref{eq:k0Formula}.  
Further, if the Siegel--zero does not exist, then all these assertions hold  uniformly in the wider range $q \le \exp(\sqrt{c_0\log x/8})$, and the constants $c_1$ and $\epsilon_0$  may be omitted. 
\end{thm}
The last assertion should also hold under a result of the form $\eta_e \gg 1-c(\theta_0)/(\log q)^{\theta_0}$. The constant $c_1$ above comes from Siegel's Theorem, and hence is ineffective.  Theorems \ref{thm:LFuncLSDVariant1} and \ref{thm:LFuncLSDVariant2} can also be used to  give more precise error terms, -- showing the dependence on the Siegel zero in the above application, -- but we have just chosen to highlight the cleanest corollaries. 
\subsection{Distributions of the least invariant factors of multiplicative groups}~\\
The algebraic structure of multiplicative group  $U_n \coloneqq (\Z/n\Z)^\times$ is of remarkable interest to number theorists. For example, its size is the Euler totient function $\phi(n)$, and $U_n$ is cyclic precisely when $n$ has a primitive root. It turns out that the invariant factor decomposition of $U_n$ also captures a lot of exotic  arithmetic information; here by the ``invariant factor decomposition", we mean the unique way of writing $U_n \cong \Z/d_1\Z \times \dots \times \Z/d_r\Z$, with  $d_1, \dots, d_r>1$ being integers  satisfying $d_1 \mid \dots \mid d_r$. (Recall that $d_1, \dots, d_r$ are called the \textsf{invariant factors} of $U_n$.) For instance, the length $r$ of this decomposition is basically the number $\omega(n)$ of distinct primes dividing $n$, while the \textit{largest} invariant factor $d_r$ is the renowned Carmichael lambda function $\lambda(n)$.   

Motivated by this last observation, Chang and Martin \cite{CM20} study the {dual} object, namely the  \textit{least} invariant factor $\lambdaSt(n) \coloneqq d_1$. Elementary linear algebra over rings shows that $\lambdaSt(n)$ must be even for $n>2$, but in fact they show that $\lambdaSt(n)$ must be equal to $2$ on a set of asymptotic density $1$. Some natural next questions are: How often is $\lambdaSt(n)$ equal to a given \textit{even} integer $q$? How often is $\lambdaSt(n)$ divisible by such a $q$? In \cite[Proposition 5.4 and Theorem 1.4]{CM20}, Chang and Martin also answer these questions for fixed $q$ or for $q$ varying roughly up to $(\log x)^{1/2-\epsilon}$.   ({Their range of $q$ is essentially the limit of their methods, even conditional on the Landau--Siegel zeros conjecture.}) 

As our next application, we qualitatively and quantitatively improve their results to allow $q$ to vary up to \textit{any} fixed power of $\log x$, with this range being further improvable under the Landau--Siegel zeros conjecture. We start by giving an asymptotic on $\#\{n \le x: q \mid \lambda^*(n)\}$, which will be one of two applications of our general main results in this paper, for a sequence $\{a_n\}_n$ which is \textit{not} defined by a multiplicative function. (The other application will be Theorem \ref{thm:MartinNguyen} below.)  
\begin{thm}\label{thm:Dqx}
Fix $K>0$ and $\epsilon_0 \in (0, 1)$. There exists a constant $c_1$ depending only on $\epsilon_0$ and $K$ such that the following estimate holds uniformly in $x \ge 4$, integers $N {\ge 0}$, and in \textbf{even} $q \le (\log x)^K$, with the implied constants  depending only on $c_0, c_1, K$ and $\epsilon_0$. 
\begin{multline*}\allowdisplaybreaks
\sum_{\substack{n \le x\\q \mid \lambdaSt(n)}} \, 1 \, = \, \frac{x}{(\log x)^{1-\frac{1}{\phi(q)}}} \sum_{j=0}^N \, \frac{r_j \, (\log x)^{-j}}{\Gamma\left(\frac{1}{\phi(q)}-j\right)} \, + \,O\left(\frac{N!\,(142 \, c_1^{-1})^N \cdot x}{(\log x)^{(N+2)(1-\epsilon_0) - \frac{1}{\phi(q)}}} \, + \, x\exp\left(-\sqrt{\frac{c_0\log x}{32}}\right)\right).
\end{multline*}
Here the $\{r_j\}_j \subset \C$ depend only on $q$ and $\mcA$, and are defined by 
\eqref{eq:App2ChangMartinrj}. Moreover, 
$$r_0 \, \asymp \, \left(\frac{\phi(q)}q \cdot \prod_{\chi \ne \chi_0} \, L(1, \chi)\right)^{1/\phi(q)} \, \ll \, (\log q) \cdot \left(\frac{\phi(q)}q\right)^{1/\phi(q)}\text{ uniformly in }q \le (\log x)^K.$$ 
If there is no Siegel--zero, then all these assertions hold uniformly for $q \le \exp(\sqrt{c_0\log x/8})$, and $c_1$ and $\epsilon_0$ do not appear. 
\end{thm}
This result is closely tied with Theorem \ref{thm:App1PrimeFactorsinAPs} via Lemma  \ref{lem:InvFactorDivisiblityCriteria}. Using Theorem \ref{thm:Dqx} and with some additional ideas from the anatomy of integers, we can also improve Chang and Martin's estimate on how often the least invariant factor of the unit group equals a given integer $q$,   
giving estimates that become sharper as $q$ grows even slightly rapidly with $x$. As such, both Theorems \ref{thm:Dqx} and \ref{thm:Eqx} also highlight a further application of Theorem \ref{thm:App1PrimeFactorsinAPs}. 
\begin{thm}\label{thm:Eqx}
Fix $K \ge 1$ and $\epsilon_0 \in (0, 1)$. There exists a constant $c_1$ depending only on $\epsilon_0$ and $K$ such that the following estimate holds uniformly in $x \ge 4$, integers $N {\ge 0}$, and in \textbf{even}  
$q \le (\log x)^K$, with the implied constants  depending only on $c_0, c_1, K$ and $\epsilon_0$. 
\begin{multline*}\allowdisplaybreaks
\#\{n \le x:~ \lambdaSt(n)=q\} \, = \, \frac{x}{(\log x)^{1-\frac{1}{\phi(q)}}} \sum_{j=0}^N \, \frac{r_j \, (\log x)^{-j}}  {\Gamma\left(\frac{1}{\phi(q)}-j\right)} \, - \frac{\li(x)+\li(x/2)}{\phi(q)}\\ \, + \,O\left(\frac{N!\,(142 \, c_1^{-1})^N \cdot x}{(\log x)^{(N+2)(1-\epsilon_0) - \frac{1}{\phi(q)}}} \, + \, 
\frac{x \, (\log \log x \cdot \log \log \log  x)^2}{q^2 (\log x)^{1-1/2\phi(q)}}\right).
\end{multline*}
Here the $r_j$ are as in Theorem \ref{thm:Dqx}, and $\li(y)= \int_2^y \, \mathrm dt/\log t$ is the logarithmic integral. 
\end{thm}
With some more work, it should be possible to remove the  $(\log \log \log x)^2$ factor. 
\subsection{Distributions of the least primary factors of multiplicative groups}\label{subsec:MartinNguyen}~\\
Another canonical form of finite abelian groups is their \textsf{primary decomposition}, which for the multiplicative group modulo $n$, amounts to writing $U_n \cong \Z/q_1\Z \times \dots \times \Z/q_k\Z$ for prime powers $q_1, \dots, q_k$. This decomposition is also unique (up to the order of the $q_j$), and we may define the \textsf{least primary factor of $U_n$} to be $\lambda'(n) \coloneqq \min_j \, q_j$. (For $n = 1, 2$, we define $\lambda'(n) \coloneqq \infty$.) Martin  and Nguyen \cite{MN24} studied the distribution of $\lambda'(n)$, obtaining asymptotics for how often $\lambda'(n) = q$ for $q=2$ and for a given \textbf{odd} prime power $q$. (Elementary number theory shows that the least primary factor of $U_n$ can never be $2^k$ for $k>1$.) 

Their key algebraic idea was to write $\#\{n \le x: \lambda'(n)=q\}$ as the difference $\mathcal A_q(x) - \mathcal A_{q^+}(x)$, where $\mathcal A_q(x) \coloneqq \#\{n \le x: \lambda'(n) \ge q\}$ and $q^+$ is the next prime power after $q$. Then they characterized $\mathcal A_q(x)$ with a set of congruence restrictions modulo primes less than $q$ (see Lemma \ref{lem:LeastPrimaryFactorCongruences}), which brought down the problem of estimating $\mathcal A_q(x)$ to the problem of estimating the number $\mathcal N(x; Q, \mathcal B)$ (defined in \cref{subsec:LandauRestrictedPrimeFactors}) with $Q$ being an integer having order of magnitude $e^{(1+o(1))q}$, and with $\mathcal B$ being a set of coprime residues modulo $Q$ of size really far away from the extreme values $1$ or $\phi(Q)$.  

Since this relies on Chang and Martin's estimate \cite[Theorem 3.4]{CM20} on $\mathcal N(x; Q, \mathcal B)$, the results in \cite{MN24} 
allow very little uniformity in $q$, something like $q = o(\log \log x)$ (because of the sizes of both $Q$ and $\mathcal B$), 
and this range is essentially the limit of their methods (even conditionally). But using the same algebraic idea as in Martin and Nguyen's work, Theorem \ref{thm:LFuncLSDVariant1} gives an unconditional range of $q \le K\log \log x$ for any fixed $K>0$, which can be improved to $q \le c\sqrt{\log x}$ conditional on the Landau--Siegel zeros conjecture. These limitations  
again come from $Q$ (but no longer from $\mathcal B$).  
\begin{thm}\label{thm:MartinNguyen}
Fix $K \ge 1$ and $\epsilon_0 \in (0, 1)$. There exists a constant $c_1$ depending only on $\epsilon_0$ and $K$ such that the following holds uniformly in $x \ge 4$, integers $N {\ge 0}$, and in \textbf{prime powers}  
$q \le K \log \log x$, with the implied constants  depending only on $c_0, c_1, K$ and $\epsilon_0$. 
\begin{multline*}\allowdisplaybreaks
\sum_{\substack{n \le x\\\lambda'(n)=q}} \, 1
 \, = \, \frac{x}{(\log x)^{1-B(q)}} \sum_{j=0}^N \, \frac{\kappa_j(q) \, (\log x)^{-j}}  {\Gamma\left(B(q)-j\right)} \, - \, \frac{x}{(\log x)^{1-B(q^+)}} \sum_{j=0}^N \, \frac{\kappa_j(q^+) \, (\log x)^{-j}}  {\Gamma\left(B(q^+)-j\right)} \\ \, + \,O\left(\frac{N!\,(142 \, c_1^{-1})^N \cdot x}{(\log x)^{(N+2)(1-\epsilon_0) - B_q}} \, + \, 
 x\exp\left(-\sqrt{\frac{c_0\log x}{32}}\right)\right).
\end{multline*}
Here, $q^+$ is the next prime power after $q$, the $\kappa_j(q)$ are as defined in \eqref{eq:kappajMartNg}, and 
\begin{equation}\label{eq:BqDef}
B(q) \, = \, \prod_{\ell < q} \, \left( \frac{\ell-2}{\ell-1} \, + \, \frac1{(\ell-1)\ell^{\lceil \log q/\log \ell\rceil - 1}} \right).
\end{equation}
If there is no Siegel--zero, then all these assertions hold uniformly for $q \le \sqrt{c_0\log x/128}$, and $c_1$ and $\epsilon_0$ do not appear.
\end{thm}
\subsection{The Sathe--Selberg Theorem in arithmetic progressions}~\\
The functions $\omega(n) = \#\{p \mid n\}$ and $\Omega(n) = \sum_{p \mid n} \, v_p(n)$  that count the number of primes dividing $n$ without and with multiplicity (respectively), have captured the interest of prominent number theorists for decades. One of the earliest renowned results on this subject, the {Hardy--Ramanujan} theorem, states that either of these functions concentrates around its mean.  
A conceptual improvement of this is the celebrated {Erd\H{o}s--Kac} theorem, which shows that either of these functions behaves like a {Gaussian random variable} with mean and variance both equal to $\log\log n$. 

The local laws of these functions (namely how often $\omega(n)$ or $\Omega(n)$ takes a given value $k$) are found in  works of L.~G.~Sathe \cite{sathe1953} and A.~Selberg \cite{selberg1954}, which form some of the cornerstones of probabilistic number theory. Their results roughly state that the functions $\omega(n)$ and $\Omega(n)$    
locally demonstrate {Poisson distribution} with parameter $\log \log x$, on the interval $[0, x]$.    

Natural generalizations of $\omega(n)$ and $\Omega(n)$ are the functions $\omega_a(n) \coloneqq \#\{p \mid n: p \equiv a \pmod q\}$ and $\Omega_a(n) \coloneqq \sum_{p \equiv a \pmod q} \, v_p(n)$, which count the number of primes that divide $n$ and lie in the progression $a$ mod $q$. Heuristic arguments suggest that both of these functions should also satisfy a ``Sathe--Selberg type" law, demonstrating local Poisson behavior with parameter $\log \log x/\phi(q)$. As our {final application} of Theorem \ref{thm:LFuncLSDVariant2} in this paper, we establish such a result in a precise quantitative form, -- allowing the modulus $q$ to vary within the Siegel--Walfisz range, -- and conditionally, -- in a still wider range. Our results below extend Theorems II.6.4 and II.6.5 in \cite{tenenbaum15}.  
\begin{thm}\label{thm:omegaa(n)}
Fix $K, \epsilon_0>0$.  
There exists a constant $c_1$ depending only on $\epsilon_0$ and $K$ such that uniformly in $x \ge 4$, $N {\in \Z_{\ge 0}}$, $k \in \NatNos$,  
in $q \le (\log x)^K$, and in coprime residues $a$ mod $q$,  
\begin{multline}\allowdisplaybreaks\label{eq:omegaa(n)Full}
\sum_{\substack{n \le x\\\omegaan = k}} \, 1 \, = \, {x} \sum_{j=0}^N \, \frac{P_{j, k}(\log \log x)}{(\log x)^{j+1/\phi(q)}} \, +  \,O\left(\frac{K^{-k} \cdot N!\,(71 \, c_1^{-1})^N \cdot x}{(\log x)^{(N+1)(1-\epsilon_0) - \frac{(K+1)}{\phi(q)}}} \, + \, \frac{x}{K^k}\exp\left(-\sqrt{\frac{c_0\log x}{16K}}\right)\right),
\end{multline}
where $P_{j, k}$ is a polynomial of degree at most $k$,  defined in \eqref{eq:omegaanPjk}. With $Y \coloneqq \log \log x$, and with $\pqa$ being the least prime in the progression $a$ mod $q$, we have uniformly in $q, k, a$ as above 
\begin{equation}\label{eq:P0k}
P_{0, k}(Y) - \left(1 - \frac1{\pqa} \right) \frac{(Y/\phi(q))^k}{k!} - \frac{(Y/\phi(q))^{k-1}}{\pqa \, (k-1)!} \, \ll \, 
\begin{cases}\vspace{4mm}
\displaystyle{\frac{\log q}{\phi(q)} \cdot \frac{(Y/\phi(q))^k}{k!}}, & \text{ if }\displaystyle{k \le \frac{KY}{\phi(q)}}\\
\displaystyle{\frac{\log q}{\phi(q)} \cdot \frac{e^{KY/\phi(q)}}{K^k}}, & \text{ for all }k.
\end{cases}
\end{equation} 
The implied constants in \eqref{eq:omegaa(n)Full} and \eqref{eq:P0k}  depend only on $c_0, c_1, K$, $\epsilon_0$. If there is no Siegel zero, then all these assertions hold uniformly for $q \le \exp(\sqrt{c_0\log x/20K})$,  
and $c_1, \epsilon_0$ don't appear. 
\end{thm} 
Theorem \ref{thm:omegaa(n)} gives a genuine asymptotic formula for $k \ll \log \log x$, as is also the range in the usual Sathe--Selberg Theorem for the function $\omega(n)$. 
For $k \le KY/\phi(q)$, we will deduce \eqref{eq:P0k} from a more precise estimate for $P_{0, k}(Y)$ that we will obtain using a variant of the saddle point method: See {Proposition \ref{prop:omegaaSPM}}. Also note that the error term in \eqref{eq:P0k} is $\ll \log q/K^k \phi(q)$ if $\phi(q) \gg KY$. 

A similar result also holds for $\Omega_a(n)$; only this time we need to make sure to stay away from the singularity of the Dirichlet series $\sum_{n=1}^\infty \, \zPowerOmegaan/n^s$ at $z = \pqa$. This is the exact same issue as the one encountered at the prime $2$ in the usual Sathe--Selberg theorem for the function $\Omega(n)$. 
\begin{thm}\label{thm:Omegaa(n)}
Fix $K, \epsilon_0>0$.  
There exists a constant $c_1$ depending only on $\epsilon_0$ and $K$ such that uniformly in $x \ge 4$, $N {\in \Z_{\ge 0}}$, $k \in \NatNos$,  
in $q \le (\log x)^K$, and in coprime residues $a$ mod $q$,  
\begin{multline}\allowdisplaybreaks\label{eq:Omegaa(n)Full}
\sum_{\substack{n \le x\\\Omegaan = k}} \, 1 \, = \, {x} \sum_{j=0}^N \, \frac{Q_{j, k}(\log \log x)}{(\log x)^{j+1/\phi(q)}} \, +  \,O\left(\frac{R^{-k} \cdot N!\,(71 \, c_1^{-1})^N \cdot x}{(\log x)^{(N+1)(1-\epsilon_0) - \frac{(K+1)}{\phi(q)}}} \, + \, \frac{x}{R^k}\exp\left(-\sqrt{\frac{c_0\log x}{16K}}\right)\right),
\end{multline}
where $R \coloneqq \min\{K, (1-\epsilon_0)\pqa\}$ and $Q_{j, k}$ is a polynomial of degree at most $k$,  defined in \eqref{eq:OmegaanQjk}. Moreover, the following estimates hold uniformly in $q, a$ as above, with $Y \coloneqq \log \log x$.\\
$\boldsymbol{\mathrm{(i)}}$ Uniformly in $k \le RY/\phi(q)$, we have
\begin{equation}\label{eq:Q0kSmallkSPM}
\displaystyle{Q_{0, k}(Y) \, = \, \frac{(Y/\phi(q))^k}{k!} \, \left\{ \left(1-\frac1{\pqa} \right) \left(1-\frac{k\phi(q)}{\pqa Y}\right)^{-1} \, + \, O\left(\frac{\log q}{\phi(q)} \, + \, \frac{k}{(\pqa Y/\phi(q))^2}\right) \right\}}.
\end{equation}
$\boldsymbol{\mathrm{(ii)}}$ Uniformly in $k \le (1-\epsilon_0) \pqa Y/\phi(q)$, we have 
\begin{equation}\label{eq:Q0kSmallk}
\displaystyle{Q_{0, k}(Y) \, =  \, \frac{(Y/\phi(q))^k}{k!} \left(1-\frac1{\pqa}\right) 
\, + \, O\left(\frac{(Y/\phi(q))^{k-1}}{(k-1)!} \, + \,  \frac{\log q}{\phi(q)} \cdot \frac{e^{RY/\phi(q)}}{R^k}\right)}. 
\end{equation} 
$\boldsymbol{\mathrm{(iii)}}$ Uniformly in $k \ge (1+\epsilon_0) \pqa Y/\phi(q)$, we have 
\begin{equation}\label{eq:Q0kLargek}
\displaystyle{Q_{0, k}(Y) \, =  \, \frac{e^{ {\pqa Y/\phi(q)} }}{(\pqa)^k} \left(1-\frac1{\pqa}\right) 
\, + \, O\left(\pqa \cdot \frac{(Y/\phi(q))^{k+1}}{(k+1)!} \, + \,  \frac{\log q}{\phi(q)} \cdot \frac{e^{RY/\phi(q)}}{R^k}\right)}. 
\end{equation}  
The implied constants in \eqref{eq:Omegaa(n)Full}, \eqref{eq:Q0kSmallkSPM} and \eqref{eq:Q0kLargek} depend only on $c_0, c_1, K$, $\epsilon_0$. If there is no Siegel zero, then all these assertions hold uniformly for $q \le \exp(\sqrt{c_0\log x/20K})$, 
and $c_1, \epsilon_0$ don't appear.
\end{thm}
In an upcoming sequel, we shall provide more useful variants of Theorems \ref{thm:LFuncLSDVariant1}--\ref{thm:LFuncLSDVariant2},  and highlight several more applications.  
\subsection*{Notation and Conventions} We do not consider zero function as  multiplicative (thus, $f(1)=1$ for any multiplicative function $f$). For $y>0$, any mention of ``$\log y$" will always mean the natural logarithm $\ln y$. All other logarithm conventions, branch cuts and  analytic continuations have been made explicit in subsection \cref{subsec:AnalyticContinuations}. Any mentions of  
the Siegel zero $\eta_e$ or the exceptional character 
$\chi_e$ are ignored if $\eta_e$ doesn't exist, i.e., if $\prod_\chi L(s, \chi)$ has no zeros inside the region $\{\sigma+it: \sigma>1-c_0/\log(q(|t|+1))\}$. Implied constants in $\ll$ and $O$-notation are allowed to depend on any parameters declared as ``fixed''; in particular, they are always allowed to depend on $K_0$, $c_0$ and $\nu$. Other dependence will be noted explicitly (for instance, with parentheses or subscripts). We use $\log_k$ to denote the $k$-th iterate of the natural logarithm.

As is commonplace, we write complex numbers $s$ as $\boldsymbol{\sigma+it}$ (with $\sigma = \Ree(s)$ and $t = \Imm(s)$). We denote a generic zero of Dirichlet $L$--functions by  $\boldsymbol{\rho = \beta+i\gamma}$,  where $\boldsymbol{\beta = \Ree(\rho)}$ and $\boldsymbol{\gamma = \Imm(\rho)}$; moreover,  $\boldsymbol{\sum_{\rho:\,L(\rho, \chi)=0}^*}~~$ denotes a sum over all zeros $\rho$ of $L(s, \chi)$ counted with appropriate multiplicity. {Other recurring notation has 
been defined in \eqref{eq:PropertP},  \eqref{eq:lambdaq,mujDef}, \eqref{eq:Lqt,DczDef} and  \eqref{eq:Hfunc}.}     
\section{Key analytic inputs: Logarithmic Derivatives, Auxiliary Functions, and the ``Inner Contour Shift"} \label{sec:AnalyticMachinery}
For any $\chi$ mod $q$, the function $\Log \Lsnuchi \coloneqq \sum_{p, r \ge 1}  \, \chi(p^r)/rp^{rs\nu}$ defines an analytic logarithm of $\Lsnuchi$ on the region $\{s: \sigma>1/\nu\}$. Hence, 
the function $\Fsnu$ is analytic on $\{s: \sigma>1/\nu\}$, and 
\begin{equation}\label{eq:FsnuDef}
\Fsnu = \prod_\chi \Lsnuchi^{\alpha_\chi} = \exp\left(\sum_\chi \alpha_\chi \, \Log \Lsnuchi\right)  = \exp\left(\sum_{p, r \ge 1} \frac1{rp^{rs\nu}} \sum_\chi \alpha_\chi {\chi(p^r)} \right)~~\text{ if $\sigma>1/\nu$}.
\end{equation}  
We analytically continue our functions into regions of interest. In what follows, we define 
\begin{equation}\label{eq:Lqt,DczDef}
\boldsymbol{\Elqt \coloneqq \log(q(|t\nu|+1))}~~\text{ and }~~\boldsymbol{\Dcz \coloneqq \left\{\sigma+it:~\sigma > \frac1\nu\left(1-\frac{c_0}{\Elqt}\right)\right\}}.
\end{equation} 
\subsection{Analytic Continuations}\label{subsec:AnalyticContinuations} Since the functions $L(s\nu, \chi_0)(s-1/\nu)$, $L(s\nu, \chi_e)(s-\eta_e/\nu)^{-1}$, and $\{\Lsnuchi\}_{\chiNeChiZeroChie}$ all continue analytically into nonvanishing functions on $\Dcz$, they have (unique) analytic logarithms $\TStschiZ$, $\TStschie$, and $\{\Tschi\}_{\chiNeChiZeroChie}$ on $\Dcz$ satisfying 
$$\mathcal T^*\left(\frac2\nu, \chi_0\right) = \sum_{p ,r \geq 1} \frac{\chi_0(p^r)}{rp^{2r}} + \ln\left(\frac2\nu - \frac1\nu\right),~~~\mathcal T^*\left(\frac2\nu, \chi_e\right) = \sum_{p ,r \geq 1} \frac{\chi_e(p^r)}{rp^{2r}} - \ln\left(\frac2\nu - \frac{\eta_e}\nu\right),$$
and $\mathcal T(2/\nu, \chi) = \sum_{p, r \ge 1} \chi(p^r)/rp^{2r}$ for all other $\chi$. (Thus $\TStschiZ$ is analytic on $\Dcz$ and satisfies $e^{\TStschiZ} = L(s\nu, \chi_0)(s-1/\nu)$ therein, etc.)  
Comparing derivatives, we see that the functions
\begin{equation}\label{eq:TschiZchieDef}
\TschiZ \coloneqq \TStschiZ - \log\left(s - \frac1\nu \right)~~\text{ and }~~\Tschie \coloneqq \TStschie + \log\left(s - \frac{\eta_e}\nu \right)
\end{equation}
define unique analytic continuations of the functions $\Log L(s\nu, \chi_0)$ and $\Log L(s\nu, \chi_e)$, 
into the regions $\DczBigReduced$ and $\DczSmallReduced$, respectively. (Here $\log z$ is the principal branch of the  logarithm, so $\log(s-1/\nu)$ is analytic on $\C \sm (-\infty, 1/\nu]$.) 
From this discussion, we see that the function $\exp(\sum_\chi ~ \alpha_\chi \Tschi) = \prod_\chi ~ e^{\alpha_\chi \Tschi}$ defines a unique analytic continuation of $\Fsnu$ in \eqref{eq:FsnuDef} into $\DczBigReduced$; hence, $\Fsnu = \exp(\sum_\chi \, \alpha_\chi \Tschi)$ for all $s$ in this region.

Note also that by the first equality in \eqref{eq:FsnuDef} and by analytic continuation, we may write  
\begin{equation}\label{eq:FsnuLogDerivAnaCont}
\frac{\Fprimesnu}{\Fsnu} = \sum_\chi \alpha_\chi \frac{\Lprimesnuchi}{\Lsnuchi}\text{ for all $s \ne 1/\nu$ s.t.~$s \ne \rho/\nu$ for any complex zero $\rho$ of $\prod_\chi L(s, \chi)$.}    
\end{equation} 
This relation is consistent with the analytic continuation of $\Fsnu$ 
in the previous paragraph. We will also need some {auxiliary functions}. By the above discussion (especially that around \eqref{eq:TschiZchieDef}), we see that if $\eta_e$ exists, then the function $s^{-1}\exp\left(\alphachiZ \TStschiZ + \alphachie \TStschie + \sum_{\chi \ne \chi_0, \chi_e}\,\alpha_\chi \Tschi\right)$ analytically continues the function 
$s^{-1} \Fsnu (s-1/\nu)^{\alphachiZ} (s-\eta_e/\nu)^{-\alphachie}$ into the region $\Dcz$. 

On the other hand, if $\eta_e$ doesn't exist (i.e.~all zeros of $\prod_\chi L(s, \chi)$ lie outside $\Dcz$), then we can define $\mathcal T(s, \chi_e)$ exactly as we defined the functions $\{\Tschi\}_{\chiNeChiZeroChie}$: In this case, $\mathcal T(s, \chi_e)$ is analytic on $\Dcz$, so that the function $s^{-1}\exp\left(\alphachiZ \TStschiZ + \sum_{\chi \ne \chi_0}\,\alpha_\chi \Tschi\right)$ analytically continues $s^{-1} \Fsnu (s-1/\nu)^{\alphachiZ}$ into $\Dcz$. Finally if $\eta_e$ exists but $\eta_e \le 1-c_0/10\lambdaqlogq$, then $\mathcal T(s, \chi_e)$ is analytic on the smaller region $\boldsymbol{\mathcal D(c_0/10\lambda_q) = \{\sigma+it: \sigma > \nu^{-1}(1-c_0/10\lambdaqElqt)\}}$; as such, the function $s^{-1} \Fsnu (s-1/\nu)^{\alphachiZ}$ continues analytically into $\mathcal D(c_0/10\lambda_q)$.       

The reader may now forget the $\mathcal T$ and $\mathcal T^*$. All that we will need  
from subsection \cref{subsec:AnalyticContinuations} are identities \eqref{eq:FsnuDef} and \eqref{eq:FsnuLogDerivAnaCont}, 
that $\Fsnu$ \textit{always} continues analytically into $\DczBigReduced$, and that 
\begin{equation}\label{eq:Hfunc}\allowdisplaybreaks
\Hfuncs \coloneqq 
\begin{cases} \vspace{2mm}
\displaystyle{\frac{\Fsnu}s \left(s-\frac1\nu\right)^{\alphachiZ} \left(s-\frac{\eta_e}\nu\right)^{-\alphachie}}\text{ cont.~an.~into }\Dcz,&\text{ if }\displaystyle{\eta_e > 1-\frac{c_0}{10\lambdaqlogq}}.\\ \vspace{2mm} 
\displaystyle{\frac{\Fsnu}s \left(s-\frac1\nu\right)^{\alphachiZ}}\text{ cont.~an.~into }\Dcz,&\text{ if }\eta_e\text{ doesn't exist.}\\ \vspace{2mm}
\displaystyle{\frac{\Fsnu}s \left(s-\frac1\nu\right)^{\alphachiZ}}\text{ cont.~an.~into }\mathcal \displaystyle{\mathcal D\left(\frac{c_0}{10\lambda_q}\right)},&\text{ if }\displaystyle{\eta_e \le 1-\frac{c_0}{10\lambdaqlogq}}.
\end{cases}
\end{equation}
Here ``cont.~an." abbreviates ``continues analytically". In what follows, we will call the three cases above as \textbf{Case 1, Case 2} and \textbf{Case 3} respectively. (It is \textit{not} necessary to assume that $\eta_e>1-c_0/10\lambdaqlogq$ for the function $s^{-1} \Fsnu (s-1/\nu)^{\alphachiZ} (s-\eta_e/\nu)^{-\alphachie}$ to continue analytically into $\Dcz$, however we do so in order to make this case trichotomy convenient for future use.)
\subsection{Analysis of Logarithmic Derivatives}
To give suitable bounds on $\Fsnu$, we will first analyze its logarithmic derivative.  
To this end, the following known results on Dirichlet $L$-functions will be useful.  
Recall that we write  $\boldsymbol{\rho = \beta+i\gamma}$ where $\boldsymbol{\beta = \Ree(\rho)}$ and $\boldsymbol{\gamma = \Imm(\rho)}$. We denote by $\boldsymbol{\sum_{\rho:\,L(\rho, \chi)=0}^*}~~$ 
a sum over all zeros $\rho$ of $L(s, \chi)$ counted with appropriate multiplicity. 
\begin{lem}\label{lem:LFuncStandardResults}
The following hold uniformly in $q \ge 2$ and in \textbf{all} Dirichlet characters $\chi$ mod $q$.\\    
$\mathrm{\boldsymbol{(1)}}$ Uniformly in all real $t$, we have $\sideset{}{^*}\sum\limits_{\substack{\rho:\,L(\rho, \chi) = 0\\0 \le \beta \le 1}}~~~~\displaystyle{\frac1{1+(t-\gamma)^2}} ~\ll~ \log(q(|t|+1)).$

$\mathrm{\boldsymbol{(2)}}$ Uniformly in all complex $s$ satisfying $\sigma \in [-1, 2]$, $|t| \ge 2$, and $t \ne \gamma$ for any of the zeros $\rho = \beta+ i \gamma$ of $L(s, \chi)$, we have $\displaystyle{\frac{L'(s, \chi)}{L(s, \chi)} ~=~ \sideset{}{^*}\sum\limits_{\substack{\rho:~L(\rho, \chi)=0\\ 0 \le \beta \le 1,~|\gamma-t|\le 1}}~\frac1{s-\rho} + O(\log(q(|t|+1)))}$.

$\mathrm{\boldsymbol{(3)}}$ We have $L'(s, \chi)/L(s, \chi) \ll \log(q|s|)$, uniformly in all complex $s$ satisfying $\sigma \le -1$ and lying outside the  disks of radius $1/4$ about the trivial zeros of $L(s, \chi)$. 

$\mathrm{\boldsymbol{(4)}}$ Uniformly in real $t \not\in (-1, 1)$, we have $\#\{\rho: 0 
\le \beta \le 1,~|\gamma-t| \le 1,~L(\rho, \chi)=0\} \ll \log(q|t|).$
\end{lem}
In most standard texts (such as \cite{davenport80, MV07, tenenbaum15}), these results are stated and proved only for primitive characters, however the generality above will be helpful here. 

\begin{figure}[h!]
\includegraphics[height=5cm]{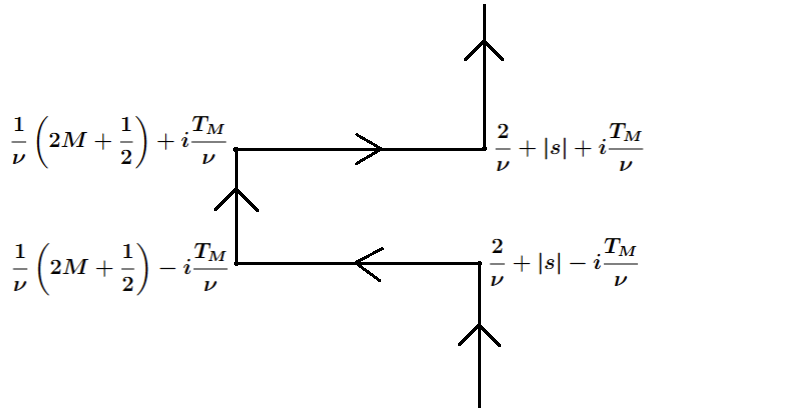}
\caption{The Contour $\omega_M$}
\label{Fig:LFuncLSDInnerContour}
\end{figure}

We now give a certain (absolutely convergent)  series expansion for the logarithmic derivative of $\Fsnu$ in terms of the zeros of the $L$-functions, with coefficients that are easy to control. 
\begin{prop}\label{prop:LogDerivSeries}
For any $s \in \C$ satisfying $s \ne 1/\nu$ and $s \ne \rho/\nu$ for any zero $\rho$ of $\prod_\chi L(s, \chi)$, 
\begin{equation}\label{eq:LogDerivSeries}
 \FsnuLogDeriv = \sum_{n \le \xi^2} \frac{\varrho(n) \Lambda(n)}{n^{s \nu}} \tau(n) ~-~ 
\frac{\alphachiZ(\xi^{1-\nu s} - \xi^{2(1-\nu s)})}{(1- \nu s)^2 \log \xi} ~+ \sum_{\chi \bmod q} \hspace{4mm} \sideset{}{^*}\sum_{\rho:\,L(\rho, \chi)=0}~~\frac{\alpha_\chi(\xi^{\rho-\nu s} - \xi^{2(\rho-\nu s)})}{(\rho - \nu s)^2 \log \xi}, 
\end{equation}
where ${\xi \coloneqq  e^{6\Elqt}}$, ${\varrho(n) \coloneqq \sum_{\chi \bmod q} ~~\alpha_\chi \cdot  \chi(n)}$, and $\tau(n) \coloneqq \bbm_{n \le \xi} ~+~ \bbm_{\xi < n \le \xi^2} ~(2-\log n/\log \xi)$.
\end{prop}
\begin{proof} Our starting point is the following identity,  which holds for any $b, y>0$ 
\begin{equation}\label{eq:GeneralContourIntegral}
\int_{b - i \infty}^{b + i \infty}~\frac{y^z}{z^2} \,\dz = \bbm_{y>1} \cdot 2 \pi i \log y    
\end{equation}
To see this, consider any $R \ge 2$, apply the residue theorem to the contour consisting of the vertical segment $[b-iR, b+iR]$ and the \textit{major} arc 
of the circle centered at the origin passing through $b \pm iR$ if $y>1$ (resp.~\textit{minor} arc for $y \le 1$), and then let $R \rightarrow \infty$. Now from $L'(s, \chi)/L(s, \chi) = -\sum_n \, \chi(n)\Lambda(n)/n^s$ and \eqref{eq:FsnuLogDerivAnaCont}, we see that ${\mathcal F'(z\nu)}/{\mathcal F(z\nu)}$ has Dirichlet series $\sum_n  \, {\varrho(n) \Lambda(n)}/{n^{z \nu}}$ on $\{z \in \C: \Ree(z)>1/\nu\}$. We claim that for all $s$ as in the statement of the proposition,  
\begin{equation}\label{eq:IntegralDirichletPolynomial}
\frac1{2\pi i} \int_{\frac2\nu + |s| - i \infty}^{\frac2\nu + |s| + i \infty}~~~\frac{\xi^{\nu(z-s)} - \xi^{2\nu(z-s)}}{(z-s)^2} \cdot \frac{\mathcal F'(z\nu)}{\mathcal F(z\nu)} \, \dz ~=~ \nu \sum_{n \le \xi^2}\frac{\varrho(n) \Lambda(n)}{n^{s\nu}} \tau(n) \log \xi.
\end{equation}
Indeed by \eqref{eq:GeneralContourIntegral}, 
the above identity is immediate if ${\mathcal F'(z\nu)}/{\mathcal F(z\nu)}$ were replaced by any finite truncation $\sum_{n \le Y}\,{\varrho(n) \Lambda(n)}/{n^{z \nu}}$ of its aforementioned Dirichlet series (for any $Y>\xi^6$). Moreover by the same Dirichlet series, the size of the integrand above is at most $2\lambda_q~\xi^{4+2\nu|s|} \left(\sum_n  \, \Lambda(n)/n^2\right) |z-s|^{-2}$, which is an $L^1$-function of $z$ since $\int_{2/\nu+|s|-i\infty}^{2/\nu+|s|+i\infty}~|\dz|/|z-s|^2 < \infty$ and $\sum_n \Lambda(n)/n^2 \ll 1$. Hence  \eqref{eq:IntegralDirichletPolynomial} follows from the Dominated Convergence Theorem. 

We will now shift contours: {This is the ``inner contour shift"  
alluded to in the introduction.} 
Note that for any $M \ge 2$, the number of zeros of $\prod_\chi\,L(s, \chi)$ in the rectangle $[0, 1] \times (M, M+1]$ is $\ll \phi(q)\log(qM)$ by Lemma \ref{lem:LFuncStandardResults}(4). Hence there exists $T_M \in (M, M+1]$ satisfying $|T_M - \gamma| \gg (\phi(q) \log(qM))^{-1}$ for \textbf{all} zeros $\rho = \beta + i\gamma$ of $\prod_\chi\,L(s, \chi)$. Since the set of zeros of $\prod_\chi\,L(s, \chi)$ is closed under complex conjugation, we have  
\begin{equation}\label{eq:TMMinusGamma}
|T_M \pm \gamma| \gg (\phi(q) \log(qM))^{-1}~~\text{for all zeros $\rho = \beta+i\gamma$ of $\prod_\chi\,L(s, \chi)$.}    
\end{equation}
With the contour $\omega_M$ as in Figure \ref{Fig:LFuncLSDInnerContour}, we claim that 
\begin{equation}\label{eq:LFuncLogDerivBound}
\frac{L'(z\nu, \chi)}{L(z\nu, \chi)} \ll \phi(q)\log^2(qM),\text{ uniformly in $q \ge 3$, $\chi$ mod $q$, $M \ge 2(1+\nu+\nu|s|)$, $z \in \omega_M$.}   
\end{equation}
If $\Ree(z) \ge 2/\nu$, this follows from the Dirichlet series of $L'(z\nu, \chi)/L(z\nu, \chi)$. If $\Ree(z) \in [-1/\nu, 2/\nu]$, then $z$ must lie on the two horizontal segments in $\omega_M$, so that by \eqref{eq:TMMinusGamma}, we have $|z\nu-\rho| \ge |\Imm(z)\nu-\gamma| = |T_M \pm \gamma| \gg (\phi(q)\log(qM))^{-1}$ for any zero $\rho = \beta+i\gamma$ of $\prod_\chi\,L(s, \chi)$. This gives  \eqref{eq:LFuncLogDerivBound} by Lemma \ref{lem:LFuncStandardResults}(2) and (4). Lastly if $\Ree(z) \le -1/\nu$, then Lemma \ref{lem:LFuncStandardResults}(3) establishes  \eqref{eq:LFuncLogDerivBound}. 

Now for any $M \ge 2\nu|s|$ and any $z \in \omega_M$, we have $|z-s| \ge |z|-|s| \ge |z|/2 \ge M/2\nu$. As such $\int_{\omega_M}~|\dz|/|z-s|^2 ~\ll_{\nu, s}~ \int_{M/2\nu}^\infty~\mathrm{d}t/t^2 ~+~ (M/2\nu)^{-2} \cdot M \ll M^{-1}$, so that \eqref{eq:FsnuLogDerivAnaCont} and \eqref{eq:LFuncLogDerivBound} yield
\begin{equation}\label{eq:OmegaMIntegralGoesTo0}
\lim_{M \rightarrow \infty}~~~~ \int_{\omega_M}~~~\frac{\xi^{\nu(z-s)} - \xi^{2\nu(z-s)}}{(z-s)^2} \cdot \frac{\mathcal F'(z\nu)}{\mathcal F(z\nu)} \, \dz ~=~ 0.
\end{equation}
Using the residue theorem to shift contours from the vertical line in \eqref{eq:IntegralDirichletPolynomial} to $\omega_M$, 
and then letting $M \rightarrow \infty$, we thus find from \eqref{eq:OmegaMIntegralGoesTo0}, \eqref{eq:IntegralDirichletPolynomial} and \eqref{eq:FsnuLogDerivAnaCont} that 
\begin{equation}\label{eq:DirPolyResidues}
\nu \sum_{n \le \xi^2}\frac{\varrho(n) \Lambda(n)}{n^{z\nu}} \tau(n) \log \xi ~=~ \left(\underset{z = s}\Res ~+~ \underset{z = 1/\nu}\Res ~+~ \sum_{\rho:~\prod_\chi L(\rho, \chi)=0} \hspace{3.5mm} \underset{z = \rho/\nu}\Res\right)~ \frac{\xi^{\nu(z-s)} - \xi^{2\nu(z-s)}}{(z-s)^2} \cdot \frac{\mathcal F'(z\nu)}{\mathcal F(z\nu)}.
\end{equation}
Finally, using \eqref{eq:FsnuLogDerivAnaCont} to compute the above residues, we obtain the proposition. For instance, note that if $\xi^{\rho-\nu s} \ne 1$ for some $\rho$ above, then \eqref{eq:FsnuLogDerivAnaCont} shows that $z=\rho/\nu$ is a simple pole of the function on the right of \eqref{eq:DirPolyResidues} of residue $\nu (\xi^{\rho-\nu s} - \xi^{2(\nu-\rho s)}) (\rho - \nu s)^{-2} \sum_\chi \,\alpha_\chi \cdot $\{multiplicity of $\rho$ in $L(s, \chi)$\}. If $\xi^{\rho-\nu s} = 1$, then $z=\rho/\nu$ is a removable singularity, so we can still give the same expression (whose value is zero) for its ``residue". The residue at $z=1/\nu$ can be computed analogously, and the residue at $z=s$ (which is always necessarily a simple pole) is equal to $-\nu (\log \xi) \Fprimesnu/\Fsnu$.   
\end{proof}
We will now use the series representation in Proposition \ref{prop:LogDerivSeries} to give a suitable bound on $\Fprimesnu/\Fsnu$. A crucial input will be provided the following zero density estimate. In what follows, we define $$\boldsymbol{N(\theta, t) ~~\coloneqq~~ \sum_{\chi}~~~ \sideset{}{^*}\sum_{\substack{\rho:~L(\rho, \chi)=0\\ \theta \le \beta \le 1,~|\gamma|\le t}}~1}.$$
\begin{lem}\label{lem:ZeroDensity}
We have $N(\theta, t) \ll (qt)^{3(1-\theta)}$, uniformly in $q \ge 3$, $\theta \in [1/2, 1]$, and $t \ge 1$.    
\end{lem}
{This may be found in works of Heath-Brown \cite{HealthBrown90} and Jutila \cite{jutila77}. (See also the classical text of Iwaniec and Kowalski  \cite{IwaniecKowalski}.) We now state the bound alluded to above.  
\begin{prop}\label{prop:LogDerivFsnuBound}
Uniformly in $q \ge 3$, and in complex numbers $s$ satisfying $\sigma \ge \nu^{-1}(1-c_0/2\Elqt)$,  
$$\left| \FsnuLogDeriv ~+~ \frac{\alphachiZ}{s\nu-1} ~-~ \frac{\alphachie}{s\nu-\eta_e} \right|~\ll~ \lambda_q\,\Elqt.$$
The term ``${\alphachie}/{(s\nu-\eta_e)}$" above is omitted if $\eta_e$ doesn't exist.   
\end{prop}
\begin{proof}
Most of the argument consists of carefully bounding the different components of the right of \eqref{eq:LogDerivSeries}. First, for all $n \le \xi^2$, we have $|n^{s\nu}| = n^{\sigma\nu} \ge n^{1-c_0/2\Elqt} \ge n \exp(-2\log \xi/2\Elqt) \gg n$, so that the first sum on the right in  \eqref{eq:LogDerivSeries} is $\ll \lambda_q \sum_{n \le \xi^2} \, \Lambda(n)/n \ll \lambda_q \, \Elqt$ by Mertens' Theorem.

Next, since the trivial zeros of any $L(s, \chi)$ are simple, the total contribution of all zeros $\{-r/\nu\}_{r \in \NatNos}$ to the right of \eqref{eq:LogDerivSeries} equals ${(\log \xi)^{-1}}~\sum_{r \ge 1}~ \left(\sum_{\chi: \, \chi(-1)=(-1)^r} \, \alpha_\chi \right)~{(\xi^{-(r+\nu s)} - \xi^{-2(r+ \nu s)})}{(r+ \nu s)^{-2}}.$
This expression is $\ll \lambda_q (\log \xi)^{-1} \sum_{r \ge 1}\,\xi^{-(r+\nu \sigma)}{(r+\nu\sigma)^{-2}} \ll \lambda_q \,  \Elqt^{-1} \sum_{r \ge 1}\,r^{-2} \ll $ $ \lambda_q \, \Elqt^{-1}$, where we noted that 
$|\sum_{\chi: \, \chi(-1)=(-1)^r} \, \alpha_\chi| \, = \, |\sum_{\chi}~\alpha_\chi (1 + \chi(-1)(-1)^r)/2| = |\varrho(1) + \varrho(-1)|/2 \le \lambda_q$.

Now, we observe that $|(\xi^{\theta - \nu s} - \xi^{2(\theta - \nu s)})(\theta - \nu s)^{-2} (\log \xi)^{-1} ~-~ (\nu s - \theta)^{-1}| \ll \Elqt$ uniformly in $\theta \in (0, 1]$ and $s$ as in the proposition. This follows by a straightforward crude bounding if $|\theta - \nu s| > (\log \xi)^{-1}$, and by the formula $\xi^{\theta - \nu s} = 1 - (\theta - \nu s) \log \xi + O\big((\theta - \nu s)^2 (\log \xi)^2\big)$ if $|\theta - \nu s| \le (\log \xi)^{-1}$. Collecting all the observations made so far, we see that this proposition would follow from \eqref{eq:LogDerivSeries}, once we show that uniformly in all $s$ with $\sigma \ge \nu^{-1}(1-c_0/2\Elqt)$, 
\begin{equation}\label{eq:SeriesZerosContrib}
\frac1{\Elqt}~~\sum_\chi ~|\alpha_{\chi}|~~~~\sideset{}{^*}\sum_{\substack{\rho:~L(\rho, \chi)=0\\0 \le \beta \le 1,~\rho \ne \eta_e}}~~~\frac{\xi^{\beta - \nu \sigma} + \xi^{2(\beta - \nu \sigma)}}{(\beta - \nu \sigma)^2 + (\gamma - \nu t)^2}~\ll~\lambda_q \Elqt.    
\end{equation}
To show this, we start by bounding the entire expression above by  $S_1+S_2+S_3+S_4$, where 
\begin{itemize}
\item $S_1$ denotes the total contribution of all $\rho$ having $\beta \le 1/2$, so that 
$$S_1 = \frac1{\Elqt}~~\sum_\chi ~|\alpha_{\chi}|~~~~\sideset{}{^*}\sum_{\substack{\rho:~L(\rho, \chi)=0\\0 \le \beta \le 1/2}}~~~\frac{\xi^{\beta - \nu \sigma} + \xi^{2(\beta - \nu \sigma)}}{(\beta - \nu \sigma)^2 + (\gamma - \nu t)^2}.$$
\item $S_2$ denotes the total contribution of all $\rho$ having $\beta \in (1/2, 1]$ and $|\gamma| \le 2|t\nu|+1$.
\item $S_3$ denotes the total contribution of all $\rho$ having $\beta \in (1/2, \sigma \nu]$ and $|\gamma| > 2|t\nu|+1$.
\item $S_4$ denotes the total contribution of all $\rho$ having $\beta \in (\sigma \nu, 1]$.
\end{itemize}
 
For any $\rho$ appearing in $S_1$, we have $\beta - \nu\sigma \le 1/2-(1-c_0/2\Elqt) \le -1/2+1/2\log q \le -1/3$, so that 
$(\beta - \nu \sigma)^2 + (\gamma - \nu t)^2 \ge (1+ (\gamma - \nu t)^2)/9$. Hence \eqref{eq:alphachilambdaq} and Lemma \ref{lem:LFuncStandardResults}(1) yield $S_1 \ll \lambda_q~\xi^{1/2 - \nu \sigma} \Elqt^{-1} \sum_\chi \sum_\rho~  (1+ (\gamma - \nu t)^2)^{-1}$ $\ll \lambda_q \cdot q\xi^{1/2 - \nu \sigma} \ll \lambda_q \cdot q\xi^{-1/2} \cdot\xi^{c_0/2\Elqt} \ll \lambda_q$. 

For any $\rho$ appearing in $S_3$, we have $\beta - \nu \sigma \le 0$ and $|\gamma - t\nu| \ge |\gamma|-|t\nu| \ge |\gamma|/2$. Thus by \eqref{eq:alphachilambdaq}, 
$$S_3 \le \frac{8 \lambda_q}{\Elqt}~~~~\sum_\chi \,~~~~~\sideset{}{^*}\sum_{\substack{\rho \ne \eta_e:~L(\rho, \chi)=0\\|\gamma|>2|t\nu|+1,~1/2<\beta\le \min\{\sigma \nu ,1\}}}~~~{\xi^{\beta-\nu\sigma}} \cdot {|\gamma|^{-2}}.$$
Partitioning the interval $(1/2, \min\{\sigma\nu, 1\}]$ into $R \coloneqq \lfloor \log \xi/2 \rfloor$ equally spaced intervals, we obtain
\begin{equation}\label{eq:S3Split}
S_3 \le \frac{8 \lambda_q}{\Elqt}~~~~\sum_{r=1}^R ~~~~\xi^{1/2 + r\mu_0/R - \nu\sigma}~~~~\sum_\chi \,~~~~~\sideset{}{^*}\sum_{\substack{\rho \ne \eta_e:~L(\rho, \chi)=0,~|\gamma|>2|t\nu|+1\\1/2 + (r-1)\mu_0/R \, < \, \beta \, \le \, 1/2 + r\mu_0/R}}~~~{|\gamma|^{-2}},
\end{equation}
where $\mu_0 \coloneqq \min\{\sigma\nu, 1\}-1/2$. Now the inner double sum (on $\chi$ and $\rho$) above is at most 
\begin{equation}\label{eq:S3RiemStiel}
\int_{2|t\nu|+1}^\infty \, \frac{\mathrm d N\left(\frac12 + \frac{(r-1)\mu_0}R, u\right)}{u^2} \ll \int_{2|t\nu|+1}^\infty~~\frac{N\left(\frac12 + \frac{(r-1)\mu_0}R, u\right)}{u^3} \, \mathrm du \ll q^{3(1/2-(r-1)\mu_0/R)} 
\end{equation}
where we have used the Stieltjes integration by parts and Lemma \ref{lem:ZeroDensity}. The last expression above is $\ll \xi^{1/4 - (r-1)\mu_0/2R} \ll \xi^{1/4-r \mu_0/2R}$, as $\mu_0 \le 1/2$, $\xi \ge q^6$ and $R \ge \log \xi/3$. Inserting these into \eqref{eq:S3Split}, we get $S_3 \ll \lambda_q\,\Elqt^{-1} \, \xi^{3/4-\nu \sigma} ~~\sum_{r=1}^R~~ \xi^{r \mu_0/2R} ~\le~ \lambda_q\,\Elqt^{-1}~ \xi^{3/4-\nu \sigma} \cdot R\xi^{\mu_0/2} \ll \lambda_q\,\xi^{1-\nu\sigma} \ll \lambda_q$.

Next, for any $\rho$ counted in $S_2$, we have $|\gamma| \le 2|t\nu|+1$ and $\rho \ne \eta_e$, so that $\beta \le 1 - c_0/\log(q(|\gamma|+1)) \le 1-c_0/\log(2q(|t\nu|+1))$. Since $\nu\sigma \ge 1-c_0/2\Elqt$, we get 
$$\nu\sigma - \beta \ge c_0\left(\frac1{\log(2q(|t\nu|+1))} - \frac1{2\log(q(|t\nu|+1))}\right) = \frac{c_0}{\Elqt}\left(1-\frac{\log 4}{\log(2q(|t\nu|+1))}\right) \ge \frac{c_0}{10\Elqt}.$$
Hence $(\beta - \nu\sigma)^2 \gg \Elqt^{-2}$. Proceeding as in \eqref{eq:S3RiemStiel} (via Lemma \ref{lem:ZeroDensity} and integration by parts), 
\begin{align*}\allowdisplaybreaks
S_2 &\ll \lambda_q\,\Elqt~~\sum_\chi \,~~~\sideset{}{^*}\sum_{\substack{\rho \ne \eta_e:~L(\rho, \chi)=0\\1/2<\beta\le 1,~|\gamma| \le 2|t\nu|+1}}~{\xi^{\beta-\nu\sigma}}~\le~ \lambda_q \, \Elqt \, \left(-\int_{1/2}^{1}~\xi^{\theta-\nu\sigma} \, \mathrm dN(\theta, 2|t\nu|+1)\right)\\
&\le~ \lambda_q \, \Elqt \,\left( \xi^{1/2-\nu\sigma} \, N(1/2, 2|t\nu|+1) \, + \, \log \xi \, \int_{1/2}^{1} \xi^{\theta-\nu\sigma} \, N(\theta, 2|t\nu|+1) \, \mathrm d\theta \right)\\ &\ll~ \lambda_q \, \Elqt \, \left(\xi^{3/4-\nu\sigma} \, + \, \log \xi \, \int_{1/2}^{1} \xi^{(1+\theta)/2-\nu\sigma} \, \mathrm d\theta \right) \, \ll \, \lambda_q \, \Elqt \, \xi^{1-\nu\sigma} \, \ll \lambda_q \, \Elqt. 
\end{align*}
For any $\rho$ ($\ne \eta_e$) in $S_4$, we have $1-c_0/\log(q(|\gamma|+1)) \ge \beta$ $> \sigma\nu \ge 1-c_0/2\Elqt$, giving $|\gamma| > q(|t\nu|+1)^2-1$. Thus also $|\gamma-t\nu| \ge |\gamma|-|t\nu| \ge |\gamma|/2$. Proceeding exactly as we did for $S_3$, 
$$S_4 \, \ll \, \frac{\lambda_q \, \xi^{2(1-\nu\sigma)}}{\Elqt}~~\sum_\chi \,~~~~\sideset{}{^*}\sum_{\substack{\rho: \, L(\rho, \chi)=0\\\sigma\nu < \beta < 1,~|\gamma| > q(|t\nu|+1)^2-1}}~~~{|\gamma|^{-2}} \, \ll \, \frac{\lambda_q}{\Elqt}~~\int_{q(|t\nu|+1)^2-1}^\infty~\frac{\mathrm dN(\sigma\nu, u)}{u^2} \, \ll \, \lambda_q.$$
Collecting all these estimates establishes \eqref{eq:SeriesZerosContrib}, completing the proof of the proposition. \end{proof}
As a consequence of Proposition \ref{prop:LogDerivFsnuBound},  
we find that  
the size of the function $\Hfuncs$ (in \eqref{eq:Hfunc}) remains roughly constant along \textit{short} \textit{horizontal} segments lying in the region $\sigma \ge \nu^{-1}(1-c_0/2\lambdaqElqt)$.    
\begin{cor}\label{cor:HfuncsVariation}
In Cases 1 and 2 of \eqref{eq:Hfunc}, we have $\Hfuncs \asymp \Hfunc(w)$, uniformly in $s, w \in \C$ having $t = \Imm(s) = \Imm(w)$ and $\nu^{-1}(1-c_0/2\lambdaqElqt) \le \Ree(s) \le \Ree(w) \le \nu^{-1}(1+1000/c_0\lambdaqElqt)$. 

In Case 3, the same assertion holds with $\nu^{-1}(1-c_0/2\lambdaqElqt)$ replaced by $\nu^{-1}(1-c_0/\boldsymbol{20}\lambdaqElqt)$.   
\end{cor}
\begin{proof}
In Cases 1 and 2 of \eqref{eq:Hfunc},  Proposition \ref{prop:LogDerivFsnuBound} directly gives $|\mathcal H'(z)/\Hfunc(z)| \ll \lambdaqElqt$ for all $z$ satisfying  $\Ree(z) \ge \nu^{-1}\big(1-c_0/2\lambda_q\,\Elq(\Imm(z))\big)$. In Case 3, we just need note that for all $z$ satisfying  $\Ree(z) \ge \nu^{-1}\big(1-c_0/20\lambda_q\,\Elq(\Imm(z))\big)$, we have $|z\nu - \eta_e| \ge \Ree(z)\nu-\eta_e \ge c_0/20\lambdaqlogq$, so that Proposition \ref{prop:LogDerivFsnuBound} still gives $|\mathcal H'(z)/\Hfunc(z)| \ll \lambdaqElqt$ for all such $z$. The corollary now follows by writing $\log|\Hfunc(w)/\Hfuncs| \le \int_{\Ree(s)}^{\Ree(w)} |\mathcal H'(u+it)/\Hfunc(u+it)| \, \mathrm{d}u$ for all $s, w$ as in the statement.   
\end{proof}
\subsection{Bounds on $\Hfuncs$ and $\Fsnu$.} 
Corollary \ref{cor:HfuncsVariation} allows us to  study the values of $\Hfuncs$ for $\Ree(s) \le 1/\nu$ using the values of $\mathcal H(w)$ for  
$\Ree(w)>1/\nu$.  
The benefit of this is that on this right half plane, $\mathcal F(w\nu)$ has an \textit{explicit expression} \eqref{eq:FsnuDef} as the exponential of a Dirichlet series. This maneuver allows us to bound the sizes of both $\Hfuncs$ and $\Fsnu$ in suitable regions,  
setting the stage for the ``outer contour shift" part of the ``nested contour shift" alluded to in the introduction.  
\begin{prop}\label{prop:HsFsnuBounds} The following assertions hold in Cases 1 and 2 of \eqref{eq:Hfunc}, with all terms involving $\eta_e$ or $\chi_e$ omitted in Case 2 (when $\eta_e$ doesn't exist).

$\mathrm{\boldsymbol{(1)}}$ We have $\Hfuncs \ll \TwoLambdaqLogqPower$, uniformly in moduli $q>4e^{1/\nu}$ and in complex numbers $s=\sigma+it$ satisfying $\sigma \ge \nu^{-1}(1-c_0/2\lambda_q \Elqt)$ and $|t| \le c_0/2\nu \lambda_q \log q$.

$\mathrm{\boldsymbol{(2)}}$ We have $\Fsnu \ll (5 \lambda_q \Elqt/4)^{\lambda_q}$, uniformly in moduli $q>e^{4+1/\nu}$ and in  $s=\sigma+it$ satisfying $\nu^{-1}(1-c_0/2\lambda_q \Elqt) \le \sigma \le {\nu^{-1}(1 + 100/c_0\lambda_q \Elqt)}$ and $\min\{|s-1/\nu|, |s-\eta_e/\nu|\} \ge c_0/40\nu \lambda_q \log q$.

In Case 3 of \eqref{eq:Hfunc}, these two assertions also hold, with all terms involving $\eta_e$ or $\chi_e$ omitted, and with all occurrences of ``$\nu^{-1}(1-c_0/2\lambda_q \Elqt)$" replaced by ``$\nu^{-1}(1-c_0/\boldsymbol{20}\lambda_q \Elqt)$". 
\end{prop}
\begin{proof}
We only give the argument in Case 1, since the arguments in the other two cases are essentially contained in it. Define $\mu(y) \coloneqq \nu^{-1}(1+4/5\lambda_q \Elq(y))$ for any $y \in \reals$.\footnote{Not to be confused with the M\"{o}bius function, which makes no appearance here.} By \eqref{eq:Hfunc} and Corollary \ref{cor:HfuncsVariation}, we have, uniformly in all $s = \sigma+it$ satisfying $\sigma \ge \nu^{-1}(1-c_0/2\lambdaqElqt)$,  
\begin{equation}\label{eq:HfuncBoundStep1}
\Hfuncs \ll |\Hfunc(\mu(t)+it)| \ll \frac{|\mathcal F\big(\nu (\mu(t)+it)\big)|}{|\mu(t)+it|} \cdot \left|\mu(t)+it - \frac1\nu\right|^{\Ree(\alphachiZ)} \cdot \left|\mu(t)+it - \frac{\eta_e}\nu\right|^{-\Ree(\alphachie)}.   
\end{equation} 
Now since $\Ree(\mu(y)+iy) > 1/\nu$, it follows from 
\eqref{eq:FsnuDef} and $|\sum_\chi \, \alpha_\chi\cdot\chi(p^r)| \le \lambda_q$ that 
\begin{align}\label{eq:FsnuDirSeriesBound} \allowdisplaybreaks
|\mathcal F\big(\nu(\mu(y)+iy)\big)| &\le \exp\left(\lambda_q\sum_{p, r \ge 1}~ \frac1{p^{r \nu \mu(y)}}\right) = \exp\Big(\lambda_q \, \log \zeta(\nu\mu(y)) \Big) \ll \left(\frac54 \lambda_q\,\Elq(y)\right)^{\lambda_q}.   
\end{align}
Here the last bound uses  
that since $\log(\nu \mu(y)) = \log(1+4/5\lambda_q \Elq(y)) \le 4/5\lambda_q \Elq(y) \le 4/5\lambda_q \log q \le 1/\lambda_q$, we have $\lambda_q \, \log \zeta(\nu\mu(y)) \le \lambda_q \, \log({\nu \mu(y)}) - \lambda_q \,\log({\nu \mu(y)-1}) \le 1+\lambda_q \, \log \big(5\lambda_q\,\Elq(y)/4\big)$. 

Inserting \eqref{eq:FsnuDirSeriesBound} into \eqref{eq:HfuncBoundStep1}, we obtain,  
uniformly in all $s=\sigma+it$ with $\sigma \ge \nu^{-1}(1-c_0/2\lambdaqElqt)$,   
\begin{equation}\label{eq:HfuncBoundStep2}
\Hfuncs \ll \frac{(5\lambdaqElqt/4)^{\lambda_q}}{|\mu(t)+it|} \cdot \left|\mu(t)+it - \frac1\nu\right|^{\Ree(\alphachiZ)} \cdot \left|\mu(t)+it - \frac{\eta_e}\nu\right|^{-\Ree(\alphachie)}.   
\end{equation}
\textbf{Completing the proof of subpart \textbf{(1)}.} We now observe that $|\mu(t)+it-1/\nu|$ and $|\mu(t)+it-\eta_e/\nu|$ both lie in the interval $(4/5\nu \lambdaqElqt, 1)$: Indeed, the lower bounds are immediate by definition of $\mu(t)$, while the upper bounds follow from the facts that $\eta_e>1-c_0/10\lambdaqlogq$, that $q>e^{1/\nu}$ and $|t| \le c_0/2\nu\lambda_q \log q$ (by the assumptions in Case 1 in \eqref{eq:Hfunc} and in the statement of subpart \textbf{(1)}). 

This observation yields $|\mu(t)+it-1/\nu|^{\Ree(\alphachiZ)} \le |\mu(t)+it-1/\nu|^{-|\Ree(\alphachiZ)|} \le (5 \nu \lambdaqElqt/4)^{|\Ree(\alphachiZ)|}$, and likewise $|\mu(t)+it-\eta_e/\nu|^{-\Ree(\alphachie)} \le (5 \nu \lambdaqElqt/4)^{|\Ree(\alphachie)|}$. 
Inserting these two bounds into \eqref{eq:HfuncBoundStep2}, we get $\Hfuncs \ll (5\lambdaqElqt/4)^{\lambda_q + |\Ree(\alphachiZ)| + |\Ree(\alphachie)|}$. 
Finally,  
$\Elqt = \log(q(|t\nu|+1)) < \log(q(c_0/2+1)) \le \log(9q/8) < (8/5)\log q$. (Here we just  used $|t|<c_0/2\nu$, $c_0<1/4$, and $q>2$.)

\textbf{Completing the proof of subpart \textbf{(2)}.} By \eqref{eq:Hfunc} and \eqref{eq:HfuncBoundStep2}, we have    
\begin{equation}\label{eq:FsnuBoundNearFinal}
\Fsnu \ll {\left(\frac54\lambdaqElqt\right)^{\lambda_q}}\cdot \frac{|s|}{|\mu(t)+it|} \cdot \left|\frac{\mu(t)+it - 1/\nu}{s-1/\nu}\right|^{\Ree(\alphachiZ)} \cdot \left|\frac{s-{\eta_e}/\nu}{\mu(t)+it - {\eta_e}/\nu}\right|^{\Ree(\alphachie)}   
\end{equation}
uniformly in all $s$ with $\sigma \ge \nu^{-1}(1-c_0/2\lambdaqElqt)$. Now for $s$ as in subpart \textbf{(2)}, we have $|s|  
\le \sqrt{4/\nu^2+t^2} \ll |\mu(t)+it|$,  
as well as $\mu(t)-\sigma \ll 1/\lambdaqElqt$, and $\sigma-\mu(t) \le \sigma - 1/\nu \ll 1/\lambdaqElqt$, and $\min\{|s-1/\nu|, |s-\eta_e/\nu|\} \gg 1/\lambdaqElqt$.  
These last three inequalities give, for both $\theta \in \{1, \eta_e\}$,   
$$\max\left\{\left|\frac{\mu(t)+it - \theta/\nu}{s-1/\nu}\right|, \, \left|\frac{s-\theta/\nu}{\mu(t)+it - 1/\nu}\right|\right\} \le 1+\max\left\{\left|\frac{\mu(t)-\sigma}{s-\theta/\nu}\right|, \, \left|\frac{\mu(t)-\sigma}{\mu(t)+it - \theta/\nu}\right|\right\} \ll 1,$$
Inserting all these observations into \eqref{eq:FsnuBoundNearFinal} completes the proof of the proposition. 
\end{proof}
\section{The LSD method for $L$--functions under  
average growth conditions:\\ Proofs of Theorem \ref{thm:LFuncLSDVariant1} and Corollary \ref{cor:LFuncLSDVariant1LogPowerInterval} }\label{sec:LfunctionsLSDVariant1Proof}
In this section, we will show the following generalizations of Theorem \ref{thm:LFuncLSDVariant1} and Corollary \ref{cor:LFuncLSDVariant1LogPowerInterval}.
\begin{thm}\label{thm:LFuncLSDVariant1Gen}
Under the conditions of Theorem \ref{thm:LFuncLSDVariant1}, the following estimates hold, with the same uniformity as in Theorem \ref{thm:LFuncLSDVariant1}, but also uniformly in $T \ge 1$ satisfying ${c_0 \lambda_q \Elq(T) \le 100\log x}$.

$\mathrm{\boldsymbol{(1)}}$ If $\eta_e$ exists and satisfies $1-c_0/10\lambdaqlogq < \eta_e < 1-3\nu/\log x$ 
then the left side of \eqref{eq:LFuncLSDVar1SiegelZero} is 
\begin{multline}\label{eq:LFuncLSDVar1SiegelZeroGen}\allowdisplaybreaks
\ll~ \frac{x^{1+1/\nu} \, \log x}{Th} ~+~ \left(\frac54 \lambda_q \Elq(T)\right)^{\lambda_q} \, \Omega(T) \left\{\frac{x^{1/\nu}}{T\log x} ~+~ x^{1/\nu - c_0/4\nu \lambda_q\Elq(T)} (\log T)\right\} 
\\~+~ \sum_{x < n \le x+h} \, |a_n| \, +~~ \frac{\Omega(1) \cdot N! \, (71(1+\nu))^N \cdot 2^{\lambda_q} \, (\lambda_q \log q)^{\lambda_q + |\Ree(\alphachiZ)| + |\Ree(\alphachie)|} \cdot x^{1/\nu}}{(1-\eta_e)^{N+1+|\Ree(\alphachie)|} \cdot (\log x)^{1-|\Ree(\alphachiZ)|} \cdot \min\{x/h, (\log x)^{N+1}\}}.
\end{multline} 
$\mathrm{\boldsymbol{(2)}}$ If $\eta_e$ does not exist or satisfies $\eta_e \le 1- c_0/10\lambdaqlogq$, then for $q < x^{c_0/80\nu\lambda_q}$, the left of \eqref{eq:LFuncLSDVar1NoSiegelZero} is bounded by the same expression as \eqref{eq:LFuncLSDVar1SiegelZeroGen}, but only with ``${-c_0/4\nu \lambda_q\Elq(T)}$" replaced by ``${-c_0/\boldsymbol{80}\nu \lambda_q\Elq(T)}$", and with the term involving $\Omega(1)$ replaced by 
$${\frac{\Omega(1) \cdot N! \, (2000(1+\nu) c_0^{-1})^N \cdot 2^{\lambda_q} \, (\lambda_q \log q)^{\lambda_q +  |\Ree(\alphachie)|} \cdot x^{1/\nu}}{(\lambda_q \log q)^{-N} \cdot (\log x)^{1-|\Ree(\alphachiZ)|} \cdot \min\{x/h, (\log x)^{N+1}\}}.}$$
{The assertions corresponding to the last sentence of Theorem \ref{thm:LFuncLSDVariant1} also holds, if $\eta_e$ does not exist.} 

Finally under \eqref{eq:Variant1GrowthConditionSpecial}, the above assertions hold exactly as stated, with $h=x/(\log x)^A$. 
\end{thm}
Theorem \ref{thm:LFuncLSDVariant1} and Corollary \ref{cor:LFuncLSDVariant1LogPowerInterval} follow from the above results by taking $T$ as in the respective statements. (The choice of $T$ in Theorem \ref{thm:LFuncLSDVariant1}(1) comes from  
writing $1/T \approx x^{-c_0/4\nu \lambda_q\Elq(T)}$ and  
$\Elq(T) \approx (\log T)^2 + \log(q\nu) \cdot (\log T)$, and then solving the quadratic in $\log T$. There may be better ways of choosing $T$ for specific $\Omega(T)$, which is why we prove Theorem \ref{thm:LFuncLSDVariant1Gen} in its additional generality.)

\subsection{Perron's Formula: Error--term control via averaging}
Our proof of Theorem \ref{thm:LFuncLSDVariant1Gen} begins by 
relating the partial sum $\sum_{n \le x}\, a_n$ with a contour integral via Perron's formula. The following lemma will allow us to control the error terms when we apply this formula.  
\begin{lem}\label{lem:PerronErrorController}
Assuming \textbf{only} \eqref{eq:Variant1GrowthConditionNew}, the following  
hold uniformly in  $x>4e^{1/\nu}$ and $h \in (0, x/2]$.

$\mathrm{\boldsymbol{(1)}}$ For all $\theta \in (0, 1]$, we have $\displaystyle{\sum_{n>1} \, {|a_n|}/{n^{1/\nu+\lambda}} \, \le \, 4\kappa/\theta}$.\\ 

$\mathrm{\boldsymbol{(2)}}$ There exists a \textbf{half--integer}\footnote{i.e. an element of the set ${\Z+1/2 = \left\{\pm 1/2, \pm 3/2, \dots\right\}}$} $X \in (x, x+h]$  satisfying $$\displaystyle{\sum_{\frac{3X}4 < n < \frac{5X}4}~~~~  \frac{|a_n|}{|\log(X/n)|} \ll \, \kappa \cdot \frac{x^{1+1/\nu} \, \log x}{h}}.$$
\end{lem}
\begin{proof}
\textbf{(1)} Indeed by \eqref{eq:Variant1GrowthConditionNew}, we have $\sum_{n \ge 2} \, {|a_n|}/{n^{1/\nu+\theta}} \, \le \, \sum_{m \ge 0} \, 2^{-m/\nu-m \theta} \, \sum_{2^m< n \le 2^{m+1}} \, {|a_n|} \, \le \, \kappa \sum_{m \ge 0} \, 2^{-m\theta} = \kappa/(1-2^{-\theta})$. Now basic calculus shows that $2^{-\theta} \le 1-\theta/4$ for all $\theta \in (0, 1]$.

\textbf{(2)} It suffices to show that uniformly in $h \in (0, x/2]$, we have 
\begin{equation}\label{eq:Averaging3X4,5X4}
\sum_{\substack{x < X \le x+h\\X \in \Z+1/2}}~~~~\sum_{\frac{3X}4 < n < \frac{5X}4}~~~~ \displaystyle{\frac{|a_n|}{\left|\log\left(X/n\right)\right|}} ~~\ll~~ \kappa \cdot x^{1+/\nu} \log x. 
\end{equation} 
Write the total double sum on the left as $S_1 + S_2$, where $S_1$ denotes the total contribution of all pairs $(X, n)$ for which $n \in (3X/4, X-1/2]$. Then for any $n$ counted in $S_1$, we can write $n = X-r$ for some half--integer $r \in [1/2, X/4) \subset [1/2, (x+h)/4)$. Moreover, $n = X-r \in (x-r, (x-r)+h]$ and $|\log(X/n)| = \log(X/(X-r)) = -\log(1-r/X) \gg r/X \gg r/x$. Combining these observations, 
\begin{equation}\label{eq:S1Reduction}
S_1 ~\stackrel{\Delta}=~~ \sum_{\substack{x < X \le x+h\\X \in \Z+1/2}}~~~~\sum_{\frac{3X}4 < n \le X-\frac12}~~~~ \displaystyle{\frac{|a_n|}{\left|\log\left(X/n\right)\right|}} ~~\ll~~ x \sum_{\substack{1/2 < r \le {(x+h)}/4\\r \in \Z+1/2}}~~~\frac1r~~\sum_{x-r < n \le (x-r)+h}~~~~ \displaystyle{{|a_n|}}.   
\end{equation}
Now for any $r$ above, $x-r>x-(x+h)/4 > h$. Hence  
using \eqref{eq:Variant1GrowthConditionNew} on each inner sum in \eqref{eq:S1Reduction}, 
$$S_1 ~\ll ~~ \kappa \cdot x \sum_{\substack{1/2 < r \le {(x+h)}/4\\r \in \Z+1/2}}~~~\frac1r \cdot {(x-r)^{1/\nu}}  ~\ll~ {\kappa  \cdot x^{1+1/\nu}} \sum_{\substack{m \le x}}~~~\frac1{m} \ll {\kappa \cdot x^{1+1/\nu}} \, \log x,$$
proving that $S_1$ is absorbed in the right of  \eqref{eq:Averaging3X4,5X4}. The argument for $S_2$ is entirely analogous, 
by writing $n=X+r$ for some $r \in [1/2, (x+h)/4)$. 
This establishes \eqref{eq:Averaging3X4,5X4}, and hence the lemma. 
\end{proof}
\begin{rmk}
We need the growth condition \eqref{eq:Variant1GrowthConditionNew} \textit{only} in the form of Lemma \ref{lem:PerronErrorController}. 
\end{rmk}
We will first show Theorem \ref{thm:LFuncLSDVariant1Gen} with $x$ replaced by 
the $X$ coming from Lemma \ref{lem:PerronErrorController}(2); for this we will 
only need that $\{a_n\}_n$ 
has property $\mathcal P(\nu, \{\alpha_\chi\}_\chi; c_0, \Omega)$ for a \textit{general non--decreasing} $\Omega: \reals_{\ge 0} \rightarrow \reals_{\ge 0}$ for all $m \ge 1$. By Perron's Formula (as stated in \cite[Theorem II.2.3]{tenenbaum15}), we have   
\begin{equation}\label{eq:PerronApp1}
\sum_{n \le X} \, a_n ~=~ \frac1{2\pi i}\int_{\frac1\nu\left(1+\frac1{\log X}\right) - iT}^{\frac1\nu\left(1+\frac1{\log X}\right) + iT} \, \frac{\Fsnu G(s) X^s}s \, \mathrm ds ~+~ O\left(\frac{X^{1/\nu}}T \, \sum_{n \ge 1} \, \frac{|a_n|}{n^{\frac1\nu + \frac1{\nu\log X}} |\log(X/n)|}\right).
\end{equation}
By Lemma \ref{lem:PerronErrorController}(2), the total contribution of all $n \in (3X/4, 5X/4)$ to the $O$-term above is $\ll \kappa  \cdot X^{1/\nu} (\log X)/T$. On the other hand, Lemma \ref{lem:PerronErrorController}(1) shows that the total contribution of all $n \not\in (3X/4, 5X/4)$ to the $O$-term is $\ll X^{1/\nu} \, T^{-1} \, (|a_1|/\log X \,  
+ \, \kappa \cdot X \log X/h)$. But now letting $s \rightarrow +\infty$ (along the real line) in \eqref{eq:PropertP} and \eqref{eq:FsnuDef}, and using that $|G(s)| \le \Omega(1)$ for all real $s>1/\nu$, we obtain $|a_1| \le \Omega(1)$.  
Inserting all these observations into \eqref{eq:PerronApp1}, we obtain 
\begin{equation}\label{eq:PerronApp2}
\sum_{n \le X} \, a_n ~=~ \frac1{2\pi i}\int_{\frac1\nu\left(1+\frac1{\log X}\right) - iT}^{\frac1\nu\left(1+\frac1{\log X}\right) + iT} \, \frac{\Fsnu G(s) X^s}s \, \mathrm ds ~+~ O\left(\frac{\Omega(1) \, X^{1/\nu}}{T\log X}~+~ \kappa  \cdot \frac{X^{1+1/\nu} \log X}{Th}\right).
\end{equation} 
The rest of the argument breaks up into the two cases, in the two subparts of Theorem \ref{thm:LFuncLSDVariant1Gen}. 
\subsection{When $\boldsymbol{1-c_0/10\lambdaqlogq < \eta_e < 1-3\nu/\log x}$: Proof of Theorem \ref{thm:LFuncLSDVariant1Gen}(1)}~\\ In this subsection, we define 
\begin{equation}\label{eq:sigmanurerDefs}
\boldsymbol{\sigmanu(t) \coloneqq \frac1\nu\left(1-\frac{c_0}{4\nu\lambdaqElqt}\right),~~r_e = \frac{1-r_e}{6\nu},~~r_1 = \frac1{2\max\{1, \nu\} \log X}}.
\end{equation}
We also define $\boldsymbol{\Gamma_0}$ to be the contour consisting of the following components. (See Figure \ref{Fig:LFuncLSDOuterContourVar1.1}.)  
\begin{itemize}
\item $\boldsymbol{\Gamma_2}$, the horizontal segment traversed from $\sigmanu(T)+iT$ to $\nu^{-1}(1+1/\log x)+iT$.
\item $\boldsymbol{\Gamma_3}$, the part of the curve $\sigmanu(T)+it$ traversed upwards from $t=0$ to $t=T$. 
\item $\boldsymbol{\Gamma_4}$, the horizontal segment traversed from $\eta_e/\nu-r_e$ to $\sigmanu(0)$ \textbf{above} the branch cut.
\item $\boldsymbol{\Gamma_5}$, the anticlockwise semicircle in the upper half plane with center $\eta_e/\nu$, radius $r_e$. 
\item $\boldsymbol{\Gamma_6}$, the horizontal segment traversed from $(2+\eta_e)/3\nu$ to $\eta_e/\nu+r_e$  \textbf{above} the branch cut.
\item $\boldsymbol{\Gamma_7}$, the horizontal segment traversed from $1/\nu-r_1$ to $(2+\eta_e)/3\nu$   \textbf{above} the branch cut.
\item $\boldsymbol{\Gamma_8}$, the circle with center at $1/\nu$, radius $r_1$, traversed anticlockwise as shown in Figure \ref{Fig:LFuncLSDOuterContourVar1.1}. 
\item $\boldsymbol{\overline\Gamma_j}$ \textbf{(for} $\boldsymbol{2 \le j \le 7}$\textbf{)}, the reflection of $\Gamma_j$ about the real line, directed as in Figure \ref{Fig:LFuncLSDOuterContourVar1.1}.
\end{itemize}
\begin{figure}[h!]
\includegraphics[height=11cm]{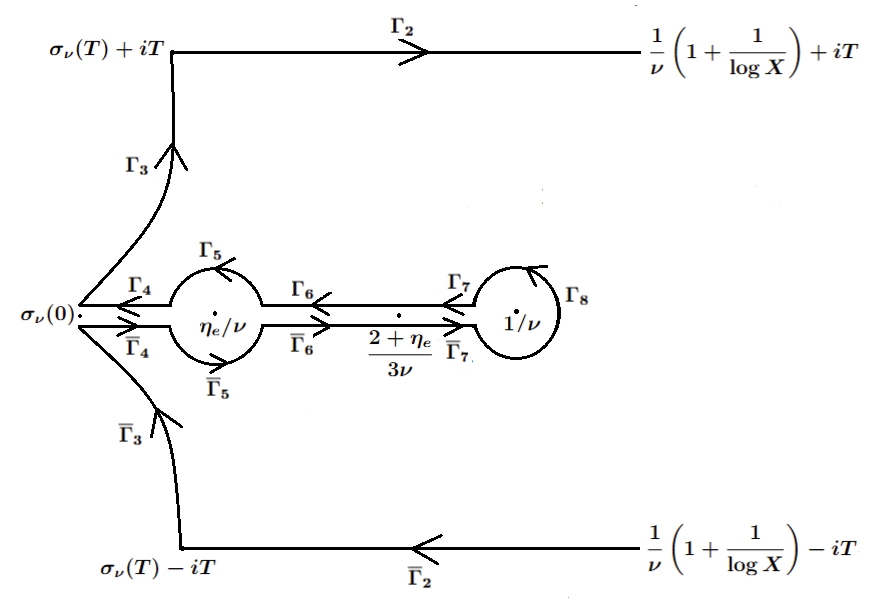}
\caption{Contour $\boldsymbol{\Gamma_0}$ for Theorem  \ref{thm:LFuncLSDVariant1Gen}(1), i.e.~when $\displaystyle{\boldsymbol{1-\frac{c_0}{10\lambdaqlogq} < \eta_e < 1-\frac{3\nu}{\log x}}}$.}
\label{Fig:LFuncLSDOuterContourVar1.1}
\end{figure}  
Since we are in the first case of \eqref{eq:Hfunc}, it follows that the function $\Fsnu G(s) X^s/s$ is analytic in a region containing the one enclosed by $\Gamma_0$ and the vertical segment $[\nu^{-1}(1+1/\log X)-iT, \nu^{-1}(1+1/\log X)+iT]$. 
As such, Cauchy's integral theorem yields, from \eqref{eq:PerronApp2},
\begin{equation}\label{eq:Var1.1PostContourShift}
\sum_{n \le X} \, a_n ~=~ \frac1{2\pi i}\int_{\Gamma_0} \, \frac{\Fsnu G(s) X^s}s \, \mathrm ds ~+~ O\left(\frac{\Omega(1) \, X^{1/\nu}}{T\log X}~+~ \kappa  \cdot \frac{X^{1+1/\nu} \log X}{Th}\right).
\end{equation}
We will now show that the contribution of all parts of $\Gamma_0$ except $\Gamma_7 + \overline\Gamma_7 + \Gamma_8$ are negligible. We will be repeatedly using the hypotheses $\boldsymbol{q \ge e^{4+5/3\nu}},$ $ \boldsymbol{\eta_e \ge 1-c_0/10\lambdaqlogq},$ and 
$\boldsymbol{\max\{|\alphachiZ|, |\alphachie|\} \le K_0}$. Our {implied constants are allowed to depend on} $\boldsymbol{c_0, \nu}$ \textbf{and} $\boldsymbol{K_0}$.  

\textbf{Contribution of $\sum_{j=2}^3 (\Gamma_j + \overline\Gamma_j)$.} Any $s$ on these four contours satisfies the conditions of Proposition \ref{prop:HsFsnuBounds}(2): The condition on $\sigma$ follows from the definition of $\sigmanu$ and the hypothesis ``${c_0\lambda_q\Elq(T) \le \log x}$" in Theorem \ref{thm:LFuncLSDVariant1Gen}(1). The other condition on $s$ in Proposition \ref{prop:HsFsnuBounds}(2) is clear if $|t| \ge 1/\nu$, whereas if $|t|< 1/\nu$, then $s \in \Gamma_3 + \overline\Gamma_3$, so $|s-1/\nu| \ge 1/\nu-\sigmanu(t) = c_0/4\nu\lambdaqElqt$ and
$$\left|s-\frac{\eta_e}\nu\right| \ge \frac{\eta_e}\nu - \sigmanu(t) \ge \frac{c_0}{4\nu \lambda_q\log(2q)} - \frac{c_0}{10\nu \lambdaqlogq} \ge \frac{c_0}{8\nu \lambdaqlogq} - \frac{c_0}{10\nu \lambdaqlogq} \ge \frac{c_0}{40\nu \lambdaqlogq}.$$
Hence Proposition \ref{prop:HsFsnuBounds}(2) yields  $\Fsnu \ll (5\lambda_q\Elq(T)/4)^{\lambda_q}$ uniformly for $s \in \Gamma_2 + \overline\Gamma_2 + \Gamma_3 + \overline\Gamma_3$. 
Moreover, $|G(s)| \le \Omega(|t|) \le \Omega(T)$ for all such $s$. Lastly,  $\int_{\Gamma_2 + \overline\Gamma_2} \, |X^s/s| \, \mathrm ds \, \le \, T^{-1}\, \int_{\sigmanu(T)}^{\frac1\nu+\frac1{\nu\log X}} \, X^\sigma \, \mathrm d\sigma \, \ll \, X^{1/\nu}/T\log X$, and $\int_{\Gamma_3 + \overline\Gamma_3} \, |X^s/s| \, \mathrm ds \, \le \, \int_0^T \, X^{\sigmanu(t)} \, \mathrm dt/(t+1) \, \le X^{\sigmanu(T)}(\log T)$, where we have noted that $|s| \gg |t|+1$ for all $s \in \Gamma_3 + \overline\Gamma_3$. Combining all the observations in this paragraph, we get
\begin{align}\allowdisplaybreaks\label{eq:Var1.1I2I3Final}
\sum_{j \in \{2, 3\}} \, \left|\int_{\Gamma_j + \overline\Gamma_j} \, \frac{\Fsnu G(s) X^s}s \, \mathrm ds\right| &\ll \left(\frac54 \lambda_q\Elq(T)\right)^{\lambda_q} \Omega(T) \left\{\frac{X^{1/\nu}}{T \log X} ~+~ X^{1/\nu - c_0/4\nu\lambda_q\Elq(T)} (\log T) \right\}.     
\end{align}
\textbf{Contribution of $\sum_{j=4}^6 \, (\Gamma_j + \overline\Gamma_j)$.} Since $r_e = (1-\eta_e)/6\nu \le c_0/60\nu\lambdaqlogq $, any $s$ on these six contours satisfies the conditions of Proposition \ref{prop:HsFsnuBounds}(1), as well as $|G(s)| \le \Omega(1)$. Any such $s$ also satisfies $\sigma \le (2+\eta_e)/3\nu$ and $|s| \ge \sigmanu(0) \gg 1$.  Hence by definition of $\Hfuncs$ and Proposition \ref{prop:HsFsnuBounds}(1),
\begin{multline}\label{eq:Var1.1I4I5I6Step1}\allowdisplaybreaks
\sum_{j=4}^6 \, \left|\int_{\Gamma_j + \overline\Gamma_j} \, \frac{\Fsnu G(s) X^s}s \, \mathrm ds\right|\\
\ll \, (2\lambda_q\log_q)^{\lambda_q+|\Ree(\alphachiZ)| + |\Ree(\alphachie)|} \, \Omega(1) \, X^{(2+\eta_e)/3\nu} \sum_{j=4}^6 \, \int_{\Gamma_j + \overline\Gamma_j} \, \left|s-\frac1\nu\right|^{-\Ree(\alphachiZ)} \cdot \left|s-\frac1\nu\right|^{\Ree(\alphachie)} \, |\mathrm ds|.
\end{multline} 
As Figure \ref{Fig:LFuncLSDOuterContourVar1.1} shows, any $s$ above satisfies $(1-\eta_e)/3\nu = 1/\nu - (2+\eta_e)/3\nu \le |s-1/\nu| \le 1/\nu - \sigmanu(0) = c_0/4\nu\lambdaqlogq < 1$ and $(1-\eta_e)/6\nu = r_e \le |s-\eta_e/\nu| \le 1/\nu-\sigmanu(0) < 1$.  
In particular, both 
$|s-1/\nu|$ and $|s-\eta_e/\nu|$ lie between $(1-\eta_e)/6\nu$ and $1$, so that $|s-1/\nu|^{-\Ree(\alphachiZ)} \cdot |s-\eta_e/\nu|^{\Ree(\alphachie)} \le |s-1/\nu|^{-|\Ree(\alphachiZ)|} \cdot |s-\eta_e/\nu|^{-|\Ree(\alphachie)|} ~ \ll_{\nu, K_0} ~~ (1-\eta_e)^{-|\Ree(\alphachiZ)| - |\Ree(\alphachie)|}$. Hence \eqref{eq:Var1.1I4I5I6Step1} yields 
\begin{equation}\label{eq:Var1.1I4I5I6Final}
    \sum_{j=4}^6 \, \left|\int_{\Gamma_j + \overline\Gamma_j} \, \frac{\Fsnu G(s) X^s}s \, \mathrm ds\right|\\
\ll \, \Omega(1) \cdot  \frac{(2\lambda_q\log_q)^{\lambda_q+|\Ree(\alphachiZ)| + |\Ree(\alphachie)|} \cdot X^{(2+\eta_e)/3\nu}}{(1-\eta_e)^{|\Ree(\alphachiZ)| + |\Ree(\alphachie)|}}.
\end{equation}
Now we extract the main term from the contribution of $\Gamma_7 + \overline\Gamma_7 + \Gamma_8$. Our method for this is partly inspired from works of Scourfield \cite{scourfield85} and Tenenbaum \cite[Chapter II.5]{tenenbaum15}. \textbf{We first claim that the disk} $\boldsymbol{\{s: |s-1/\nu| \le 2(1-\eta_e)/3\nu\}}$ \textbf{is contained} 
\textbf{in the region}  
$$\boldsymbol{\left\{s=\sigma+it: ~~ \sigma > \frac1\nu\left(1-\frac{c_0}{4\nu\lambdaqElqt}\right), ~~ |t| \le \frac{c_0}{15\nu \lambdaqlogq},  ~~ s \not\in \left(-\infty, \frac{\eta_e}\nu\right]\right\}}.$$ 
Indeed for any $s$ in this disk, we have $\sigma - \eta_e/\nu \ge \big(1/\nu-2(1-\eta_e)/3\nu\big) - \eta_e/\nu > 0$, and $|t| \le 2(1-\eta_e)/3\nu < c_0/15\nu \lambda_q \log q < 1/\nu$. (Recall $1-\eta_e < c_0/10\lambdaqlogq$.) The claim now follows from 
$$\sigma \ge \frac1\nu - \frac{2(1-\eta_e)}{3\nu} > \frac1\nu - \frac{c_0}{15\nu \lambda_q \log q} > \frac1\nu\left(1-\frac{c_0}{4\lambda_q\log(2q)}\right) > \frac1\nu\left(1-\frac{c_0}{4\lambda_q\log\big(q(|t\nu|+1)\big)}\right),$$
Two consequences of this claim are that any $s$ in the above disk 
satisfies the hypotheses of Proposition \ref{prop:HsFsnuBounds}(1), 
and that the function $\Hfuncs G(s) (s-\eta_e/\nu)^{\alphachie}$ is analytic on this disk. Hence 
\begin{align}\allowdisplaybreaks\label{eq:Var1.1HsGsPowerSeries}
&\Hfuncs G(s) \left(s-\frac{\eta_e}\nu\right)^{\alphachie} = \sum_{j=0}^\infty \mu_j \, \left(s-\frac1\nu\right)^j \text{ for all }s\text { satisfying }\left|s-\frac1\nu\right| \le \frac{2(1-\eta_e)}{3\nu}.
\end{align} 
where by Cauchy's integral formula and  Proposition \ref{prop:HsFsnuBounds}(1), we have 
\begin{equation}\label{eq:mujBound}
|\mu_j| = \frac1{2\pi i}\int_{\left|z-\frac1\nu\right| = \frac{2(1-\eta_e)}{3\nu}} ~~ \frac{\Hfunc(z) G(z) \left(z-\eta_e/\nu\right)^{\alphachie}}{(z-1/\nu)^{j+1}} \, \mathrm dz \, \ll \, \Omega(1) \cdot \frac{\TwoLambdaqLogqPower}{(1-\eta_e)^{|\Ree(\alphachie)|} \cdot \big(2(1-\eta_e)/3\nu\big)^j}   
\end{equation}
\text{uniformly in all }$j \ge 0$. (Here, we also observed that if $|z-1/\nu| = 2(1-\eta_e)/3\nu$, then $|z-\eta_e/\nu|$ lies in the interval $[(1-\eta_e)/3\nu, ~4(1-\eta_e)/3\nu]$ by the triangle inequality.) 

Now for $s \in \Gamma_7 + \overline\Gamma_7 + \Gamma_8$, we have $|s-1/\nu| \le 1/\nu - (2+\eta_e)/3\nu = (1-\eta_e)/3\nu$. 
Hence 
\eqref{eq:mujBound} yields  
\begin{equation}\label{eq:Var1.1HsGsPowerSeriesTail}
\sum_{j \ge N+1}\, \left|\mu_j \left(s-\frac1\nu\right)^j \right| \, \ll \, \Omega(1) \cdot \frac{\TwoLambdaqLogqPower}{(1-\eta_e)^{|\Ree(\alphachie)|}} \cdot \left(\frac{|s-1/\nu|}{2(1-\eta_e)/3\nu}\right)^{N+1}
\end{equation}
uniformly in all such $s$ and in $N \ge 0$;  here we noted that   
$\displaystyle{\sum_{j \ge N+1} \, \left(\frac{|s-1/\nu|}{2(1-\eta_e)/3\nu}\right)^{j-(N+1)}} \le \displaystyle{\sum_{j \ge N+1} \, \left(1/2\right)^{j-(N+1)} \le 2}$. 
Hence by definition of $\Hfuncs$, along with \eqref{eq:Var1.1HsGsPowerSeries} and \eqref{eq:Var1.1HsGsPowerSeriesTail}, we obtain 
\begin{multline}\label{eq:Var1.1PreHankel}\allowdisplaybreaks 
\frac1{2\pi i} \int_{\Gamma_7 + \overline\Gamma_7 + \Gamma_8} \, \frac{\Fsnu G(s) X^s}s \, \mathrm ds \,=\, \sum_{j=0}^N \, \frac{\mu_j}{2\pi i} \, \int_{\Gamma_7 + \overline\Gamma_7 + \Gamma_8} \, X^s \left(s-\frac1\nu\right)^{j-\alphachiZ} \, \mathrm ds \\ ~+~ O\left(\Omega(1) \cdot \frac{\TwoLambdaqLogqPower \cdot (3\nu/2)^N}{(1-\eta_e)^{N+1+|\Ree(\alphachie)|}} \cdot \int_{\Gamma_7 + \overline\Gamma_7 + \Gamma_8} \, X^\sigma \left|s-\frac1\nu\right|^{N+1-\Ree(\alphachiZ)} \, |\mathrm ds|\right).   
\end{multline}
Now setting $w \coloneqq (s-1/\nu)\log X$, we see that the total main term on the right hand of \eqref{eq:Var1.1PreHankel} is 
\begin{multline}\label{eq:Var1.1PostHankel}\allowdisplaybreaks
\sum_{j=0}^N \, \frac{\mu_j}{2\pi i} \, \int_{\Gamma_7 + \overline\Gamma_7 + \Gamma_8} \, X^s \left(s-\frac1\nu\right)^{j-\alphachiZ} \, \mathrm ds \, = \, \sum_{j=0}^N \, \frac{\mu_j \, X^{1/\nu}}{(\log X)^{j+1-\alphachiZ}} \cdot \frac1{2\pi i} \, \int_{\mathcal W} \, e^w \, w^{j - \alphachiZ} \, \mathrm dw \,\\ 
= \, \frac{X^{1/\nu}}{(\log X)^{1-\alphachiZ}} \sum_{j=0}^N \, \frac{\mu_j (\log X)^{-j}}{\Gamma(\alphachiZ-j)} \, + \, O\left(\frac{X^{1/\nu + (\eta_e-1)/6\nu}}{(\log X)^{1-\Ree(\alphachiZ)}} \sum_{j=0}^N \, |\mu_j| \, \Gamma(j+1+|\alphachiZ|) \, \left(\frac{47}{\log X}\right)^j \right),     
\end{multline}
where $\mathcal W$ is the truncated Hankel contour in Figure \ref{Fig:Var1.1HankelContour}, and we have used \cite[Corollary II.0.18]{tenenbaum15}. (The hypothesis $\eta_e < 1-3\nu/\log x$ in Theorem \ref{thm:LFuncLSDVariant1Gen}(1)  guarantees that \cite[Corollary II.0.18]{tenenbaum15} applies.) 
\begin{figure}[h!]
\includegraphics[height=3cm]{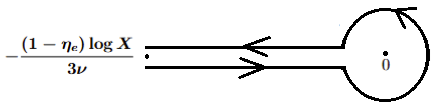}
\caption{The truncated Hankel contour $\mathcal W$ after the substitution $\displaystyle{w=\left(s-\frac1\nu\right)\log X}$.}~\\ 
\label{Fig:Var1.1HankelContour}
\end{figure} 
By \eqref{eq:mujBound}, the entire $O$--term in \eqref{eq:Var1.1PostHankel} is 
\begin{align*}\allowdisplaybreaks
&\ll \frac{\TwoLambdaqLogqPower \, \Omega(1) \, X^{1/\nu + (\eta_e-1)/6\nu}}{(1-\eta_e)^{|\Ree(\alphachie)|} \cdot (\log X)^{1-\Ree(\alphachiZ)}} \, \sum_{j=0}^N \, \frac{\Gamma(j+1+|\alphachiZ|)}{\big(2(1-\eta_e)\log X/141\nu\big)^j}
\end{align*}
Since $\Gamma(N+1+|\alphachiZ|) \, = \,\Gamma(j+1+|\alphachiZ|) \cdot \prod_{i=j+1}^N \, (i+|\alphachiZ|) \, \ge \, \Gamma(j+1+|\alphachiZ|) \cdot (N-j)!$, we get 
$$\sum_{j=0}^N \, \frac{\Gamma(j+1+|\alphachiZ|)}{\big(2(1-\eta_e)\log X/141\nu\big)^j} \, = \, \frac{(141\nu/2)^N \, \Gamma(N+1+|\alphachiZ|)}{(1-\eta_e)^N \cdot (\log x)^N} \, \sum_{j=0}^N \, \frac{\big(2(1-\eta_e)\log X/141\nu\big)^{N-j}}{(N-j)!},  
$$
and the last sum above is at most $X^{2(1-\eta_e)/141\nu}$. 
Inserting this into \eqref{eq:Var1.1PostHankel}, we obtain 
\begin{multline}\label{eq:Var1.1PostHankelMainTerm}\allowdisplaybreaks
\sum_{j=0}^N \, \frac{\mu_j}{2\pi i} \, \int_{\Gamma_7 + \overline\Gamma_7 + \Gamma_8} \, X^s \left(s-\frac1\nu\right)^{j-\alphachiZ} \, \mathrm ds \, = \, \frac{X^{1/\nu}}{(\log X)^{1-\alphachiZ}} \sum_{j=0}^N \, \frac{\mu_j (\log X)^{-j}}{\Gamma(\alphachiZ-j)} \,\\ + \, O\left(\Omega(1) \, \TwoLambdaqLogqPower \cdot \frac{\left({141\nu}/2\right)^N \, \Gamma(N+1+K_0) \cdot X^{1/\nu + {43(\eta_e-1)}/{282\nu}}}{(1-\eta_e)^{N+|\Ree(\alphachie)|} \cdot (\log X)^{N+1-\Ree(\alphachiZ)}} \right),     
\end{multline}
Next, we bound the integral in the $O$--term of \eqref{eq:Var1.1PreHankel}. Note that for all $s \in \Gamma_7 + \overline\Gamma_7 + \Gamma_8$, we have $r_1 \le |1/\nu - s| < 1/\nu - (2+\eta_e)/3\nu = (1-\eta_e)/3\nu < c_0/30 \lambdaqlogq < 1$. As such $|1/\nu-s|^{-\Ree(\alphachiZ)} \le |1/\nu-s|^{-|\Ree(\alphachiZ)|} < (r_1)^{-|\Ree(\alphachiZ)|} \ll (\log x)^{|\Ree(\alphachiZ)|}$ via \eqref{eq:sigmanurerDefs}. Parametrizing the circle $\Gamma_8$ as $s=1/\nu+r_1 e^{i\theta}$, $-\pi < \theta \le \pi$, we thus find that the integral in the $O$--term of \eqref{eq:Var1.1PreHankel} is 
$$\ll \, (\log x)^{|\Ree(\alphachiZ)|} \left\{\int_{\frac{2+\eta_e}{3\nu}}^{\frac1\nu-r_1} \, X^\sigma \left(\frac1\nu-\sigma\right)^{N+1} \, \mathrm d\sigma \, + \, X^{1/\nu+r_1} \, r_1^{N+2}\right\} \, \ll \, \frac{(N+1)! \, X^{1/\nu}}{(\log X)^{N+2-|\Ree(\alphachiZ)|}}.$$
To get the last bound, we have again used \eqref{eq:sigmanurerDefs}, and we have made the change of variable $\sigma = 1/\nu-u/\log X$ in the last integral above to see that its value is at most $\Gamma(N+2) X^{N+1}/(\log X)^{N+2}$. Inserting the above bound into \eqref{eq:Var1.1PreHankel}, we find that the total $O$--term in \eqref{eq:Var1.1PreHankel} is 
\begin{equation*}
\ll \, \Omega(1) \cdot \frac{\TwoLambdaqLogqPower \cdot (3\nu/2)^N \, (N+1)! \, X^{1/\nu} }{(1-\eta_e)^{N+1+|\Ree(\alphachie)|} \cdot (\log X)^{N+2-|\Ree(\alphachiZ)|}}.
\end{equation*}
Inserting this observation along with \eqref{eq:Var1.1PostHankelMainTerm} into \eqref{eq:Var1.1PreHankel}, we now obtain 
\begin{multline}\allowdisplaybreaks\label{eq:Var1.1I7I8Final}
\frac1{2\pi i} \int_{\Gamma_7 + \overline\Gamma_7 + \Gamma_8} \, \frac{\Fsnu G(s) X^s}s \, \mathrm ds \, - \, \frac{X^{1/\nu}}{(\log X)^{1-\alphachiZ}} \sum_{j=0}^N \, \frac{\mu_j (\log X)^{-j}}{\Gamma(\alphachiZ-j)}\\
\,\ll \,   \Omega(1) \cdot \frac{\TwoLambdaqLogqPower \cdot N! \, (71\nu)^N \,  X^{1/\nu} }{(1-\eta_e)^{N+1+|\Ree(\alphachie)|} \cdot (\log X)^{N+2-|\Ree(\alphachiZ)|}};   
\end{multline}
to absorb the $O$--term in \eqref{eq:Var1.1PostHankelMainTerm} into the right side above, we noted that  
$\Gamma(N+1+K_0)/(N+1)! \, \le \Gamma(1+K_0) \cdot \prod_{i=1}^{N+1} \,  (1+K_0/i) \, \ll_{K_0} \, \exp(K_0 \sum_{i=1}^{N+1} \, 1/i) \, \ll \, (N+1)^{K_0} \ll (N+1)^{-1} \cdot (142/141)^N$.

Combining \eqref{eq:Var1.1I2I3Final}, \eqref{eq:Var1.1I4I5I6Final} and \eqref{eq:Var1.1I7I8Final}, we deduce that Theorem \ref{thm:LFuncLSDVariant1Gen} holds with 
$X \in (x, x+h]$ in place of $x$. Now  $|\sum_{n \le X} a_n - \sum_{n \le x} a_n| \, \le \, \sum_{x<n \le x+h} \, |a_n|$ is absorbed in the right of \eqref{eq:LFuncLSDVar1SiegelZeroGen}. Hence, to complete the proof of Theorem \ref{thm:LFuncLSDVariant1Gen}(1), it only  remains to show that the difference between $\sum_{j=0}^N \, \mu_j \, X^{1/\nu} (\log X)^{-j-1+\alphachiZ}/\Gamma(\alphachiZ-j)$ and its ``$x$--analogue" is absorbed in the right of \eqref{eq:LFuncLSDVar1SiegelZeroGen}. 

To this end, note that since $x<X \le x(1+h/x)$ and $h<x/2$, we have  
$\boldsymbol{X^{1/\nu} = x^{1/\nu}(1+O(h/x))}$ 
and $\log x < \log X  
\le \log x + h/x$, so that $\boldsymbol{|\log X - \log x| \le h/x}$. This in turn shows that the difference $|(\log X)^{-j-1+\alphachiZ} - (\log x)^{-j-1+\alphachiZ}|$ is at most $|j+1-\alphachiZ| \cdot \int_{\log x}^{\log X} \, t^{-j-2+\Ree(\alphachiZ)} \, \mathrm dt \, \ll \, (jh/x) \cdot (\log x)^{-j-2+\Ree(\alphachiZ)}$, 
uniformly in \textit{all} $j \ge 0$. (Here the last ``$\ll$" bound is tautological for $j+2 > \Ree(\alphachiZ)$, while for $j+2 \le \Ree(\alphachiZ) \le K_0$, it follows from the fact that $\log X = (\log x) \, (1+h/x\log x)$.)   
As such, $\boldsymbol{(\log X)^{-j-1+\alphachiZ} \, = \, (\log x)^{-j-1+\alphachiZ} \cdot \big(1+O(jh/x\log x)\big)}$. 

We will also need to upper bound $1/|\Gamma(\alphachiZ-j)|$ uniformly in all $j \in \{0, \dots, N\}$. If $j \le 2K_0+1$, then the fact that $1/\Gamma$ is entire yields $1/|\Gamma(\alphachiZ-j)|\, \ll_{K_0} \, 1$. On the other hand, if $j>2K_0+1$, then from $|\Gamma(\lceil K_0 \rceil + 1 - \alphachiZ)| \, \ll_{K_0} \, 1$ 
and the reflection identity, we find that 
\begin{align*}\allowdisplaybreaks
\frac1{|\Gamma(\alphachiZ-j)|} &\ll |\sin(\pi\alphachiZ)| \cdot |\Gamma(j+1-\alphachiZ)| \, \ll_{K_0} \, |\Gamma(\lceil K_0 \rceil + 1 - \alphachiZ)| \cdot \prod_{i=\lceil K_0 \rceil + 1}^j \, \left|i - \alphachiZ\right|\\
&\ll_{K_0} \, 
\prod_{i=1}^j \, \left(i + K_0\right) \le N! \cdot \prod_{i=1}^{N+1} \, \left(1+\frac{K_0}i\right) \, \le \, N! \cdot \exp \left({K_0}\sum_{i=1}^{N+1} \, \frac1i\right) \ll N! \cdot (N+1)^{K_0}
\end{align*}
Hence, $\boldsymbol{1/|\Gamma(\alphachiZ-j)| \ll N! \, (N+1)^{K_0}}$ uniformly in all $j \in \{0, \dots, N\}$. Collecting all the observations in bold in this paragraph and the last, and using \eqref{eq:mujBound}, we find that 
\begin{multline*}\allowdisplaybreaks
\left|\sum_{j=0}^N \, \frac{\mu_j}{\Gamma(\alphachiZ-j)} \left\{ \frac{X^{1/\nu}}{(\log X)^{j+1-\Ree(\alphachiZ)}} -  \frac{x^{1/\nu}}{(\log x)^{j+1-\Ree(\alphachiZ)}}\right\}\right| \, 
\\\ll \, \frac hx \cdot \Omega(1) \cdot \frac{\TwoLambdaqLogqPower \cdot N! \, (N+1)^{K_0} \cdot x^{1/\nu}}{(1-\eta_e)^{|\Ree(\alphachie)|} \cdot (\log x)^{1-\Ree(\alphachiZ)}} \sum_{j=0}^N \, j{\left(\frac{2(1-\eta_e)\log x}{3\nu}\right)^{-j}},   
\end{multline*}
which is absorbed in the right of \eqref{eq:LFuncLSDVar1SiegelZeroGen}, since by the hypothesis $(1-\eta_e)\log x > 3\nu$,   the sum on $j$ above is at most $\sum_{j \ge 0} \, j/2^j \ll 1$. This concludes the proof of subpart \textbf{(1)} of Theorem \ref{thm:LFuncLSDVariant1Gen}. \hfill \qedsymbol
\subsection{When $\boldsymbol{\eta_e}$ does not exist or $\boldsymbol{\eta_e \le 1-c_0/10\lambdaqlogq}$: Proof of Theorem \ref{thm:LFuncLSDVariant1Gen}(2)} \label{subsec:Variant1NoSiegelZeroProof}~\\  
We are in cases 2 and 3 of \eqref{eq:Hfunc}. 
We just mention the main changes from the above arguments. First, we redefine $\sigmanu(t) \coloneqq \nu^{-1}\big(1-c_0/\boldsymbol{8} \lambdaqElqt\big)$ for case 2, and ${\sigmanu(t) \coloneqq \nu^{-1}\big(1-c_0/\boldsymbol{80} \lambdaqElqt\big)}$ for {case 3}. Second, we redefine $\Gamma_0$ by replacing all the contours $\sum_{j=4}^7 \, (\Gamma_j + \overline\Gamma_j)$ by $\Gamma_4' + \overline{\Gamma_4'}$, where $\Gamma_4'$ is  the horizontal segment traversed from $1/\nu-r_1$ to $\sigmanu(0)$ above the branch cut. Here $\Gamma_3 + \overline\Gamma_3$ have been  automatically redefined using the respective $\sigmanu$'s, and $r_1$ is still as defined in \eqref{eq:sigmanurerDefs}. 

In the rest of the subsection, we continue only with case 3 of \eqref{eq:Hfunc} (namely, that $\eta_e \le 1-c_0/10\lambdaqElqt)$), since case 2 will be entirely analogous and simpler. Note that by the ``case 3" assertion of Proposition \ref{prop:HsFsnuBounds}, the analogue of \eqref{eq:Var1.1I2I3Final} continues to hold, only with the ``$X^{-c_0/4\nu \lambda_q \Elq(T)}$" term replaced 
``$X^{-c_0/80\nu \lambda_q \Elq(T)}$" 
Our main terms come from $\Gamma_4' + \overline{\Gamma_4'} + \Gamma_8$.  
The 
analogue of \eqref{eq:Var1.1HsGsPowerSeries} is that $\Hfuncs G(s) = \sum_{j=0}^\infty \, \mu_j (s-1/\nu)^j$ for all $s$ satisfying $|s-1/\nu| \le c_0/40 \nu \lambda_q \log q$. As such, the respective analogues of \eqref{eq:mujBound} and  \eqref{eq:Var1.1HsGsPowerSeriesTail}  continue to hold, with all instances of ``$\eta_e$" and ``$\alphachie$" removed, and with all instances of ``$2(1-\eta_e)/3\nu$" replaced by ``$c_0/40\nu\lambda_q \log q$". Finally, the analogues of  \eqref{eq:Var1.1PreHankel},   \eqref{eq:Var1.1PostHankelMainTerm} and \eqref{eq:Var1.1I7I8Final} hold with ``$\Gamma_4'+\overline{\Gamma_4'}$" playing the role of ``$\Gamma_7 + \overline\Gamma_7$";   
the error term in the analogue of \eqref{eq:Var1.1I7I8Final} is $\ll \Omega(1) \cdot (2\lambda_q \log q)^{\lambda_q+|\Ree(\alphachiZ)|} \, N!\, (2000 \, \nu \, c_0^{-1} \cdot \lambdaqlogq)^{N+1} \, X^{1/\nu}/(\log X)^{N+2-|\Ree(\alphachiZ)|}$. \hfill \qedsymbol
\begin{figure}[h!]
\includegraphics[height=11cm]{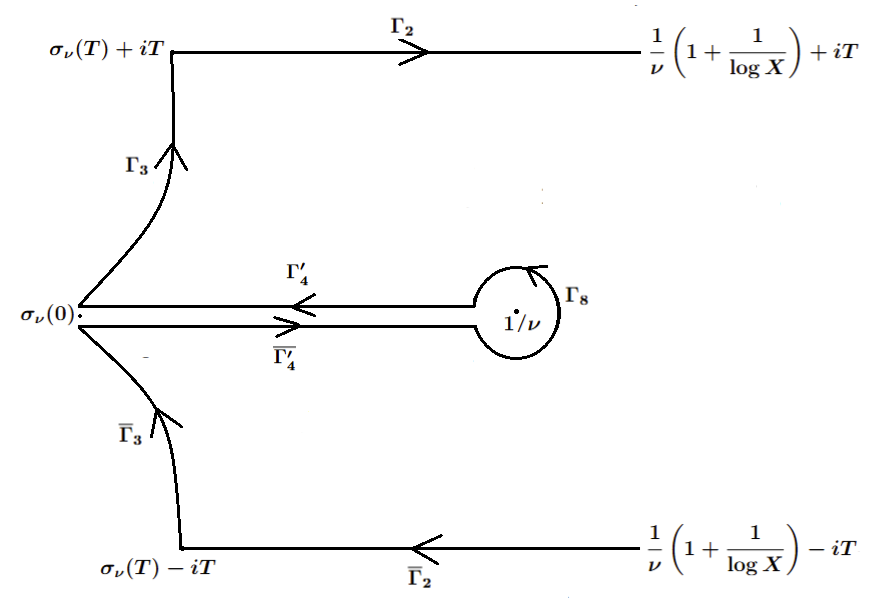}
\caption{Contour $\boldsymbol{\Gamma_0}$  when $\boldsymbol{\displaystyle{\eta_e \le 1-\frac{c_0}{10\lambdaqlogq}}}$. Here $\displaystyle{\boldsymbol{\sigmanu(t) = \frac1\nu\left(1-\frac{c_0}{80\lambdaqElqt}\right)}}$.} 
\label{Fig:LFuncLSDOuterContourVar1.2}
\end{figure} 
\subsection{When the stronger growth condition \eqref{eq:Variant1GrowthConditionSpecial} is available} To complete the proof of Theorem \ref{thm:LFuncLSDVariant1Gen}, it thus only remains to show its very last assertion (which would then also establish Corollary \ref{cor:LFuncLSDVariant1LogPowerInterval}). To this end, we just need to show that 
condition \eqref{eq:Variant1GrowthConditionSpecial}, (assumed for some $A \ge 1$ and $\kappa_A \ge 2$)    
implies   
condition \eqref{eq:Variant1GrowthConditionNew} with $\kappa \coloneqq {\kappa_A \cdot 2^{A+1/\nu}}$.  

Indeed, given any $x \ge 2$, consider the sequence  $\{y_m\}_{m=0}^\infty$ defined by $y_0 \coloneqq x$ and $y_{m+1} \coloneqq y_m + y_m/(\log y_m)^A$. Let $M \coloneqq M(x)$ be the unique index satisfying $y_M \le x <  y_{M+1}$. Since $\{y_m\}_m$ is clearly increasing, we have $x < y_m \le 2x$ for all $m \in \{0, \dots, M\}$, so that by \eqref{eq:Variant1GrowthConditionSpecial} yields
$$\sum_{x < n \le 2x} \, |a_n| \, \le \, \sum_{0 \le m \le M} \, \sum_{y_m < x \le y_{m+1}} \, |a_n| \, \le \kappa_A \sum_{0 \le m \le M} \, \frac{y_m^{1/\nu}}{(\log y_m)^A} \le \kappa_A \cdot \frac{M \, (2x)^{1/\nu}}{(\log x)^A}.$$
But again the fact that $y_m \in (x, 2x]$ 
and the recurrence defining $y_m$ yield $2x \ge y_M \ge y_{M-1}+ x/\big(\log (2x) \big)^A\ge y_{M-2}+ 2x/\big(\log (2x) \big)^A \ge \dots \ge x+ Mx/\big(\log (2x) \big)^A$, leading to $M \le \big(\log (2x) \big)^A$. Inserting this into the above display yields  \eqref{eq:Variant1GrowthConditionNew} with the desired value of $\kappa$.  
\hfill \qedsymbol
\section{The LSD method for $L$--functions via a similarly--behaving bounding sequence:\\ Proof of Theorem \ref{thm:LFuncLSDVariant2} }\label{sec:LfunctionsLSDVariant2Proof}
In this section, we are in the setting of Theorem \ref{thm:LFuncLSDVariant2}, so we only assume that $\{a_n\}_n$ has property $\mathcal P(\nu, \{\alpha_\chi\}_\chi; c_0, \Omega)$ with $\Omega(t) = \CalM (1+t)^{1-\delta_0}$, and that there exist $\{b_n\}_n$ satisfying $|a_n| \le b_n$ for all $n$, 
such that $\{b_n\}_n$ has property $\mathcal P(\nu, \{\beta_\chi\}_\chi; c_0, \Omega)$ for some $\{\beta_\chi\}_\chi \subset \C$. In the entire section, $\max\{|\alphachiZ|, |\alphachie|, |\beta_{\chi_0}|, |\beta_{\chi_e}|\} \le K_0$, and the implied constants depend only on $c_0, \nu, \delta_0$ and $K_0$. The parameter $\boldsymbol{\Lambda_q \coloneqq 1 + \max\limits_{\substack{a \bmod q}}~\max\left\{\left|\sum_\chi \,  \alpha_\chi \cdot \chi(a)\right|, \left|\sum_\chi \,\beta_\chi \cdot \chi(a)\right|\right\}}$ plays the role of ``$\lambda_q$" from section  \ref{sec:AnalyticMachinery}: As such, 
the analogue of Proposition \ref{prop:HsFsnuBounds} (which we will use without further reference) holds, with $\lambda_q$ replaced by $\Lambda_q$, and with ``cases 1-3" from \eqref{eq:Hfunc} redefined accordingly. 

We will only highlight the main changes required from the arguments in the previous section. The first main change is that we use a different version of Perron's formula \cite[Theorem II.2.5]{tenenbaum15} to write 
\begin{equation}\label{eq:Variant2.1Perron}
\int_0^x \, A(t) \, \mathrm dt \, = \, \frac1{2\pi i}  \int_{\frac1\nu\left(1+\frac1{\log x}\right)-i\infty}^{\frac1\nu\left(1+\frac1{\log x}\right)+i\infty} \, \frac{\Fsnu G(s) x^{s+1}}{s(s+1)} \, \ds, 
\end{equation}
where $\boldsymbol{A(x) \coloneqq \sum_{n \le x} \, a_n}$. Note that there is no error term here, and ``$x$" itself will play the role of the ``$X$" in the previous section. 
Next, all our contours $\Gamma_0$ from the previous section are redefined, with the following additional specifications (everything else being as before). 
\begin{itemize}
\item $T \ge 1$ is a parameter satisfying $\boldsymbol{c_0 \Lambda_q \Elq(T) \le 100 \log (x/2)}$.
\item $\sigmanu(t)$ is as in section \ref{sec:LfunctionsLSDVariant1Proof}, defined according to the (analogues of the) three cases of \eqref{eq:Hfunc}.    
\item ``$X$" and ``$\lambda_q$" are replaced by ``$x$" and ``$\Lambda_q$", respectively (in all three cases). 
\item We add the infinite vertical line  $\boldsymbol{\Gamma_1 \coloneqq [\nu^{-1}(1+\log x)+iT, ~ \nu^{-1}(1+\log x)+i\infty)}$ and its reflection about the real axis, both of these traversed upwards (in all three cases). 
\item $r_1$, the radius of the circle $\Gamma_8$, is just taken to be any positive parameter satisfying $r_1 \le (4\max\{1, \nu\}\log x)^{-1}$. We do not choose $r_1$ until the very end (to avoid technical issues).  
\end{itemize} 
In all cases, our integrand on the right of \eqref{eq:Variant2.1Perron} is holomorphic on the region enclosed by (the respective) $\Gamma_0$ and the infinite vertical line $(\nu^{-1}(1+\log x)+iT, ~ \nu^{-1}(1+\log x)+i\infty)$, so that 
\begin{equation}\label{eq:Variant2.1ContourShift}
\int_0^x \, A(t) \, \mathrm dt \, = \, \frac1{2\pi i}  \int_{\Gamma_0} \, \frac{\Fsnu G(s) x^{s+1}}{s(s+1)} \, \ds. 
\end{equation}
The third nontrivial change from the previous section is that we will use that 
\begin{equation}\label{eq:Lq(t)Lambdaq}
\Elqt^{\Lambda_q} \, \le \left(2\delta_0^{-1} \, \Lambdaqlogq\right)^{\Lambda_q} \, (1+|t\nu|)^{\delta_0/2} \, \ll_{\nu, \delta_0} \,  \left(2\delta_0^{-1} \, \Lambdaqlogq\right)^{\Lambda_q} \, (1+|t|)^{\delta_0/2}   
\end{equation}
uniformly in \textbf{all real} $\boldsymbol{t}$. Here the first inequality can be seen by calculus:  
If $\Lambda_q \le \delta_0(\log q)/2$,  the function $\Elqt^{\Lambda_q}/(1+t\nu)^{\delta_0/2}$ is strictly decreasing on $(0, \infty)$;  in the other case, 
the function attains its global maximum (over all positive reals) at the unique $t_{\max}>0$ satisfying $\Elq(t_{\max}) = 2\Lambda_q/\delta_0$. 

We henceforth restrict to (the redefined) case 1 of \eqref{eq:Hfunc} and establish subpart \textbf{(1)} of Theorem \ref{thm:LFuncLSDVariant2}; the changes required for the other cases are entirely analogous to those described in subsection \cref{subsec:Variant1NoSiegelZeroProof}. Hence, we are assuming that $1-c_0/10\Lambdaqlogq < \eta_e < 1-3\nu/\log x$, and we have $\sigmanu(t) = \nu^{-1}(1-c_0/4\LambdaqElqt)$. The analogue of Proposition \ref{prop:HsFsnuBounds}(2), combined with \eqref{eq:Lq(t)Lambdaq}, yields  $\Fsnu G(s) \ll \CalM \, (5\Lambda_q^2 \log q/2\delta_0)^{\Lambda_q} \cdot (1+t)^{1-\delta_0/2}$ uniformly in $s \in \sum_{j=1}^3 \, (\Gamma_j + \GammajBar)$. Hence the analogue of \eqref{eq:Var1.1I2I3Final} is 
\begin{equation}\label{eq:Var2I1I2I3Final}
\sum_{j=1}^3 \, \left|\int_{\Gamma_j +\GammajBar} \, \frac{\Fsnu G(s) x^{s+1}}{s(s+1)} \, \ds\right| \, \ll \, \left(\frac{5\Lambda_q^2 \log q}{2\delta_0}\right)^{\Lambda_q} \, \CalM \, x \left\{\frac{x^{1/\nu}}{T^{\delta_0/2}} + x^{\sigmanu(T)}\right\}.     
\end{equation}
Likewise, following the method leading to \eqref{eq:Var1.1I4I5I6Final}, we find that its analogue here is 
\begin{equation}\label{eq:Var2I4I5I6Final}
\sum_{j=4}^6 \, \left|\int_{\Gamma_j +\GammajBar} \, \frac{\Fsnu G(s) x^{s+1}}{s(s+1)} \, \ds\right| \, \ll \, \left(2\Lambdaqlogq\right)^{\Lambda_q+2K_0} \, (1-\eta_e)^{-2K_0} \, \CalM \, x^{1 \, + \, (2+\eta_e)/3\nu}.      
\end{equation}
From this point on, we closely follow part of the argument given for \cite[Theorem II.5.2]{tenenbaum15}. Inserting \eqref{eq:Var2I1I2I3Final} and \eqref{eq:Var2I4I5I6Final} into \eqref{eq:Variant2.1ContourShift}, we obtain   $\int_0^x \, A(t) \, \mathrm dt \, = \, \Phi(x) + O\big(\err \big)$, where $$\boldsymbol{\Phi(x) \coloneqq} \boldsymbol{\frac1{2 \pi i} \int_{\Gamma_7 + \overline\Gamma_7 + \Gamma_8} \, \frac{\Fsnu G(s) x^{s+1}}{s(s+1)} \, \ds},$$ and where  $\boldsymbol{\err}$ is the sum of the two  expressions on the right hand sides of \eqref{eq:Var2I1I2I3Final} and \eqref{eq:Var2I4I5I6Final}. With $\boldsymbol{u \in [-x/2, x/2]}$ a parameter to be specified later, it thus follows that 
\begin{equation}\label{eq:AxFiniteDifference}
\int_x^{x+u} \, A(t) \, \mathrm dt \, = \, \Phi(x+u) - \Phi(x) + O(\err) \, = \, u \, \Phi'(x) + u^2\int_0^1 \, (1-t)\Phi''(x+tu) \, \mathrm + O(\err).    
\end{equation}
Differentiating under the integral sign, we have $\Phi'(x) = (2\pi i)^{-1} \int_{\Gamma_7 + \overline\Gamma_7 + \Gamma_8} \, {\Fsnu G(s) x^{s}}/{s} \, \ds$ and $\Phi'(x) = (2\pi i)^{-1} \int_{\Gamma_7 + \overline\Gamma_7 + \Gamma_8} \, {\Fsnu G(s) x^{s-1}} \, \ds$. Now uniformly in $y \in [x/2, 2x]$, we see that 
\begin{align*}\allowdisplaybreaks
|\Phi''(y)| \, &\ll \, \int_{\Gamma_7 + \overline\Gamma_7 + \Gamma_8} \, |\Hfuncs G(s)| \, y^{\sigma-1} \, \left|s-\frac1\nu\right|^{-\Ree(\alphachiZ)} \cdot \left|s-\frac{\eta_e}\nu\right|^{\Ree(\alphachie)} \, |\ds| \\
&\ll \, (2\Lambdaqlogq)^{\Lambda_q + 2K_0} \, \big(r_1(1-\eta_e)\big)^{-K_0}\cdot \CalM x^{1/\nu-1}, 
\end{align*}
where we have used the definition of $\Hfuncs$ and the analogue of Proposition \ref{prop:HsFsnuBounds}(1), along with the facts that $r_1 \le |s-1/\nu| < 1/\nu - \sigmanu(0) < 1$, that $2(1-\eta_e)/3\nu \, \le \, |s-\eta_e/\nu| \, \le \, 7(1-\eta_e)/6\nu$, that $\max\{|\Ree(\alphachiZ)|, |\Ree(\alphachie)|\} \le K_0$, and that $y^{\sigma - 1} \, \asymp \, x^{\sigma -1} \, \le \, x^{1/\nu+r_1-1} \, \ll \, x^{1/\nu-1}$ as $r_1 < (4 \max\{1, \nu\} \log x)^{-1}$. Inserting the bound in the above display into \eqref{eq:AxFiniteDifference}, we obtain 
\begin{equation}\label{eq:AxAverage}
\frac1u\int_x^{x+u} \, A(t) \, \mathrm dt \, = \, \Phi'(x) \, + \, O\left(\frac{\err}u \, + \, u \cdot (2\Lambdaqlogq)^{\Lambda_q + 2K_0} \, \big(r_1(1-\eta_e)\big)^{-K_0}\cdot \CalM x^{1/\nu-1} \right).    
\end{equation}
Now we use the properties of $\{b_n\}_n$. Since $|a_n| \le b_n$ for all $n$, we find that for all $u>0$,
\begin{align}\allowdisplaybreaks\label{eq:AxAverageMinusAx}
\nonumber \left| \frac1u\int_x^{x+u} \, A(t) \, \mathrm dt \, - \, A(x)\right| &\le  \frac1u\int_x^{x+u} \, \left|\sum_{x < n \le t} \, a_n\right| \, \mathrm dt \le  \frac1u\int_x^{x+u} \, \left(\sum_{x < n \le t} \, b_n\right) \, \mathrm dt\\
 &= \frac1u\int_x^{x+u} \, \left(B(t) - B(x)\right) \, \mathrm dt \le \left|\frac1u\int_x^{x+u} \, B(t) \, \mathrm dt \, -\, \frac1u\int_{x-u}^x \, B(t) \, \mathrm dt \right|,
\end{align}
where we have noted that $B(t) \coloneqq \sum_{n \le t} \, b_n$ is an increasing function of $t$. Now since $\{b_n\}_n$ has property $\mathcal P(\nu, \{\beta_\chi\}_\chi; c_0, \Omega)$ and $\max\{|\beta_{\chi_0}|, \, |\beta_{\chi_e}|\} \le K_0$, we see that \textbf{all} the arguments given for \eqref{eq:AxAverage} go through for $\{b_n\}_n$, so that the analogue of \eqref{eq:AxAverage} holds for $\{b_n\}_n$ as well. Hence uniformly in $u \in (0, x/2]$, the last difference in \eqref{eq:AxAverageMinusAx} is absorbed into the $O$--term in \eqref{eq:AxAverage}. Thus 
from \eqref{eq:AxAverage}, 
\begin{equation*}
A(x) = \frac1{2 \pi i} \int_{\Gamma_7 + \overline\Gamma_7 + \Gamma_8} \, \frac{\Fsnu G(s) x^s}s \, \ds \, + \, O\left(\frac{\err}u  + u \cdot (2\Lambdaqlogq)^{\Lambda_q + 2K_0} \big(r_1(1-\eta_e)\big)^{-K_0} \, \CalM x^{1/\nu-1} \right),
\end{equation*}
where we have also used the expression for $\Phi'(x)$ given after \eqref{eq:AxFiniteDifference}. Finally, the analogue of \eqref{eq:Var1.1I7I8Final} (with ``$X$" and ``$\Lambda_q$" replaced by ``$x$" and ``$\lambda_q$" respectively) estimates the last integral above, showing that the left hand side of \eqref{eq:LFuncLSDVar2NoSiegelZero} is bounded by the expression  
\begin{multline*}
\frac{(2 \Lambda_q \log q)^{\Lambda_q+|\Ree(\alphachiZ)|+|\Ree(\alphachie)|} \, N! \, (71 \nu)^N \, \CalM \, x^{1/\nu}}{(1-\eta_e)^{N+1+|\Ree(\alphachie)|} \cdot (\log x)^{N+2-|\Ree(\alphachiZ)|}} \, + \, u \, (2\Lambdaqlogq)^{\Lambda_q + 2K_0} \big(r_1(1-\eta_e)\big)^{-K_0} \, \CalM x^{1/\nu-1}\\ + 
u^{-1} \, \left(\frac{5\Lambda_q^2 \log q}{2\delta_0}\right)^{\Lambda_q} \, \CalM \, x \left\{\frac{x^{1/\nu}}{T^{\delta_0/2}} + x^{\sigmanu(T)}\right\}\, + \, u^{-1} \, \left(2\Lambdaqlogq\right)^{\Lambda_q+2K_0} \, (1-\eta_e)^{-2K_0} \, \CalM \, x^{1 \, + \, (2+\eta_e)/3\nu},
\end{multline*}
uniformly in all $u \in (0, x/2]$, and in $T \ge 1$ satisfying $c_0\lambda_q\Elq(T) \le 100\log(x/2)$. Theorem \ref{thm:LFuncLSDVariant2} subpart \textbf{(1)}   follows by taking $r_1 \coloneqq \big(4 \max\{1, \nu\} \log x\big)^{-1}$, choosing $T$ as in its statement (which ensures that ${x^{1/\nu}}{T^{-\delta_0/2}} \, \asymp \, x^{\sigmanu(T)}$),  and lastly, by choosing $u$ so as to equate the first two terms in the above display that contain $u$ (namely, the term with ``$u$" and the first term with ``$u^{-1}$"). \hfill \qedsymbol
\section{Distribution of integers with prime factors restricted to progressions with varying moduli: Proof of Theorem \ref{thm:App1PrimeFactorsinAPs}}\label{sec:App1PrimeFactorsinAPProof}
We intend to apply Theorem \ref{thm:LFuncLSDVariant1} with $\nu=1$ and with $a_n$ being the indicator function of the property that all the prime factors of $n$ lie in the residue classes in $\mcA$. Note that \eqref{eq:Variant1GrowthConditionNew} is tautological. Moreover, since $\{a_n\}_n$ is a multiplicative sequence, we can use the Euler product  
to write
\begin{equation}\label{eq:App1ModifiedEulerProduct}
\sum_n \, \frac{a_n}{n^s} \, = \, \prod_{p \bmod q \, \in \, \mcA} \, \left(1 + \sum_{r \ge 1}, \frac1{p^{rs}}\right) \, = \, \prod_{p \bmod q \, \in \, \mcA} \, \left(1-\frac1{p^s}\right)^{-1} \end{equation}
for all complex $s$ having $\Ree(s)>1$.  Now by 
the orthogonality of Dirichlet characters, we see that for any $b \in U_q$, we have $\sum_{p \equiv b \pmod q} \, 1/p^{s} \, = \, \phi(q)^{-1} \sum_\chi \, \chibar(b) \, \sum_p \, \chi(p)/p^s$. 

Moreover, the Euler product of $L(s, \chi)$ shows that $\log L(s, \chi) = \sum_p \, \chi(p)/p^s \, + \, \sum_{p, r \ge 2} \, \chi(p^r)/rp^{rs}$.  
Eliminating $\sum_p \, \chi(p)/p^s$ from these two identities, and once again invoking orthogonality, we obtain
\begin{equation}\label{eq:pPowersInAP}
\sum_{p \equiv b \pmod q} \, \, \frac1{p^s} \, = \, \frac1{\phi(q)} ~~~ \sum_\chi \, \, \chibar(b) \, \log L(s, \chi) ~~~~~ - \sum_{\substack{p, r \ge 2\\p^r \equiv b \pmod q}} \, \frac1{rp^{rs}}
\end{equation}
for all $b \in U_q$ and all complex numbers $s$ having $\Ree(s)>1$. Now writing    
$$-\sum_{p \bmod q \, \in \, \mcA} \, \log\left(1-\frac1{p^s}\right) \, = \, -\sum_{b \in \mcA} \, \sum_{p \equiv b \pmod q} \, \log\left(1-\frac1{p^s}\right) \, = \,  \sum_{b \in \mcA} \, \sum_{\substack{p, r \ge 1\\p^r \equiv b \bmod q}} \, \frac1{rp^{rs}},$$
and using \eqref{eq:pPowersInAP} replace the contribution of $r=1$, we find that  
\begin{equation}\label{eq:pLogsInAP}
-\sum_{p \bmod q \, \in \, \mcA} \, \log\left(1-\frac1{p^s}\right) \, = \, \frac1{\phi(q)} ~~~ \sum_\chi \, \, \chibar(b) \, \log L(s, \chi) ~~~~~ + \sum_{\substack{p, r \ge 2\\p \bmod q \, \in \, \mcA}} \, \frac1{rp^{rs}} - \sum_{\substack{p, r \ge 2\\p^r \bmod q \, \in \, \mcA}} \, \frac1{rp^{rs}}.
\end{equation}
Inserting \eqref{eq:pLogsInAP} into \eqref{eq:App1ModifiedEulerProduct}, we find that $\sum_{n \ge 1} \, a_n/n^s \, = \, \mathcal F(s) G(s)$, where $\mathcal F(s) \, \coloneqq \, \prod_\chi \, L(s, \chi)^{\alpha_\chi}$, where $\alpha_\chi = \phi(q)^{-1} \sum_{b \in \mcA} \, \chibar(b)$, and where $G(s) \coloneqq \exp\left(\sum_{\substack{p, r \ge 2\\p \bmod q \, \in \, \mcA}} \, 1/{rp^{rs}} - \sum_{\substack{p, r \ge 2\\p^r \bmod q \, \in \, \mcA}} \, 1/{rp^{rs}}\right)$. 
Note that $\sum_{p, r \ge 2} \, |1/rp^{rs}| \, \le \, \sum_{p, r \ge 2} \, p^{-3r/4} $ $ \, \ll \, \sum_p \, p^{-3/2} \, \ll \, 1$ uniformly in $s$ with $\Ree(s) \ge 3/4$. Hence $G(s)$ is analytic on the half--plane $\Ree(s) \ge 3/4$ and satisfies $G(s) \ll 1$ uniformly therein (with the last implied constant being absolute). This shows that the sequence $\{a_n\}_n$ satisfies property $\mathcal P(1, \{\alpha_\chi\}_\chi; c_0, \Omega)$, with $\alpha_\chi$ as defined above and with $\Omega(t) \ll 1$. 

We now observe one of the first concrete applications showing the benefit of introducing the parameter $\lambda_q$ in our main result Theorem \ref{thm:LFuncLSDVariant1}: Note that for any coprime residue class $a$ mod $q$,  
$$\sum_\chi \, \alpha_\chi \cdot \chi(a) \, = \,  \frac1{\phi(q)}\sum_{b \in \mcA} \, \sum_\chi \,  \chibar(b) \chi(a) \, = \, \sum_{b \in \mcA} \, \bbm_{b \equiv a \pmod q} \, = \, \bbm_{a \in \mcA}.$$
This shows that $|\lambda_q| \le 2$. The asymptotic formula \eqref{eq:App1PrimeFactorsinAPs} now follows from Theorem \ref{thm:LFuncLSDVariant1}(1) in the case  $q \in [e^6, (\log x)^{K_0}]$,  
upon using Siegel's Theorem to note that $1-\eta_e \, \ge \, c_1 \cdot q^{-\epsilon_0/K_0}$ for some constant $c_1>0$ depending only on $\epsilon_0$ and $K_0$. Likewise, \eqref{eq:App1PrimeFactorsinAPs} follows from the last assertion of Theorem \ref{thm:LFuncLSDVariant1}(2) in the case when the Siegel zero does not exist.

Finally, we use  \eqref{eq:lambdaq,mujDef} to see  that the coefficients $k_j$ in \eqref{eq:App1ModifiedEulerProduct} are given by 
\begin{equation}\label{eq:App1PrimeFactorsinAPkj}
k_j  \, = \, \frac1{j!} \cdot \frac{\mathrm d^j}{\mathrm{d}s^j}\Bigg\vert_{s=1}\frac{\mathcal F(s) G(s)}s\left(s-1\right)^{\alphachiZ},
\end{equation}
with $\mathcal F(s)$, $G(s)$ and $\alpha_\chi$ as defined after \eqref{eq:pLogsInAP}. In particular, it is worth noting that 
$$k_0 \, = \, \lim_{s \rightarrow 1+} \, \frac{\mathcal F(s) G(s)}s\left(s-1\right)^{\alphachiZ} \, = \,  G(1) \cdot \left(\prod_{\chi \ne \chi_0} L(1, \chi)^{\alpha_\chi}\right) \cdot \lim_{s \rightarrow 1+} \, \big(L(s, \chi_0)(s-1)\big)^{\alphachiZ},$$
Now since $L(s, \chi_0) = \zeta(s) \cdot \prod_{\ell \mod q} (1-1/\ell^s)$, and since $\lim_{s \rightarrow 1+} \, \zeta(s)(s-1) = 1$, we obtain the following explicit formula for the first coefficient $k_0$ of the asymptotic series \eqref{eq:App1PrimeFactorsinAPs}.    
\begin{equation}\label{eq:k0Formula}
 k_0 \, = \, \left(\frac{\phi(q)}q\right)^{\frac{|\mcA|}{\phi(q)}} \left(\prod_{\chi \ne \chi_0} L(1, \chi)^{\alpha_\chi}\right) \cdot  
 \exp\Bigg(\sum_{\substack{p, r \ge 2\\p \bmod q \, \in \, \mcA}} \, \frac1{rp^{r}} - \sum_{\substack{p, r \ge 2\\p^r \bmod q \, \in \, \mcA}} \, \frac1{rp^{r}}\Bigg). 
\end{equation}
 
\section{Distribution of the least invariant factor of the multiplicative group:\\ Proof of Theorems \ref{thm:Dqx} and \ref{thm:Eqx}}\label{sec:App2ChangMartin} 
The following algebraic result characterizes when an even integer $q>2$ divides the least invariant factor $\lambdaSt(n)$ of the multiplicative group modulo $n$. This is a restatement of \cite[Proposition 5.3]{CM20}. 
In this section, $\boldsymbol{P(q)}$ denotes the \textbf{largest} prime divisor of $q$ and $\boldsymbol{e_p = v_{P(q)}(q)}$ is the exponent (highest power) of $P(q)$ in the prime decomposition of $q$.  
\begin{lem}\label{lem:InvFactorDivisiblityCriteria}
Let $q \ge 4$ be an even integer.\\ 
$\mathrm{(i)}$ If $q \nmid \phi(P(q)^{e_p+1})$, then the positive integers $n$ for which $q \mid \lambdaSt(n)$ are precisely those of the form $n = 2^{e_2} \cdot m$, where $e_2 \in \{0, 1\}$ and $m$ is only divisible by primes $\equiv 1 \pmod q$. 

$\mathrm{(ii)}$ If $P(q)>2$ and $q \mid \phi(P(q)^{e_p+1})$, then the positive integers $n$ for which $q \mid \lambdaSt(n)$ are precisely those of the form $n = 2^{e_2} \cdot P(q)^v \cdot m$, where $e_2 \in \{0, 1\}$, where $m$ is only divisible by primes $\equiv 1 \pmod q$, and where $v$ is either $0$ or at least $e_p+1$.
\end{lem}
To see that the above statement is indeed equivalent to  \cite[Proposition 5.3]{CM20}, 
we just need to note that if $q \mid \ell^{v_\ell(q)}(\ell-1)$ for some prime $\ell \mid q$, then we 
have $P(q) \mid \ell^{v_\ell(q)}(\ell-1)$, forcing $\ell = P(q)$.  
In order to establish Theorem \ref{thm:Dqx}, we use the above characterization to see that for all complex $s$ with $\sigma = \Ree(s)>1$, we may write 
\begin{align*}\allowdisplaybreaks
\sum_n \, \frac{\bbm_{q \mid \lambdaSt(n)}}{n^s} \, &= \, 
\begin{multlined}[t]
\sum_{e_2 \in \{0, 1\}} \, \frac1{2^{e_2s}} \, \sum_{m:\, p \mid m \implies p \equiv 1 \pmod q} \, \frac1{m^s} \,\\ + \, \bbm_{q \, \mid \, \phi(P(q)^{e_p+1})} \, \sum_{e_2 \in \{0, 1\}} \, \frac1{2^{e_2s}} \, \sum_{r \ge e_p+1} \, \frac1{P(q)^{rs}} \, \sum_{m:\, p \mid m \implies p \equiv 1 \pmod q} \, \frac1{m^s}
\end{multlined} \\
&= \, \left(1+\frac1{2^s}\right) \cdot \left(1\, + \, \bbm_{q \mid \phi(P(q)^{e_p+1})}\, \sum_{r \ge e_p+1} \, \frac1{P(q)^{rs}}\right) \cdot \prod_{p \equiv 1 \pmod q} \left(1+\sum_{r \ge 1} \frac1{p^{rs}}\right)
\end{align*}
Invoking \eqref{eq:pLogsInAP} with $\mathcal A \coloneqq \{1\} \subset U_q$, we get $\sum_n \bbm_{q \mid \lambdaSt(n)}/n^s = \left(\prod_\chi \, L(s, \chi)\right)^{1/\phi(q)} \cdot \mathcal G(s)$, where 
\begin{multline*}\allowdisplaybreaks
\mathcal G(s) \coloneqq \left(1+\frac1{2^s}\right) \cdot \left(1 \, + \, 
\frac{P(q)^{-se_p}}{P(q)^s-1}\right)^{\bbm_{q \mid \phi(P(q)^{e_p+1})}} \cdot \exp\left(\sum_{\substack{p, r \ge 2\\p \equiv 1 \pmod q}} \, \frac1{rp^{rs}} - \sum_{\substack{p, r \ge 2\\p^r \equiv 1 \pmod q}} \, \frac1{rp^{rs}}\right).
\end{multline*}
Once again, since $\sum_{p, r \ge 2} \, |1/rp^{rs}| \ll 1$ uniformly in the half plane $\Ree(s) \ge 3/4$, we find that $\mathcal G(s)$ is analytic and of size $O(1)$ uniformly on this half place. Hence, our sequence $\{\bbm_{q \mid \lambdaSt(n)}\}_n$ has property $\mathcal P(1, \{\alpha_\chi\}_\chi; c_0, \Omega)$, with   $\Omega(t) \ll 1$ and with all $\alpha_\chi = 1/\phi(q)$. As such, $\sum_\chi \, \alpha_\chi \cdot \chi(a) \, = \, \bbm_{a \equiv 1 \bmod q}$ for any coprime residue class $a$ mod $q$, showing that once again $|\lambda_q| \le 2$ in Theorem \ref{thm:LFuncLSDVariant1}. The rest of the argument for Theorem \ref{thm:Dqx} goes through exactly as in the previous section. 

The coefficients $r_j$ (and in particular $r_0$) have analogous explicit formulas to \eqref{eq:App1PrimeFactorsinAPkj} and \eqref{eq:k0Formula} respectively, and they are proven analogously as well. We spell them out below, for completeness. 
\begin{align}\label{eq:App2ChangMartinrj}\allowdisplaybreaks
r_j \, = \, \frac1{j!} \cdot \frac{\mathrm d^j}{\mathrm{d}s^j}\Bigg\vert_{s=1}\frac{\mathcal G(s)}s \cdot \left((s-1) \cdot \prod_\chi L(s, \chi)\right)^{1/\phi(q)}.
\end{align}
In particular, since $L(1, \chi) \ll \log q$ for all nontrivial characters $\chi$ mod $q$, we see that    
$$r_0 \, = \, \mathcal G(1) \cdot   \left(\frac{\phi(q)}q \, \prod_{\chi \ne \chi_0} L(1, \chi)\right)^{1/\phi(q)} \asymp \,  \left(\frac{\phi(q)}q \, \prod_{\chi \ne \chi_0} |L(1, \chi)|\right)^{1/\phi(q)} \ll \log q.$$
This completes the proof of Theorem \ref{thm:Dqx}. \hfill \qedsymbol
\begin{rmk}
This approach of directly  applying our Theorem \ref{thm:LFuncLSDVariant1} on the non--multiplicative sequence $a_n = \bbm_{q \mid \lambdaSt(n)}$ also simplifies Chang and Martin's approach in \cite[Proposition 5.4]{CM20}  of first estimating the count $\mathcal N(x; q, 1)$, and then passing to the count $\sum_{n \le x} \, \bbm_{q \, \mid \, \lambdaSt(n)}$ via an additional combinatorial argument.  
\end{rmk}
We now come to Theorem  \ref{thm:Eqx}; \textbf{in the rest of this section, we assume that $\boldsymbol{q \le (\log x)^K}$ with $\boldsymbol{K>0}$ fixed}. We will refine some of the ideas of Chang and Martin with an additional input from the anatomy of integers, where we split off the largest prime factor of $n$ and carefully study the case when this factor grows somewhat rapidly with $x$. To this end, we define $\boldsymbol{z \coloneqq x^{1/\log_2 x}}$, and let $\boldsymbol{\Dqx}$ denote the number of positive integers $n \le x$ satisfying all the following properties 
\begin{itemize}
\item $q \mid \lambdaSt(n)$.
\item $P(n)>z$, \, $P(n)^2 \, \nmid \, n$. 
\item $\#\{p \mid n:~ p \equiv 1 \pmod q\} \ge 2$.
\end{itemize}
We define $\boldsymbol{\Eqx}$ the same way, only with the first (divisbility) condition replaced by the equality $q=\lambdaSt(n)$.  We start by giving the following upper bounds on $\Dqx$. 
\begin{lem}\label{lem:DQxUniformBounds} Uniformly in $x \ge 16$, we have 
\begin{equation}\label{eq:DQxUniformBounds}
\DQx \, 
\begin{cases}\vspace{1em}
\ll \, \displaystyle{\frac{x}{(\log x)^{1-1/\phi(Q)}} \cdot \left(\frac{\log_2 x}{\phi(Q)}\right)^2}, &\text{ uniformly in even }Q \le (\log x)^{100K}.\\
\le \, \displaystyle{x \, \left(\frac{\log_2 x}{\phi(Q)} \, + \, O\left(\frac{\log Q}Q\right)\right)^2}, &\text{ uniformly in even }Q \in [2, x].
\end{cases}
\end{equation}
\end{lem}
\begin{proof}[Proof of Lemma \ref{lem:DQxUniformBounds}]
To see the first bound, we consider any $Q \le (\log x)^K$, and 
write $\DQx = \mathcal D_{\mathrm{odd}}(x) + \mathcal D_{\mathrm{even}}(x)$, where $\mathcal D_{\mathrm{odd}}(x)$ and $\mathcal D_{\mathrm{even}}(x)$ denote the contribution of all odd and even  positive integers counted in $\DQx$, respectively. By Lemma \ref{lem:InvFactorDivisiblityCriteria} and the definition of $\DQx$, any positive integer $n$ counted in $\mathcal D_{\mathrm{odd}}(x)$ can be written as $m p P$, where $P = P(n) > \max\{z, P(mp)\}$, where $P \equiv p \equiv 1 \pmod Q$, where any prime $\ell$ dividing $m$ must either be equal to $P(Q)$ or must be $1$ mod $\ell$. 
(Here we have noted that $P = P(n)>z>q \ge P(Q)$.)       

Now given $m$ and $p$, the number of primes $P \in (z, x/mp]$ satisfying $P \equiv 1 \pmod Q$ is $\ll x/\phi(Q) mp \log z$ by the Brun--Titchmarsh inequality. Moreover, by Brun--Titchmarsh and partial summation, we have $\sum_{\substack{p \le x\\p \equiv 1 \pmod Q}} \, 1/p \, \ll \, \log_2 x/\phi(Q)$. This bound also shows that  
\begin{align*}\allowdisplaybreaks
\sum_{\substack{m \le x\\ \ell \mid m \implies \ell = P(Q) \text{ or }\ell \equiv 1 \pmod Q}} \,&\,\, \frac1m \, \le \, \left(1+\sum_{r \ge 1} \, \frac1{P(Q)^r}\right) \cdot \prod_{\substack{\ell \le x\\\ell \equiv 1 \pmod Q}} \left(1 \, + \, \sum_{r \ge 1} \, \frac1{\ell^r}\right)\\
&\le \, \left(1+\sum_{r \ge 1} \, \frac1{2^r}\right) \cdot \exp\Bigg(\sum_{\substack{\ell \le x\\\ell \equiv 1 \pmod Q}} \, \frac1\ell \, + \, O\left(\sum_{\ell, r \ge 2} \, \frac1{\ell^r}\right)\Bigg) \, \ll \, (\log x)^{1/\phi(Q)}.
\end{align*}
Collecting all these observations, we find that $\mathcal D_{\mathrm{odd}}(x) \, \ll \, x \, (\log x)^{\phi(Q)-1} \cdot (\log_2 x/\phi(Q))^2$. An entirely analogous argument shows that $\mathcal D_{\mathrm{even}}(x)$ satisfies the same bound: We just starting by writing all the $n$ counted in $\mathcal D_{\mathrm{even}}(x)$ as $2m p P$, with $m, p, P$ satisfying similar conditions as above. (Here we recalled that $v_2(n) \in \{0, 1\}$ by Lemma \ref{lem:InvFactorDivisiblityCriteria}.)  This shows the first bound in \eqref{eq:DQxUniformBounds}. 

For the second bound, we simply note that by the last condition in the definition of $\DQx$, 
$$\DQx \, \le \, \sum_{\substack{p_1, p_2 \le x\\p_1 \equiv p_2 \equiv 1 \pmod Q}} ~~~~ \sum_{\substack{n \le x\\p_1p_2 \mid n}} \, 1 \, \, \le \, \, \sum_{\substack{p_1, p_2 \le x\\p_1 \equiv p_2 \equiv 1 \pmod Q}} ~~~ \frac{x}{p_1 p_2} ~\le~ x\, \Bigg(\sum_{\substack{p \le x\\p \equiv 1 \pmod Q}} \, \frac1p\Bigg)^2.$$
The second bound in \eqref{eq:DQxUniformBounds} now follows by using an estimate due independently to Pomerance (see Remark 1 of \cite{pomerance77}) and Norton (see the Lemma on p.\ 699 of \cite{norton76}). 
\end{proof}
Using Lemma \ref{lem:DQxUniformBounds}, we will now provide an estimate on $\Eqx$ (also defined before the lemma) uniformly for $q \le (\log x)^K$. To do this, we adapt the argument given for Lemma 6.2 in \cite{CM20}. Note that by definition of $\Dqx$ and $\Eqx$, we have $\Dqx = \sum_{m \ge 1} \, \err_{mq}(x)$, so that by a version of the M\"obius Inversion Formula, we have 
\begin{equation}\label{eq:EqxMobousInversion}
\Eqx = \sum_{m \ge 1} \, \mu(m) \mathcal D_{mq}(x).     
\end{equation}
(Here both the seemingly infinite sums are actually finite, and only go up to $m \le x/q$.) Now from the first bound in \eqref{eq:DQxUniformBounds} and the well--known estimate $\phi(mq) \gg mq/\log_2 (mq)$, we see that 
\begin{align}\allowdisplaybreaks\label{eq:DmqSmallm}
\sum_{2 \le m \le (\log x)^{3K+3}} \, \mathcal D_{mq}(x) \, &\ll \, \sum_{2 \le m \le (\log x)^{3K+3}} \, \frac{x}{(\log x)^{1-1/\phi(mq)}} \cdot \left(\frac{\log_2 x \cdot \log_2(mq)}{mq}\right)^2 \nonumber\\
&\ll \, \frac{x \, (\log_2 x \cdot \log_3 x)^2}{q^2 \, (\log x)^{1-1/2\phi(q)}} \, \sum_{m \ge 2} \, \frac1{m^2} \ll \, \frac{x \, (\log_2 x \cdot \log_3 x)^2}{q^2 \, (\log x)^{1-1/2\phi(q)}}, 
\end{align}
where in the second line, we have used that $mq \le (\log x)^{4K+3}$,  and that $\phi(mq) \ge 2\phi(q)$ for any $m \ge 2$ as $q$ is even. On the other hand, using the second bound in \eqref{eq:DQxUniformBounds}, we find that 
\begin{align}\allowdisplaybreaks\label{eq:DmqLargem}
\sum_{(\log x)^{3K+3} < m \le x/q} \, \mathcal D_{mq}(x) \, &\ll \, \sum_{(\log x)^{3K+3} < m \le x/q} \, \, x \cdot \left(\frac{\log_2 x \cdot \log_2(mq)}{mq} \, + \, \frac{\log (mq)}{mq}\right)^2 \nonumber\\
&\ll \, x \, \sum_{m> (\log x)^{3K+3}} \,  \left(\frac{\log x}{mq}\right)^2 \, \ll \, \frac{x}{q^2\, (\log x)^{3K}}. 
\end{align}
Inserting \eqref{eq:DmqSmallm} and \eqref{eq:DmqLargem} into \eqref{eq:EqxMobousInversion}, we obtain 
\begin{equation}\label{eq:EqxVsDqxExplicit}
\Eqx \, = \, \Dqx \, + \, O\left(\frac{x \, (\log_2 x \cdot \log_3 x)^2}{q^2 \, (\log x)^{1-1/2\phi(q)}}\right).    
\end{equation}
Now by known estimates on smooth numbers (see \cite[p.~15]{CEP83} or \cite[Theorem 5.13\text{ and }Corollary 5.19, Chapter III.5]{tenenbaum15}), the number of $n \le x$ having $P(n) \le z = x^{1/\log_2 x}$ is $\ll x/(\log x)^{100K}$. Hence  
\begin{align*}\allowdisplaybreaks
\sum_{\substack{n \, \le \, x\\q \, \mid \, \lambdaSt(n)}} \, 1 ~~~ = ~~~ \Dqx ~~~~ + \sum_{\substack{n \, \le \, x: ~~ q \, \mid \, \lambdaSt(n)\\P(n)>z, ~~ P(n)^2 \, \nmid \, n\\\#\{p \mid n: ~ p\equiv 1 \pmod q\} \, \in \,\{0, 1\} }} 1 ~~~~+~~~~ O\left(\frac x{(\log x)^{100K}}\right).\\
\sum_{\substack{n \, \le \, x\\q \, = \, \lambdaSt(n)}} \, 1 ~~~ = ~~~ \Eqx ~~~~ + \sum_{\substack{n \, \le \, x: ~~ q \, = \, \lambdaSt(n)\\P(n)>z, ~~ P(n)^2 \, \nmid \, n\\\#\{p \mid n: ~ p\equiv 1 \pmod q\} \, \in \,\{0, 1\} }} 1 ~~~~+~~~~ O\left(\frac x{(\log x)^{100K}}\right).
\end{align*}
Subtracting the second estimate from the first, and using \eqref{eq:EqxVsDqxExplicit}, we obtain 
\begin{equation}\label{eq:EqVsDqNearFinal}
\sum_{\substack{n \le x\\q \, = \, \lambdaSt(n)}} \, 1 ~~~~=~~~~ \sum_{\substack{n \le x\\q \, \mid \, \lambdaSt(n)}} \, 1 ~~~~-  \sum_{\substack{n \, \le \, x: ~~ q \, \mid \, \lambdaSt(n), ~~ \lambdaSt(n) \, > \, q\\P(n)>z, ~~ P(n)^2 \, \nmid \, n\\\#\{p \mid n: ~ p\equiv 1 \pmod q\} \, \in \,\{0, 1\} }} 1 ~~+~~ O\left(\frac{x \, (\log_2 x \cdot \log_3 x)^2}{q^2 \, (\log x)^{1-1/2\phi(q)}}\right).
\end{equation} 
We now estimate the second sum on the right hand above, which we call $\mathcal D_{\mathrm{rem}}$.  

\textbf{Case 1.} If $q \nmid \phi(P(q)^{e_p+1})$, then Lemma \ref{lem:InvFactorDivisiblityCriteria}(i) shows that any $n$ counted in the sum $\mathcal D_{\mathrm{rem}}$ must necessarily be of the form $n=P$ or $n=2P$ where $P \in (z, x)$ satisfies $P \equiv 1 \pmod q$. This condition is also sufficient since from the isomorphisms $U_P \cong \Z/(P-1)\Z$ and $U_{2P} \cong U_2 \times U_P \cong \Z/(P-1)\Z$, we have $\lambdaSt(2P) = \lambdaSt(P) = P-1 \equiv 0 \pmod q$. Hence by the Siegel--Walfisz Theorem,
\begin{equation}\label{eq:DremCase1}
\sum_{\substack{n \, \le \, x: ~~ q \, \mid \, \lambdaSt(n), ~~ \lambdaSt(n) \, > \, q\\P(n)>z, ~~ P(n)^2 \, \nmid \, n\\\#\{p \mid n: ~ p\equiv 1 \pmod q\} \, \in \,\{0, 1\} }} 1 ~~=~~ \sum_{\substack{z<P \le x\\P \equiv 1 \pmod q}}~1~~+~~\sum_{\substack{z<P \le x/2\\P \equiv 1 \pmod q}}~1 ~~=~~ \frac{\li(x)+\li(x/2)}{\phi(q)} \, + \, O\big(x\exp(-\sqrt{\log x})\big).
\end{equation}
\textbf{Case 2.} Now assume that $q \mid \phi(P(q)^{e_p+1})$. By Lemma \ref{lem:InvFactorDivisiblityCriteria}(ii), any $n \le x$ satisfying $q \mid \lambdaSt(n)$ but having no prime factor $\equiv 1 \pmod q$, must be of the form $n = 2^{e_2} P(q)^v$. Hence the number of such $n$ is $\ll \, \sum_{0 \le x \le \log x} \, 1 \ll \log x$. On the other hand, if any $n$ counted in $\mathcal D_{\mathrm{rem}}$ has exactly one prime factor $\equiv 1 \pmod q$, then Lemma \ref{lem:InvFactorDivisiblityCriteria}(ii) shows that $n$ must be of the form $2^{e_2} P(q)^v P$ with $e_2 \in \{0, 1\}$, with $P \in (z, x]$ satisfying $P \equiv 1 \pmod q$, and with $v \in \{0, e_p+1, e_p+2, \dots\}$. Hence the total number of $n$ counted in $\mathcal D_{\mathrm{rem}}$ exactly one prime factor in the residue class $1$ mod $q$ is
$$\sum_{\substack{z<P \le x\\P \equiv 1 \pmod q}}~1~~+~~\sum_{\substack{z<P \le x/2\\P \equiv 1 \pmod q}}~1 ~~ + O\Bigg(\sum_{e_2 \in \{0, 1\}} ~~ \sum_{v \ge e_p+1} ~~ \sum_{\substack{z \, < \, P \, \le \, x/2^{e_2}P(q)^v\\P \equiv 1 \pmod q}}~~1\Bigg),$$
where the first two sums come from the case $v=0$ (and then correspond to whether $e_2 = 0$ or $e_2 = 1$). Now the first two sums above can again be estimated by the Siegel--Walfisz Theorem, whereas by the Brun--Titchmarsh theorem, the total $O$--term in the above display is 
$$\ll \, \frac{x}{\phi(q)\log z} \sum_{e_2 \in \{0, 1\}} \frac1{2^{e_2}} ~\sum_{v \ge e_p+1} \frac1{P(q)^v} ~\ll~ \frac{x}{P(q)^{e_p+1}\phi(q)\log z} ~\ll~ \frac{x \, \log_2 x \cdot \log_3 x}{q^2 \, \log x}$$
where we have noted that $q \mid \phi(P(q)^{e_p+1})$ forces $q \le P(q)^{e_p+1}$. By the discussion under Case 2, 
\begin{equation}\label{eq:DremCase2}
\sum_{\substack{n \, \le \, x: ~~ q \, \mid \, \lambdaSt(n), ~~ \lambdaSt(n) \, > \, q\\P(n)>z, ~~ P(n)^2 \, \nmid \, n\\\#\{p \mid n: ~ p\equiv 1 \pmod q\} \, \in \,\{0, 1\} }} 1 ~~=~~ \frac{\li(x)+\li(x/2)}{\phi(q)} \, + \, O\left(\frac{x \, \log_2 x \cdot \log_3 x}{q^2 \, \log x}\right).
\end{equation}
Finally, inserting \eqref{eq:DremCase1} and \eqref{eq:DremCase2} into \eqref{eq:EqVsDqNearFinal}, we find that 
$$\sum_{\substack{n \le x\\q \, = \, \lambdaSt(n)}} \, 1 ~~~~=~~~~ \sum_{\substack{n \le x\\q \, \mid \, \lambdaSt(n)}} \, 1 ~~-~~\frac{\li(x)+\li(x/2)}{\phi(q)} ~~+~~ O\left(\frac{x \, (\log_2 x \cdot \log_3 x)^2}{q^2 \, (\log x)^{1-1/2\phi(q)}}\right),$$
whereupon Theorem \ref{thm:Eqx} follows from the very first assertion in Theorem \ref{thm:Dqx}. \hfill \qedsymbol
\section{Distribution of the least primary factor of the multiplicative group: Proof of Theorem \ref{thm:MartinNguyen}}
As alluded to in subsection \cref{subsec:MartinNguyen}, our algebraic input will be the same as in Martin and Nguyen's work \cite{MN24}, namely we will start by writing $\#\{n \le x: \lambda'(n) = q\} = \mathcal A_q(x) - \mathcal A_{q^+}(x)$, and then estimating $\mathcal A_q(x)$. The following characterization of $\mathcal A_q(x)$ will be really useful for this purpose: This is a restatement of \cite[Proposition 3.6]{MN24}, but we also summarize their argument. 
\begin{lem}\label{lem:LeastPrimaryFactorCongruences}
Given a prime power $q \ge 3$, the integers $n \ge 3$ having $\lambda'(n) \ge q$ are precisely those of the form $2^{e_2}m$, with $e_2 \in \{0, 1\}$, and with the odd integer $m$ satisfying the following property:
\begin{equation}\label{eq:PrimFactOddPart}
\text{For all primes }p \mid m\text{ and all primes }\ell < q,\text{ either }\ell \nmid (p-1)\text{ or }p \equiv 1 \pmod{\ellceilpower}.
\end{equation}
In particular, any such $n$ cannot be divisible by any odd prime at most $q$.
\end{lem}
\begin{proof}[Proof of Lemma \ref{lem:LeastPrimaryFactorCongruences}] Since $U_n \cong \prod_{p^k \parallel n} \, U_{p^k}$, we have $\lambda'(n) = \min_{p^k \parallel n} \, \lambda'(p^k)$, so that the condition $\lambda'(n) \ge q$ amounts to $\lambda'(p^k) \ge q$ for all prime powers $p^k \parallel n$. Now if $p = 2$ but $k \ge 2$, then from $U_{2^k} \cong U_2 \times U_{2^{k-2}}$, we see that $\lambda'(2^k) = 2 < q$. This forces $4 \nmid n$. (Conversely, $4 \nmid n$ also guarantees that $\lambda'(2^{v_2(n)}) = \infty$ by our convention.) On the other hand, if $p>2$, then from $U_{p^k} \cong \Z/(p-1)\Z \times \Z/p^{k-1}\Z$, we see that $\lambda'(p^k)$ is the smallest prime power $\ell^j \parallel (p-1)$, so that the condition $\lambda'(p^k) \ge q$ amounts to having $\ell^j \ge q$ for all $\ell^j \parallel (p-1)$, which in turn is equivalent to having $v_\ell(p-1) \ge \log q/\log \ell$, i.e., to $p \equiv 1 \pmod{\ellceilpower}$.  

This proves the equivalence in the first assertion of the lemma. For the final assertion, note that if an odd prime $p \le q$ divides $n$, then any prime divisor $\ell$ of $p-1$ violates \eqref{eq:PrimFactOddPart}. 
\end{proof}
In this section, we define $Q \, \coloneqq \,  \prod_{\ell < q} \, \ellceilpower$. By the prime number theorem, we have $e^{q/2} \le Q \le e^{4q}$ for all sufficiently large $q$. Note that condition \eqref{eq:PrimFactOddPart} amounts to throwing each odd prime factor of $n$ into $\prod_{\ell<q} \,\big(1 \, + \, (\ell-2) \cdot \ell^{\ceillogqlogell - 1} \big)$ many coprime residue classes modulo $Q$. We use  $\mathcal B_Q$ to denote this set of residues mod $Q$, so that the quantity $B(q)$ in \eqref{eq:BqDef} is exactly $\#\mathcal B_Q/\phi(Q)$. (Note that the residue classes in $\mathcal B_Q$ are precisely those that for each prime $\ell \in (2, q)$, either leave a remainder of $1$ modulo $\ellceilpower$, or are obtained by lifting the residues $\{2, \dots, \ell-1\}$ from modulus $\ell$ to modulus $\ellceilpower$.)

Hereafter proceeding exactly as we did in section \ref{sec:App1PrimeFactorsinAPProof} with the integer $Q$ playing the role of ``$q$" in that section, we find that with $\bbm_{(7.1)}$ being the indicator function of condition \eqref{eq:PrimFactOddPart}, 
$$\sum_n \, \frac{\bbm_{\lambda'(n) \ge q}}{n^s} \, = \, \sum_{2 \, \nmid \, n} \, \frac{\bbm_{(7.1)}}{n^s} \, + \, \frac1{2^s}\sum_{2 \, \nmid \, n} \, \frac{\bbm_{(7.1)}}{n^s} \, = \, \left(1+\frac1{2^s}\right) \cdot \prod_{p \bmod q \, \in \, \mathcal B_Q} \, \left(1-\frac1{p^s}\right)^{-1}$$
 can be written in the form $\left(\prod_{\chi \bmod Q} \, L(s, \chi)^{\alpha_\chi} \right) \cdot G_Q(s)$, with $\alpha_\chi \coloneqq \phi(Q)^{-1} \sum_{b \in \mathcal B_Q} \, \chibar(b)$ (so that $\alpha_{\chi_0} = B(q)$), and with the function 
 $$G_Q(s) \, \coloneqq \, \left(1+\frac1{2^s}\right) \left(1-\frac1{2^s}\right)^{\bbm_{2 \in \mathcal B_Q}} \cdot \exp\Bigg(\sum_{\substack{p, r \ge 2\\p \bmod q \, \in \, \mathcal B_Q}} \, \frac1{rp^{rs}} - \sum_{\substack{p, r \ge 2\\p^r \bmod q \, \in \, \mathcal B_Q}} \, \frac1{rp^{rs}}\Bigg)$$
 being analytic of size $O(1)$ uniformly in the half plane $\sigma \ge 3/4$. This shows that the sequence $\{\bbm_{S(n) \ge q})_n$ has property $\mathcal P(1, \{\alpha_\chi\}_\chi; c_0, \Omega)$ with $\Omega(t) \ll 1$. Theorem \ref{thm:LFuncLSDVariant1} thus gives
 \begin{equation}\label{eq:AqxLSD}
    \mathcal A_q(x) = \, \frac{x}{(\log x)^{1-B(q)}} \sum_{j=0}^N \, \frac{\kappa_j(q) \, (\log x)^{-j}}  {\Gamma\left(B(q)-j\right)} \, + \,O\left(\frac{N!\,(142 \, c_1^{-1})^N \cdot x}{(\log x)^{(N+2)(1-\epsilon_0) - \frac{1}{\phi(q)}}} \, + \,  x\exp\left(-\sqrt{\frac{c_0\log x}{32}}\right) \right),
 \end{equation} 
uniformly in $Q \le e^{4K} \le (\log x)^{4K}$ (and uniformly in $Q \le \exp(\sqrt{c_0\log x/8}$ if the Siegel zero does not exist), with the coefficients $\kappa_j$ being defined by 
\begin{equation}\label{eq:kappajMartNg}
\kappa_j = \frac1{j!} \cdot \frac{\mathrm d^j}{\mathrm{d}s^j}\Bigg\vert_{s=1}\frac{G_Q(s)}s \cdot \left(\prod_{\chi \bmod Q} \, L(s, \chi)^{\alpha_\chi}\right) \cdot \left(s-1\right)^{B(q)}.
\end{equation}
Note that the exponents $\alpha_\chi \, = \, \phi(Q)^{-1} \sum_{b \in \mathcal B_Q} \, \chibar(b)$ have been computed explicitly in \cite[Section 6]{MN24}. Finally, Theorem \ref{thm:MartinNguyen} follows from \eqref{eq:AqxLSD} by writing $\#\{n \le x: \lambda'(n)\} = \mathcal A_q(x) -  \mathcal A_{q^+}(x)$, and then noting that $B(q^+) \le B(q)$ by \eqref{eq:BqDef} as well as $q^+ \le 2q$ by  Bertrand's Postulate. \hfill \qedsymbol
\section{The Sathe--Selberg Theorem in Arithmetic Progressions to varying moduli: Proofs of Theorems \ref{thm:omegaa(n)} and \ref{thm:Omegaa(n)}}
We first prove Theorem \ref{thm:omegaa(n)} and then mention the changes required for Theorem \ref{thm:Omegaa(n)}. 
{In this entire section, we may assume that $\boldsymbol{q> \max\{e^6, \,  25(K+2)^2\}}$, for otherwise all assertions can be deduced by following our arguments below, and by replacing our upcoming uses of Theorem \ref{thm:LFuncLSDVariant2} by the standard form of the Landau--Selberg--Delange method as in \cite[Theorem II.5.2]{tenenbaum15}.} 

Inspired by Selberg's idea in \cite{selberg1954}, we will identify the number $\#\{n \le x: \, \omegaan = k\}$ as the coefficient of $\zPoweromegaan$ in the polynomial $\sum_{n \le x} \, \zPoweromegaan$. We thus start by giving an estimate on $\sum_{n \le x} \, \zPoweromegaan$ for complex numbers $z$ having $|z| \le K+1$. To this end, note that since $\omegaan$ is an additive function, the function $n \mapsto \zPoweromegaan$ is multiplicative.  
Hence for all $s$ with $\sigma = \Ree(s)>1$, we have  
\begin{align}\label{eq:omegaaEulerProdModif}\allowdisplaybreaks
\sum_n \, \frac{\zPoweromegaan}{n^s} \, &= \, \prod_{p} \, \left(1+\sum_{r \ge 1} \, \frac{z^{\omega_a(p^r)}}{p^{rs}}\right) \, = \, \prod_{p \equiv a \pmod q} \, \left(1+z\sum_{r \ge 1} \, \frac1{p^{rs}}\right) \cdot \prod_{p \not\equiv a \pmod q} \, \left(1+\sum_{r \ge 1} \, \frac{1}{p^{rs}}\right)\\
&= \, \zeta(s) \cdot \prod_{p \equiv a \pmod q} \, \left(1+\frac{z}{p^s} \cdot \left(1-\frac1{p^s}\right)^{-1} \right) \, \left(1-\frac1{p^s}\right) \, = \, \zeta(s) \cdot \prod_{p \equiv a \pmod q} \, \left(1+\frac{z-1}{p^s} \right)
\end{align}
where we have used the Euler product of $\zeta(s)$ and the fact that $\omega_a(p^r) = \bbm_{p \equiv a \pmod q}$. Continuing, 
\begin{equation*}
\sum_n \, \frac{\zPoweromegaan}{n^s} \, = \, \zeta(s) \cdot \exp\left(\sum_{p \equiv a \pmod q} \frac{z-1}{p^s}\right) \cdot \left(\prod_{p \equiv a \pmod q} \, \left(1+\frac{z-1}{p^s} \right) \, e^{(1-z)/p^s}\right).
\end{equation*}
Invoking \eqref{eq:pPowersInAP} for the sum on the right hand side above and using the fact that $e^{(1-z)/p^s} \, = \, 1 \, - \, (z-1)/p^s \, + \, \sum_{r \ge 2} \, (1-z)^r/r! \, p^{rs}$ and that $\zeta(s) = L(s, \chi_0) \, \prod_{\ell \mid q} \, (1-1/\ell^s)^{-1}$, we thus obtain 
\begin{equation}\label{eq:omegaanLSDPrep}
\sum_n \, \frac{\zPoweromegaan}{n^s} \, = \, \left(\prod_{\chi} \, L(s, \chi)^{\bbm_{\chi = \chi_0} \, + \, (z-1)\chibara/\phi(q)}\right) \cdot \Gzs,
\end{equation}
where the function $\Gzs$ is holomorphic on the half plane $\sigma \ge 3/4$, and has the formula
\begin{multline}\label{eq:Gzs}\allowdisplaybreaks
\Gzs = \exp\Bigg((1-z)\, \sum_{\substack{p, r \ge 2\\p^r \equiv a \pmod q}} \, \frac1{rp^{rs}}\Bigg) \cdot \left(\prod_{\ell \mid q} \, \left(1-\frac1{\ell^s}\right)^{-1}\right) \cdot \left(\prod_{p \equiv a \pmod q} \, \left(1+\frac{z-1}{p^s} \right) \, e^{(1-z)/p^s}\right). 
\end{multline}
We now bound $|\Gzs|$. For $\sigma \ge 3/4$, we have $\sum_{p, r \ge 2} \, 1/|rp^{rs}|$ $\, \le \, \sum_p \, \sum_{r \ge 2} \, 1/p^{3r/4} \, \le \, \sum_p \, 1/p^{3/2} \, \ll \, 1$. Thus the exponential factor in  
\eqref{eq:Gzs} has size $O_K(1)$   
for all $|z| \le K+1$.  

Moreover, for $\sigma \ge 3/4$ and $|z| \le K+1$, 
the last infinite product in \eqref{eq:Gzs} is 
\begin{equation}\label{eq:ModifEulerProdBdd}
\prod_{p \equiv a \pmod q} \, \left(1+\frac{z-1}{p^s} \right) \, \left(1+\frac{1-z}{p^s} + O\left(\frac{|z-1|^2}{p^{2\sigma}}\right) \right) \, = \, \prod_{p} \, \left(1+O_K\left(\frac1{p^{3/2}}\right) \right) \, \ll_K \, 1,
\end{equation}
Lastly, for all $s$ with $\sigma \ge 1-c_0/\log q$, we have $\ell^\sigma \, \ge \, \ell^{1-c_0/\log q} \, \ge \, e^{-c_0} \cdot \ell \, \ge \, 2\ell/3$, so that
$$\sum_{\ell \mid q} \, \left|\log\left(1-\frac1{\ell^s}\right)\right| \, \le \, \sum_{\ell \mid q} \, \frac1{\ell^\sigma} \, + \, O\left(\sum_{\ell \ge 2} \, \frac1{\ell^{3/2}} \right) \, \le \, \frac32\sum_{\ell \mid q} \, \frac1{\ell} \, + \, O(1) \, \le \, \frac32\sum_{\ell \le \omega(q)} \, \frac1{\ell} \, + \, O(1),$$
and by Mertens' Theorem, this last sum is $\le \, (3/2) \log_3 q \, + \, O(1)$. Inserting all these observations into \eqref{eq:Gzs}, we get $\Gzs \, \ll \, (\log_2 q)^{3/2}$, uniformly for all $s$ with $\sigma \ge 1-c_0/\log q$. 

Hence by \eqref{eq:omegaanLSDPrep}, we find that the sequence $\{\zPoweromegaan\}_n$ has property $\mathcal P(1, \{\alpha_\chi\}_\chi; c_0, \Omega)$, where $\alpha_\chi \, = \, \bbm_{\chi = \chi_0} \, + \, (z-1)\chibara/\phi(q)$ and $\Omega(t) \ll (\log_2 q)^{3/2}$. An entirely analogous argument shows that the sequence $\{|z|^{\omegaan}\}_n$ has property $\mathcal P(1, \{\beta_\chi\}_\chi; c_0, \Omega)$, with $\beta_\chi \, = \, \bbm_{\chi = \chi_0} \, + \, (|z|-1)\chibara/\phi(q)$. This neatly places us in the setting of Theorem \ref{thm:LFuncLSDVariant2}, with $\nu = 1$, $\CalM \ll (\log_2 q)^{3/2}$, $\delta_0=1$, and with  
$\max\{|\alphachiZ|, |\alphachie|, |\beta_{\chi_0}|, |\beta_{\chi_e}|\} \, \le \, 1+(K+1)/\phi(q) \le 2$. (Note that this is one of the situations where the hypotheses of Theorem \ref{thm:LFuncLSDVariant2} are easier to verify than those of Theorem \ref{thm:LFuncLSDVariant1}.) Since 
$\sum_\chi \, \alpha_\chi \cdot \chi(b) \, = \, 1 \, + \, {(z-1)}\sum_\chi \, \chibar(b) \chi(b)/{\phi(q)} = z, \text{ and likewise } \sum_\chi \, \beta_\chi \cdot \chi(b) = |z|, \text{ for any }b\text{ coprime to }q,$
it follows that $\Lambda_q = |z| \le K$. Theorem \ref{thm:LFuncLSDVariant2}(1) yields, uniformly in $q \le (\log x)^K$ and $|z| \le K$, 
\begin{equation}\label{eq:zomegaanLSD}
\sum_{n \le x} \, z^{\omegaan} \, - \, \frac{x}{(\log x)^{(1-z)/\phi(q)}} \, \sum_{j=0}^N \, \frac{C_j(z)}{(\log x)^j} \, \ll \, \frac{N!\,(71 \, c_1^{-1})^N \cdot x}{(\log x)^{(N+1)(1-\epsilon_0) - \frac{K+1}{\phi(q)}}} \, + \, x\exp\left(-\sqrt{\frac{c_0\log x}{16K}}\right).
\end{equation}
Here we have used $1-\eta_e \ge c_1\,q^{-\epsilon_0/(K+1)}$ to see that $(1-\eta_e)^{-|\Ree(\alphachie)|} \ll q^{\epsilon_0/\phi(q)} \ll 1$. Moreover, $C_j(z)$ is a holomorphic function (of $z$) on the disk $|z| \le K$, defined by  
\begin{equation}\label{eq:CjzDef}
\boldsymbol{C_j(z) = \frac{1/j!}{\Gamma\left(1-j+\frac{z-1}{\phi(q)}\right)} \cdot \frac{\mathrm d^j}{\mathrm{d}s^j}\Bigg\vert_{s=1}\frac{\Gzs}s \cdot \left(\prod_{\chi} L(s, \chi)^{\bbm_{\chi = \chi_0} \, + \, \frac{(z-1)\chibara}{\phi(q)}}\right) \cdot (s-1)^{1+\frac{z-1}{\phi(q)}}}.   
\end{equation} 

Now since $\#\{n \le x: \omegaan=k\}$ is the coefficient of $z^k$ in the polynomial $\sum_{n \le x} \, z^{\omegaan}$, estimate \eqref{eq:omegaa(n)Full} follows from \eqref{eq:zomegaanLSD} by applying Cauchy's integral formula on the circle of radius $K$ centered at the origin. 
The polynomials $P_{j, k}(T)$ appearing in \eqref{eq:omegaa(n)Full} are given by the following identities. 
\begin{equation}\label{eq:omegaanPjk}
P_{j, k}(T) \coloneqq \frac1{2\pi i} \int_{|z| = K} \, \frac{C_j(z) e^{zT/\phi(q)}}{z^{k+1}} \, \mathrm dz \, = \, \frac1{k!} \cdot \frac{\mathrm d^k}{dz^k}\Bigg\vert_{z=0} C_j(z) e^{zT/\phi(q)} \, = \, \sum_{r=0}^{k} \, \frac{C_j^{(r)}(0)}{r!} \cdot \frac{(T/\phi(q))^{k-r}}{(k-r)!},   
\end{equation} 
The last equality above uses the generalized product formula, with $C_j^{(r)}$ being the $r$-th derivative of $C_j$. 
To show \eqref{eq:P0k}, we will need 
the following lemma, which also seems to be of general interest. 
\begin{lem}\label{lem:PrimeRecipDeviationLimit}
Uniformly in all moduli $q \ge 1$ and in coprime residues $a$ mod $q$, we have
$$\displaystyle{\lim_{s \rightarrow 1} \, \sum_{p \equiv a \pmod q} \, \frac1{p^{s}} \, - \, \frac1{\phi(q)} \log\left(\frac1{s-1}\right) \, = \, \frac1{\pqa}  \, + \, O\left(\frac{\log(2q)}{\phi(q)}\right)}.$$
\end{lem}
\begin{proof}[Proof of Lemma \ref{lem:PrimeRecipDeviationLimit}]
It is well known that the limit exists, hence it suffices to show the lemma with $s \rightarrow 1+$ along the real numbers. Write $s = 1+1/\log X$ with $X \rightarrow \infty$. For all $p \le X$, we can write $p^{-1-1/\log X} \, = \, p^{-1} \exp(-\log p/\log X) \, =\, 1/p \, + \,O(\log p/p\log X)$. Now by \cite[Lemma (6.3), p.\,699]{norton76} or \cite[Remark 1]{pomerance77}, we have $\sum_{p \le X:~ p \equiv a \pmod q} \, 1/p \, = \, \log_2 X/\phi(q) + 1/\pqa + O(\log(2q)/\phi(q))$ for all $X \ge 3q$. Moreover by partial summation, the last estimate also yields $\sum_{p \le X:~ p \equiv a \pmod q} \, 1/p \, = \, \log(\pqa) /\pqa + O(\log(2q) \cdot \log X/\phi(q))$ for all such $X$. Combining these three observations, we get 
\begin{equation}\label{eq:PrimeRecipAPSmallp}
\sum_{\substack{p \le X\\ p \equiv a \pmod q}} \, \frac1{p^{1+1/\log X}} \, = \, \frac{\log_2 X}{\phi(q)} \, + \, \frac1{\pqa}  \, + \, O\left(\frac{\log(2q)}{\phi(q)} \, + \, \frac{\log(\pqa) /\pqa} {\log X}\right),   
\end{equation}
uniformly in all $X \ge 2q$. Next, by partial summation and the Brun--Titchmarsh theorem, 
\begin{equation*}
 \sum_{\substack{p > X\\ p \equiv a \pmod q}} \, \frac1{p^{1+1/\log X}} \, \ll \, \int_X^\infty \, \Bigg(\sum_{\substack{X < p \le t\\p \equiv a \pmod q}} \, 1\Bigg) \, \frac{\mathrm dt}{t^{2+1/\log X}} \, \ll \, \frac1{\phi(q)\log X} \int_X^\infty \,  \frac{\mathrm dt}{t^{1+1/\log X}} \, \ll  \, \frac1{\phi(q)},
\end{equation*}
uniformly in all $X \ge q^2$. Letting $X \rightarrow \infty$ in this estimate and in \eqref{eq:PrimeRecipAPSmallp}, we deduce that 
$$\displaystyle{\lim_{s \rightarrow 1+} \, \sum_{p \equiv a \pmod q} \, \frac1{p^{s}} \, - \, \frac1{\phi(q)} \log\left(\frac1{s-1}\right) \, = \, \lim_{X \rightarrow \infty} \, \sum_{p \equiv a \pmod q} \, \frac1{p^{1+1/\log X}} \, - \, \frac{\log_2 X}{\phi(q)}}$$
is equal to $1/\pqa + O(\log(2q)/\phi(q))$. This completes the proof of Lemma \ref{lem:PrimeRecipDeviationLimit}. 
\end{proof}
We now proceed to estimate $C_0(z)$ uniformly in all complex $z$ with $|z| \le K+1$.  
Since $L(s, \chi_0) = \zeta(s) \prod_{\ell \mid q} (1-1/\ell^s)$, we see from \eqref{eq:Gzs} and \eqref{eq:pPowersInAP} that for all $s$ with $\sigma>1$, we have 
\begin{multline*}\allowdisplaybreaks
 \Gzs \cdot \left(\prod_{\chi} L(s, \chi)^{\bbm_{\chi = \chi_0} \, + \, {(z-1)\chibara}/{\phi(q)}}\right) \cdot (s-1)^{1+(z-1)/\phi(q)} \\
 = \zeta(s)(s-1) \cdot \exp\left((z-1)\left\{\sum_{p \equiv a \pmod q} \, \frac1{p^{s}} \, + \, \frac{\log(s-1)}{\phi(q)}\right\}\right) \cdot \left(\prod_{p \equiv a \pmod q} \, \left(1+\frac{z-1}{p^s} \right) \, e^{(1-z)/p^s}\right).
\end{multline*}
Letting $s \rightarrow 1+$ above and using $\lim_{s \rightarrow 1+} \, \zeta(s)(s-1) = 1$, we obtain from \eqref{eq:CjzDef} and Lemma  \ref{lem:PrimeRecipDeviationLimit},   
\begin{align}\allowdisplaybreaks\label{eq:C0zIntermediate}
C_0(z) \, = \,  e^{(z-1)/\pqa} \, \left(1 +  O\left(\frac{\log q}{\phi(q)}\right)\right) \cdot \left(\prod_{p \equiv a \pmod q} \, \left(1+\frac{z-1}{p} \right) \, e^{(1-z)/p}\right),  
\end{align}
uniformly in $|z| \le K+1$. Here we also used the fact that $1/\Gamma(1+(z-1)/\phi(q)) \, = \, 1+O_K(1/\phi(q))$, which holds true as $(K+2)/\phi(q) \le 5(K+2)/2q^{1/2} \le 1/2$. 

Now since $|1-z|/p \le (K+2)/q \le 1/(K+2)$ for all $p>q$, we have  
$$\sum_{\substack{p>q\\p \equiv a \pmod q}} \, \left\{\log \left(1+\frac{z-1}{p} \right) + \frac{1-z}p \right\} \, \ll_K \, \sum_{p>q} \frac1{p^2} \, \ll \, \frac1{q\log q},$$
which shows that the contribution of all primes $p>q$ to the last infinite product in \eqref{eq:C0zIntermediate} is $1+O(1/q\log q)$. Hence from \eqref{eq:C0zIntermediate}, we deduce that uniformly in $|z| \le K+1$, we have
\begin{equation}\label{eq:C0z}
\displaystyle{C_0(z) \, = \, e^{(z-1)/\pqa} \, \left(1 +  O\left(\frac{\log q}{\phi(q)}\right)\right) \cdot \left(1+\frac{z-1}{\pqa} \right) \, e^{(1-z)/\pqa} \, = \,  \left(1 - \frac1{\pqa}\right) + \frac{z}{\pqa} +  O\left(\frac{\log q}{\phi(q)}\right)}.   
\end{equation}
(Note that if $\pqa > q$, then $1/\pqa$ gets absorbed in the error term in all these computations.)
As a consequence, we also obtain uniformly in all positive integers $r$, 
\begin{equation}\label{eq:C0r0}
\displaystyle{C_0^{(r)}(0) \, = \, \frac{r!}{2\pi i} \, \int_{|z|=K} \, \frac{C_0(z)}{z^{r+1}} \, \mathrm dz \,  = \,  \frac{\bbm_{r=1}}{\pqa} \, + \, O_K\left(\frac{r!}{K^r} \cdot \frac{\log q}{\phi(q)}\right)}.   
\end{equation}
Inserting \eqref{eq:C0z} and \eqref{eq:C0r0} into the last equality in \eqref{eq:omegaanPjk}, we obtain with $Y \coloneqq \log_2 x$, 
\begin{equation}\label{eq:P0kNearFinal}
P_{0, k}(Y) \, = \, \left(1 - \frac1{\pqa} \right) \frac{(Y/\phi(q))^k}{k!} + \frac{(Y/\phi(q))^{k-1}}{\pqa \, (k-1)!} \, + \, O\left(\frac{\log q}{K^k \, \phi(q)} \, \sum_{r=0}^k \, \frac{(KY/\phi(q))^{k-r}}{(k-r)!}\right).
\end{equation}
Since the sum in the $O$--term above is at most $e^{KY/\phi(q)}$, we obtain the second bound of \eqref{eq:P0k}. 
We will now show that the (stronger) first bound holds uniformly in the range $k \le KY/\phi(q)$. 

It is worth noting that in the more restricted range $k \le (K-\epsilon_0)Y/\phi(q)$, 
estimate \eqref{eq:P0kNearFinal} already yields the first bound in \eqref{eq:P0k}: This is because for all such $k$, the sum in the $O$--term of \eqref{eq:P0kNearFinal} is at most $(KY/\phi(q))^k \cdot (k!)^{-1} \, \sum_{r=0}^k \, (k\phi(q)/KY)^r \, \le \, (KY/\phi(q))^k \cdot (k!)^{-1} \, \sum_{r \ge 0} \, (1-\epsilon_0/K)^r$. However in the \textit{full range} $k \le KY/\phi(q)$, the first bound in \eqref{eq:P0k} follows immediately from \eqref{eq:C0z} and   
 
\begin{prop}\label{prop:omegaaSPM} With $Y \coloneqq \log_2 x$, we have uniformly in all $k \le KY/\phi(q)$, 
\begin{equation}\label{eq:P0kSPM}
P_{0, k}(Y) \, = \, \frac{(Y/\phi(q))^k}{k!} \, \left\{C_0\left(\frac{k\phi(q)}Y\right) \, + \, O\left(\frac{k \phi(q) \log q}{Y^2}\right)\right\}.
\end{equation}    
\end{prop}
\begin{proof}[Proof of Proposition \ref{prop:omegaaSPM}]
We use a variant of the saddle point method, adapting the argument in \cite[Theorem II.6.3]{tenenbaum15}. We start by using the 
first equality in \eqref{eq:omegaanPjk} to write    
\begin{equation}\label{eq:omegaaSPMStep1}
P_{0, k}(Y) \, = \, \frac{C_0(\tau)}{2\pi i} \, \int_{|z|=\tau} \, \frac{e^{zY/\phi(q)}}{z^{k+1}} \, \mathrm dz \, + \, \frac1{2\pi i} \, \int_{|z|=\tau} \, \left\{C_0(z)-C_0(\tau)-(z-\tau)C_0'(\tau)\right\} \cdot \frac{e^{zY/\phi(q)}}{z^{k+1}} \, \mathrm dz,    
\end{equation}
where we have chosen $\tau \coloneqq k \phi(q)/Y \le K$, so as to ensure that 
$$\int_{|z|=\tau} \, (z-\tau)\frac{e^{zY/\phi(q)}}{z^{k+1}} \, \mathrm dz \, = \, \int_{|z|=\tau} \, \frac{e^{zY/\phi(q)}}{z^{k}} \, - \, \tau \, \int_{|z|=\tau} \, \frac{e^{zY/\phi(q)}}{z^{k+1}} \, \mathrm dz \, = \, 0$$
via Cauchy's integral formula. As such, it follows from \eqref{eq:omegaaSPMStep1} that 
\begin{equation}\label{eq:omegaaSPMStep2}
P_{0, k}(Y) = C_0\left(\tau\right) \, \frac{(Y/\phi(q)^k}{k!} \, + \, \frac1{2\pi i} \, \int_{|z|=\tau} \, (z-\tau)^2 \left(\int_0^1 \, (1-u)\, C_0''(uz + (1-u)\tau) \, \mathrm du\right)  \, \frac{e^{zY/\phi(q)}}{z^{k+1}} \mathrm dz;  
\end{equation}
here the identity $C_0(z)-C_0(\tau)-(z-\tau)C_0'(\tau) \, = \, (z-\tau)^2 \, \int_0^1 \, (1-u)\, C_0''(\tau + u(z-\tau)) \, \mathrm du$ can be verified by integrating by parts twice. To complete the proof of Proposition \ref{prop:omegaaSPM}, it only remains to show that the entire double integral above is absorbed in the error term of \eqref{eq:P0kSPM}.

Now we observe that $C_0''(w) \, \ll \, \log q/\phi(q)$ uniformly in all complex $w$ with $|w| \le K$: This follows by using \eqref{eq:C0z} in Cauchy's integral formula on the disk of radius $1$ centered at $w$. By this observation, we have 
$C_0''(uz + (1-u)\tau) \, \ll \, \log q/\phi(q)$ uniformly in all $u$ and $z$ appearing on the right of \eqref{eq:omegaaSPMStep2}. 
Writing $z = \tau e^{i\theta}$, we thus find that the entire double integral in \eqref{eq:omegaaSPMStep2} has size 
\begin{equation*}  
\ll \, \frac{\log q}{\tau^{k-2} \, \phi(q)} \, \int_0^{2\pi} \, |e^{i \theta} - 1|^2 \cdot e^{\tau Y\cos \theta/\phi(q)} \, \mathrm d\theta \, = \, \frac{\log q}{\phi(q)} \cdot \left(\frac Y{k\phi(q)}\right)^{k-2}
\, \int_0^{2\pi} \, (1-\cos\theta) \, e^{k\cos \theta} \, \mathrm d\theta
\end{equation*}
The integral above is at most $2 \pi \, + \, 2\int_0^{\pi/2} \, (1-\cos\theta) \, e^{k\cos \theta} \, \mathrm d\theta \, \le \, 2\pi \, + \, 2\int_0^1 \, e^{ku}\sqrt{1-u} \, \mathrm du$ $\le \, 2\pi \, + \, 2 e^k \, k^{-3/2} \, \int_0^k \, e^{-v}v^{1/2} \, \mathrm dv \ll e^k \, k^{-3/2}$. (In this chain of inequalities, we first took $u \coloneqq \cos \theta$, noted that $\sqrt{1+u} \ge 1$, and then let $v \coloneqq k(1-u)$, noting that $\int_0^\infty \, \, e^{-v}v^{1/2} \, \mathrm dv = \Gamma(3/2)$.) 
Collecting these observations, we find that the double integral in \eqref{eq:omegaaSPMStep2} is $\ll\, (\log q/\phi(q)) \cdot (Y/\phi(q))^{k-2} \cdot ke^k/k^{k+1/2}$. By Stirling's formula, this expression is absorbed in the error term of \eqref{eq:P0kSPM}, as desired.
\end{proof}
This establishes all the assertions of Theorem \ref{thm:omegaa(n)}, except the very last one. But the last assertion  
(corresponding to the case when the Siegel zero does not exist) can be dealt with by all the same arguments, but only by replacing the use of subpart \textbf{(1)} of Theorem \ref{thm:LFuncLSDVariant2} by subpart \textbf{(2)}. \hfill \qedsymbol  
\subsection{Proof of Theorem \ref{thm:Omegaa(n)}} 
 
We just mention the main changes from the proof of Theorem \ref{thm:Omegaa(n)}. This time we start by estimating $\sum_{n \le x} \, z^{\Omega_a(n)}$ for $\boldsymbol{|z| \le \min\{K+2, \, \pqa(1-\epsilon_0/100)\}}$. The analogue of \eqref{eq:omegaanLSDPrep}
is 
$\sum_n \, {\zPowerOmegaan}/{n^s} \, = \,  \left(\prod_{\chi} \, L(s, \chi)^{\bbm_{\chi = \chi_0} \, + \, (z-1)\chibara/\phi(q)}\right) \cdot \GzTils$, 
where $\GzTils$ is defined the same as $\Gzs$ in \eqref{eq:Gzs}, but only with the last infinite product changed to 
\begin{align}\allowdisplaybreaks\label{eq:OmegaaEulerProdModif}
\prod_{p \equiv a \pmod q} \, \left(1- \frac{z}{p^s} \right)^{-1} \left(1- \frac{1}{p^s} \right) \, e^{(1-z)/p^s}.  
\end{align}
Now since $K$ and $\epsilon_0$ are fixed throughout, we may choose our constant $c_0 \coloneqq c_0(K, \epsilon_0)$ (defined at the start of subsection \cref{subsec:GeneralMainResults}) to be small enough that $\boldsymbol{e^{c_0} < \min\{K, \, (1-\epsilon_0/100)^{-1}\}}$. For all primes $p>q$ and for all complex numbers $s$ with $\sigma = \Ree(s) \ge 1-c_0/\log q$, we have $|z/p^s| \le (K+1)/q^{1-c_0/\log q} \le e^{c_0}/(K+2) < 1$, where we have recalled that $q>(K+2)^2$. Moreover, if $\pqa \le q$, then for all such $s$, we also have $|z/(\pqa)^s| \le (|z|/\pqa) \cdot \exp(c_0 \log(\pqa)/\log q) \le e^{c_0}(1-\epsilon_0/100) < 1$. These two observations show that the product \eqref{eq:OmegaaEulerProdModif} defines an analytic function of size $O(1)$ on the half plane $\sigma \ge 1-c_0/\log q$. (Here we wrote $(1-z/p^s)^{-1} = 1+z/p^s + O(1/p^{2\sigma})$.)

We can now proceed as we did for Theorem \ref{thm:omegaa(n)} to see that 
once again Theorem \ref{thm:LFuncLSDVariant2} applies with the \textbf{same} $\alpha_\chi, \beta_\chi, \nu, \Omega, \CalM$ and $\delta_0$ as before. Hence, the analogue of \eqref{eq:zomegaanLSD} for $\sum_{n \le x} \, \zPowerOmegaan$ holds true, uniformly in $q \le (\log x)^K$ and in $|z| \le \min\{K+2, \, \pqa(1-\epsilon_0/100)\}$: The only difference is that  
$C_j(z)$ 
is replaced by the function $\Cjtil(z)$ 
which is defined exactly as in \eqref{eq:CjzDef}, with the obvious replacement of $\Gzs$ by $\GzTils$. This proves the first assertion \eqref{eq:Omegaa(n)Full}, with 
\begin{equation}\label{eq:OmegaanQjk}
Q_{j, k}(T) \coloneqq \frac1{2\pi i} \int_{|z| = R} \, \frac{\Cjtil(z) e^{zT/\phi(q)}}{z^{k+1}} \, \mathrm dz \, = \, \frac1{k!} \cdot \frac{\mathrm d^k}{dz^k}\Bigg\vert_{z=0} \Cjtil(z) e^{zT/\phi(q)} \, = \, \sum_{r=0}^{k} \, \frac{\Cjtil^{(r)}(0)}{r!} \cdot \frac{(T/\phi(q))^{k-r}}{(k-r)!},
\end{equation}
and with $\boldsymbol{R \coloneqq \min\{K, \, \pqa(1-\epsilon_0)\}}$ as defined in the statement of Theorem \ref{thm:Omegaa(n)}. 

Now proceeding exactly as we did for \eqref{eq:C0z} and \eqref{eq:C0r0}, but replacing all instances of ``$|z| \le K+1$" by ``$|z| \le \min\{K+2, \, \pqa(1-\epsilon_0/100)\}$",  we see that the analogues of these two estimates are 
\begin{equation}\label{eq:C0Til}
    \CZeroTil(z) \, =\, \left(1-\frac z{\pqa}\right)^{-1} \left(1-\frac1{\pqa}\right) \, + \, O\left(\frac{\log q}{\phi(q)}\right),\\  
\end{equation}\\
$\text{ uniformly in all complex }z\text{ with }|z| \le \min\left\{K+2, \, \pqa\left(1-{\epsilon_0}/{100}\right)\right\},$ and 
\begin{align}\allowdisplaybreaks
\CZeroTil^{(r)}(0) \, = \, \frac{r!}{(\pqa)^r} \left(1-\frac1{\pqa}\right) \, + \, O\left(\frac{r!}{R^r} \cdot \frac{\log q}{\phi(q)}\right), \text{ uniformly in }r \in \NatNos \cup \{0\}.\label{eq:CrTil}
\end{align}
Inserting these two estimates into the last expression for $Q_{0, k}$ in \eqref{eq:OmegaanQjk}, we find that
\begin{align}\allowdisplaybreaks
Q_{0, k}(Y) \, &= \, \frac1{(\pqa)^k} \cdot \left(1-\frac1{\pqa}\right) \, \sum_{r=0}^k \, \frac{(\pqa Y/\phi(q))^{k-r}}{(k-r)!} \, + \, O\left(\frac{\log q}{R^k \, \phi(q)} \, \sum_{r=0}^k \, \frac{(R Y/\phi(q))^{k-r}}{(k-r)!}\right)\nonumber\\ 
& = \, \frac1{(\pqa)^k} \cdot \left(1-\frac1{\pqa}\right) \, \sum_{r=0}^k \, \frac{(\pqa Y/\phi(q))^r}{r!} \, + \, O\left(\frac{\log q}{\phi(q)} \cdot \frac{e^{RY/\phi(q)}}{R^k}\right) \label{eq:Q0kIntermed} 
\end{align}
with $Y \coloneqq \log_2 x$. Now for $k \le (1-\epsilon_0)\pqa Y/\phi(q)$, we see that 
$$\sum_{r=0}^{k-1} \, \frac{(\pqa Y/\phi(q))^r}{r!} \, \le \, \frac{(\pqa Y/\phi(q))^{k-1}}{(k-1)!} \, \sum_{r=0}^{k-1} \,  \left(\frac{k\phi(q)}{\pqa Y}\right)^{k-1-r} \, \le \, \frac{(\pqa Y/\phi(q))^{k-1}}{(k-1)!} \, \sum_{m \ge 0} \, (1-\epsilon_0)^m,$$
which is $\ll \, (\pqa Y/\phi(q))^{k-1}\Big/(k-1)!$. Inserting this into \eqref{eq:Q0kIntermed}, we  \eqref{eq:Q0kSmallk}. On the other hand, for $k \ge (1+\epsilon_0)\pqa Y/\phi(q)$, the difference between $e^{\pqa Y/\phi(q)}$ and the last sum in \eqref{eq:Q0kIntermed} is   
$$\sum_{r \ge k+1} \, \frac{(\pqa Y/\phi(q))^r}{r!} \, \le \, \frac{(\pqa Y/\phi(q))^{k+1}}{(k+1)!} \cdot \sum_{r \ge k+1} \, \left(\frac{\pqa Y}{k\phi(q)}\right)^{r-k-1} \, \ll_{\epsilon_0} ~ \frac{(\pqa Y/\phi(q))^{k+1}}{(k+1)!},$$
which also establishes estimate \eqref{eq:Q0kLargek} in Theorem \ref{thm:Omegaa(n)}. 
Finally \eqref{eq:Q0kSmallkSPM} follows from \eqref{eq:C0Til} and the following analogue of Proposition \ref{prop:omegaaSPM}. 
\begin{prop}
 With $Y = \log_2 x$ and $R = \min\{K, \, \pqa(1-\epsilon_0)\}$ as above, we have uniformly in all $k \le RY/\phi(q)$, 
\begin{equation}\label{eq:Q0kSPM}
\displaystyle{Q_{0, k}(Y) \, = \, \frac{(Y/\phi(q))^k}{k!} \, \left\{ \CZeroTil\left(\frac{k\phi(q)}Y\right) \, + \, O\left(\frac{k}{(\pqa Y/\phi(q))^2} \, + \, \frac{k \phi(q) \log q}{Y^2}\right) \right\}.}
\end{equation} 
\end{prop}
The proof of this result is entirely analogous to that of Proposition \ref{prop:omegaaSPM}: The only main difference 
is that this time we use the bound $\CZeroTil''(w) \, \ll \, 1/(\pqa)^2 \, + \, \log q/\phi(q)$,  which holds uniformly in $|w| \le R$, and is obtained by an application of \eqref{eq:C0Til} in conjunction with Cauchy's integral formula on the disk of radius $\epsilon_0$ centered at $w$. This concludes the proof of Theorem \ref{thm:Omegaa(n)}. \hfill \qedsymbol  
\section*{Acknowledgements}
I thank Paul Pollack for useful discussions. I am sincerely  grateful to the Department of Mathematics at the University of Georgia for their support and hospitality. The application extending the Sathe--Selberg Theorem  
to study the local laws of $\omega_a(n)$ and $\Omega_a(n)$ arose out of a discussion with Larry Guth (MIT), hence I would like to thank him for visiting UGA to give the Cantrell Lectures in March 2025, and for expressing interest in my work.   
\bibliographystyle{amsplain}
\bibliography{LfuncLSD}
\bigskip 
\end{document}